\documentclass{amsart}      % This tells LaTeX what kind of 

\usepackage{amsmath, amssymb, amsthm, amscd}   % AMS maths styles

\usepackage[dvips]{graphicx}
\usepackage{enumerate}
\usepackage{amscd}
\newtheorem{lemma}{Lemma}[section]
\newtheorem{theorem}[lemma]{Theorem}
\newtheorem{corollary}[lemma]{Corollary}
\newtheorem{proposition}[lemma]{Proposition}
\theoremstyle{definition}
\newtheorem{definition}{Definition}[section]

\newtheorem{claim}[lemma]{Claim}

\newcommand{\blackboard}[1]{\ensuremath{\mathbb{#1}}}

\newcommand{\complexes}{\blackboard{C}}
\newcommand{\hyperbolic}{\blackboard{H}}

\newcommand{\integers}{\blackboard{Z}} %
\newcommand{\reals}{\blackboard{R}}
\newcommand{\naturals}{\blackboard{N}}

\newcommand{\PSL}{\ensuremath{\mathrm{PSL}}}

\newcommand{\Int}{\mathrm{Int}}
\newcommand{\pt}{\{\mathrm{pt.}\}}
\newcommand{\Fr}{\mathrm{Fr}}
\newcommand{\length}{\mathrm{length}}
\newcommand{\PML}{\mathbb P\mathcal{ML}}
\newcommand{\PM}{\mathbb P \mathcal{M}}
\newcommand{\Area}{\mathrm{Area}}

\newcommand{\name}[1]{\noindent \small #1}
\newcommand{\ie}{i.e.\ }

\newcommand{\vol}[1]{\textbf{#1}}

{%
\settowidth\labelwidth{\@biblabel{#1}}
\addtolength{\leftskip}{\labelwidth}
\addtolength{\leftskip}{\labelsep} 
\addtolength{\leftskip}{-\labelwidth} 
\addtolength{\leftskip}{\labelwidth}
}{}

\newcommand{\classmark}[1]{\noindent \small \textbf{AMS Classification:} #1}
\renewcommand{\labelenumi}{{\rm(\theenumi)}}

\begin{document}

\title[Realising end invariants]{Realising end invariants by limits of minimally parabolic, geometrically
finite groups}
\author{Ken'ichi Ohshika}

 \date{}

\maketitle
\begin{abstract}
We shall show that for a given homeomorphism type and a set of end 
invariants (including the parabolic locus) with necessary topological conditions
which a topologically 
tame Kleinian  group with that homeomorphism type must satisfy, there is 
an algebraic limit of minimally parabolic, geometrically finite Kleinian groups which has exactly that 
homeomorphism type and  end invariants.
This shows that the Bers-Sullivan-Thurston density conjecture follows from Marden's 
conjecture proved by Agol, Calegari-Gabai combined with Thurston's uniformisation
theorem and the ending lamination conjecture proved by Minsky, partially collaborating
with Masur, Brock and Canary.

\end{abstract}

\section{Introduction}
A Kleinian group is said to be geometrically finite when the corresponding hyperbolic 3-orbifold has a convex core which has finite volume.
It was conjectured by Bers in \cite{Bers} that every b-group, \ie a Kleinian  group with a unique simply-connected invariant component in the domain of discontinuity, is an algebraic limit of quasi-Fuchsian group.
The conjecture was generalised by Sullivan and Thurston to the one saying that every finitely generated Kleinian group is an algebraic limit of quasi-conformal deformations of a geometrically finite group.
This generalised conjecture is called the Bers-Sullivan-Thurston density conjecture today.
The original  conjecture of Bers was solved by Bromberg \cite{Brom} using deformations of cone manifolds and Minsky's solution of the ending lamination conjecture for hyperbolic 3-manifolds with bounded geometry in \cite{MiT} and \cite{MiJ} which was prior to the solution of the general case.
(A hyperbolic manifold is said to have bounded geometry if there is a positive lower bound for the lengths of the closed geodesics.)
This result is then generalised to freely indecomposable Kleinian groups in Brock-Bromberg \cite{BB}.
The purpose of the present paper is to prove this density conjecture in general, including the case of freely decomposable Kleinian groups.
Before explaining what is the difficulty in the freely decomposable case, we shall review the ending lamination conjecture and its resolution, on which our proof of the main result depends.

The ending lamination conjecture, which is due to Thurston and has been one of the most famous conjectures in the modern Kleinian group theory, says that every topologically tame hyperbolic $3$-manifold is determined up to isometries by its homeomorphism type and   end invariants consisting of the conformal structures at infinity and the ending laminations.
Here a hyperbolic 3-manifold is said to be topologically tame when it is homeomorphic to the interior of a compact 3-manifold.
We shall explain what the end invariants are more in details.
Let $M=\hyperbolic^3/\Gamma$ be a topologically tame hyperbolic 3-manifold, and $M_0$ its non-cuspidal part.
We choose a compact core $C$ such that for each component $T$ of $\partial M_0$, the intersection $C \cap T$ is a core of $T$, and each component of $M \setminus C$ is a product open-interval bundle over a component of $\Fr C$.
We shall call such a compact core nice.
For a nice compact core $C$, we see that $(C, C \cap \partial M_0)$ has a structure of a pared manifold.
Each component of $\partial C \setminus  P$ faces an end of $M_0$.
The ends of $M_0$ have invariants as follows.
If an end of $M_0$ is geometrically finite, \ie has a neighbourhood  containing no closed geodesics, then there is a component of the quotient $\Omega_\Gamma/\Gamma$ of the region of the  discontinuity which can be regarded as the points at infinity of the end, and its marked conformal structure constitutes an invariant for the end.
When an end is geometrically infinite, it has an invariant called the
ending lamination, which is represented by a measured lamination contained  in
the Masur domain of the frontier component of $C$ facing the end, determined 
uniquely up to transverse measures if we fix a marking on the frontier component.
(The marking of $\Fr C$ is not uniquely determined by the marking of $M$ when $\Fr C$ has a compressible component.)
The ending lamination conjecture says that these pieces of information, the parabolic locus, the marked conformal structure
for the  geometrically finite
ends and the ending laminations for the geometrically infinite ends both up to the action of the auto-diffeomorphisms of the frontier component homotopic to the
identity in $C$, together with the homeomorphism type,  are
sufficient to determine the isometry type of the hyperbolic $3$-manifold.

The first essential contribution to solving the ending lamination conjecture was due to Minsky \cite{MiT}, \cite{MiJ} for hyperbolic 3-manifolds with freely indecomposable fundamental groups and bounded geometry.
After some work on the unbounded geometry case in special situations as in Minsky \cite{MiA} and \cite{MiG}, the conjecture for the general case when hyperbolic manifolds may have unbounded geometry was finally solved in Minsky \cite{Mi} and Brock-Canary-Minsky \cite{BCM} based on the work on the geometry of  the curve complex developed in Masur-Minsky \cite{MM1}, \cite{MM2}.
Although Brock-Canary-Minsky \cite{BCM} only contains a proof for freely indecomposable case, it can be generalised to freely decomposable case, as was written in the first section of \cite{BCM}, using the technique of constructing negatively curved metric on branched covers due to Canary \cite{Ca}.
This method in the case of manifolds with bounded geometry corresponding to \cite{MiT} was explained in \cite{OhT} and its generalisation to the case of manifolds having the same end invariants as those with bounded geometry, corresponding to \cite{MiG}, was sketched in Ohshika-Miyachi \cite{OM}.

%On the other hand, it has been conjectured, first by Bers in the
% case of boundary groups, and then by Thurston in the general
%form, that every finitely generated Kleinian group is an algebraic
%limit of quasi-conformal deformations of a geometrically finite group.
%In the case of freely indecomposable Kleinian groups, this conjecture
%follows from the ending lamination conjecture combined with the
%result in Ohshika \cite{OhL} where it was shown how to realise end invariants
%by limits of geometrically finite groups.
%Also, this Bers-Thurston density conjecture for freely indecomposable Kleinian
%groups without parabolic elements was proved by Brock and Bromberg \cite{BB}
%only using the ending lamination conjecture for hyperbolic $3$-manifolds with bounded geometry, which was proved by Minsky prior to the solution in general, and making use of deformations of cone-manifolds.

For the Bers-Sullivan-Thurston density conjecture in the freely decomposable case,  it is
first necessary to show that every finitely generated Kleinian group is topologically tame so that ending laminations can be defined.
This  is exactly what  Marden's tameness conjecture says, which was recently solved by
Agol and Calegari-Gabai independently (\cite{Ag}, \cite{CG}).
Once Marden's conjecture is solved, to show that every (finitely generated and
torsion-free) Kleinian
group is a limit of quasi-conformal deformations of a geometrically finite group, using the resolution of the ending lamination conjecture by Brock-Canary-Minsky, what we need to do is  to prove that any possible system of a homomorphism type, a parabolic locus, and end invariants can be realised in a Kleinian group which is a limit of minimally parabolic geometrically finite groups. 
In the case of freely decomposable groups, this consists of two main steps.
The first step is to provide a
convergence theorem for freely decomposable group generalising results of Kleineidam-Souto \cite{KS}, Lecuire \cite{Lecuire} and Ohshika \cite{Oh1}, in the last of which we only dealt with Kleinian groups without parabolic elements.
The second step is to show that the limit group obtained by the convergence theorem has the desired properties.
This in particular necessitates to show that the support of  an arational measured lamination contained in the Masur domain, which we know to be unrealisable, is homotopic in $M$ to an ending lamination.
This latter step was easy for the freely indecomposable case by work of Thurston and Bonahon, but in our general case, the argument is rather complicated.

We note that Namazi and Souto also have also given a proof of this latter step in  \cite{NS}, and have proved the Bers-Sullivan-Thurston density conjecture independently of our work.

We now state in the form of a theorem an affirmative solution of the Bers-Sullivan-Thurston density conjecture for general (topologically tame)  Kleinian  groups as explained above.

\begin{theorem}
Let $\Gamma$ be a finitely generated, torsion-free Kleinian group.
Then there is a geometrically finite Kleinian group $G$ without
infinite cyclic maximal parabolic subgroups such that $\Gamma$ 
is an algebraic limit of quasi-conformal deformations of $G$.
\end{theorem}

Kleinian groups without infinite cyclic parabolic groups are sometimes called minimally parabolic, as in our title.
It should be noted that minimally parabolic geometrically finite groups are convex cocompact unless they have  rank-$2$ abelian subgroups.

This theorem is derived as a corollary from the main theorem of the present paper, which we shall state below, combined with the  solution of the ending
lamination conjecture and Thurston's uniformisation theorem for compact atoroidal
$3$-manifolds with boundary (see for instance \cite{Mo}).
Before stating the theorem formally, we summarise what the theorem says.
We consider a geometrically finite Kleinian group $G$ and a nice core $C$ of $(\hyperbolic^3/G)_0$.
We give on non-torus components of $\partial C$ a union $P$ of disjoint non-parallel essential simple closed curves such that $(C,P)$ is a pared manifold, and on each non-torus component of $\partial C \setminus P$ either a conformal structure or a arational measured lamination contained in the Masur domain, with some conditions which are evidently necessary to make them a parabolic locus and end invariants.
Then what we shall claim in the theorem is that there is a topologically tame Kleinian group $\Gamma$ which is an algebraic limit of quasi-conformal deformations of $G$ such that $(\hyperbolic^3/\Gamma)_0$ has $P$ as the parabolic locus and the given conformal structures and laminations as end invariants under the natural marking.

\begin{theorem}
\label{main}
Let $G$ be a torsion-free geometrically finite Kleinian group without infinite cyclic maximal parabolic subgroups.
Let $C$ be a nice compact core of $(\hyperbolic^3/G)_0$.
Let $T$ denote the union of the torus components of $\partial C$.
Let $P$ be a disjoint union of annular neighbourhoods of essential simple closed curves on 
$\partial C$ such that $(C, P\cup T)$ is a pared manifold.
Let $\Sigma_1, \dots , \Sigma_m$ be the components of $\partial C
\setminus (P \cup T)$.
Among $\Sigma_1, \dots , \Sigma_m$, suppose that on each
$\Sigma_j$ with $j=1, \dots , n$ (possibly $n=0$), 
a marked conformal structure $m_j $ making the components of the frontier 
punctures is given, and that on each $\Sigma_j$ with $j=n+1,
\dots, m$ (possibly $m=n$),   an arational measured lamination $\mu_j$
contained in the
 Masur domain of
$\Sigma_j$ is given.
When $(C,P)$ is a trivial $I$-bundle (as a pared manifold) and $n=0$, we further
assume that the supports of  $\mu_1$ and $\mu_2$ are not  homotopic in
$C$.
When $(C,P)$ is a twisted $I$-bundle and $n=0$, we further assume that $\mu_1$ is
not a lift of a measured lamination in the base space of the $I$-bundle, which is
a non-orientable surface.

Then there is a sequence of quasi-conformal deformations $G_i$ of $G$ with isomorphisms  $\phi_i: G \rightarrow G_i$  converging algebraically to an isomorphism $\psi: G \rightarrow \Gamma$ to a topologically tame Kleinian group $\Gamma$  such that
$(\hyperbolic^3/\Gamma)_0$ has a nice compact core $C'$ with a
homeomorphism $\Phi : C \rightarrow C'$ inducing $\phi$ between
the fundamental groups, such that 
\begin{enumerate}
\item 
$\Phi(P \cup T)$ coincides with the parabolic locus of $C'$ regarded as a pared
manifold,
\item
the end of $(\hyperbolic^3/\Gamma)_0$ facing $\Phi(\Sigma_j)$ with $j=1, \dots , n$ is geometrically finite and has marked conformal structure at infinity equal to the one represented by $\Phi_*(m_j)$,
\item
and the end  of $(\hyperbolic^3/\Gamma)_0$ facing $\Phi(\Sigma_j)$ with $j=n+1, \dots , m$ is geometrically infinite with ending lamination represented by $\Phi(\mu_j)$.
\end{enumerate}
\end{theorem}

We should note that this theorem can be also regarded as a generalisation 
of the main theorem of Ohshika \cite{OhP}, where we 
generalised Maskit's theorem in \cite{Mas} on function groups to one on general
geometrically finite groups.

The proof of this theorem will proceed as follows.
We shall first construct quasi-conformal deformations $(G_i, \phi_i)$ of $G$
so that the conformal structure at infinity of $\hyperbolic^3/G_i$ restricted
to 
$\Sigma_j$  converges to $m_j$ for $j=1, \dots , n$ and   diverges towards
$[\mu_j]$ in the Thurston compactification of the Teichm\"{u}ller space for $j=n+1, \dots m$, and the length of
$\mu_j$ with respect to the hyperbolic metric compatible with the conformal structure $m_j$ is bounded.
The first of the two  main steps is to show that such a sequence of quasi-conformal structures converges algebraically after passing to a subsequence.
The proof of this fact relies on work of Otal \cite{OtT}, Kleineidam-Souto
\cite{KS} and a more recent result of Lecuire \cite{Lecuire} on the extension of the Masur domain.
%
%One may imagine that it is possible to prove the convergence part of the theorem
%above by simply combining the geometrically finite
%case dealt with in \cite{OhP} and a relative version of the argument we
%developed in
%\cite{Oh1}.
%Although this kind of argument may actually work, this is not the approach which
%we have chosen in this paper.
%That is because our argument presented here is simpler than proving a relative
%version of
%the convergence theorem, including that of  Kleineidam-Souto, which was essentially used in the
%argument in \cite{Oh1}.Also, our argument in the present paper has a merit that in the process of the
%proof, we get a
%concrete sequence of quasi-conformal deformations which
%is guaranteed to converge to any given  group that we want to realise as a
%limit.

The second main step is to show that the limit group has properties as we wanted.
It follows from results in Ohshika \cite{OhM}, Brock-Souto \cite{BS}  (or the general solution of Marden's conjecture) that the limit $\Gamma$ is topologically tame, and using some geometric argument, we can show that  the $\mu_j$ cannot be realised there.
A most difficult part of the second step is to show that (the images by a
homeomorphism from $\hyperbolic^3/G$ to $\hyperbolic^3/\Gamma$ of)  the
$\mu_j$ actually represent ending laminations.
 (Proposition \ref{ending lamination}.)
Section 6 will be entirely devoted to the proof of this fact.

\section{Preliminaries}
Throughout this paper, Kleinian groups are assumed to be finitely generated
and torsion free except for the case when we consider geometric limits, which may be infinitely generated.
Similarly, we always consider hyperbolic 3-manifolds with finitely generated fundamental groups except for geometric limits.
For a Kleinian group $G$, we consider the corresponding hyperbolic
$3$-manifold $M_G=\hyperbolic^3/G$.
Throughout this paper, we use the symbol $M$ with the name of a Kleinian group added as a subscript to denote the corresponding hyperbolic 3-manifold.
For a constant $\epsilon >0$ less than the three-dimensional Margulis
constant, we define the {\sl non-cuspidal part}, denoted by $({M_G})_0$,
to be
the complement
of the $\epsilon$-cusp neighbourhoods of $M_G$, that is,
neighbourhoods of cusps consisting of points where the injectivity radii are
less than $\epsilon/2$.

Since $G$ is assumed to be finitely generated, by Scott's core theorem
\cite{Sc},
there
is a compact $3$-manifold $C_G$, which we call a {\sl compact core}, embedded in
$M_G$
such that the inclusion from $C_G$ to $M_G$ is a homotopy
equivalence.
When $G$ has parabolic elements, it is more convenient to consider a relative
compact core of the non-cuspidal part $({M_G})_0$, whose existence
was proved by McCullough \cite{McC}.
A relative compact core intersects the boundary of the non-cuspidal part at tori
corresponding to $\integers \times \integers$-cusps and annuli  which are cores
of open annulus  corresponding to $\integers$-cusps, one
annulus for each $\integers$-cusp.
The intersection with the boundary of the non-cuspidal part  induces a structure
of pared manifold on a relative compact core, which we shall explain below.

A Kleinian group $G$ and the corresponding hyperbolic 3-manifold $M_G$ are said to be topologically tame when $M_G$ is homeomorphic to the interior of a compact $3$-manifold.
In this case, we can choose a relative compact core $C_G$ so that each component $E$ of $M_G \setminus C_G$ is homeomorphic to $F \times \reals$, where $F$ is the component of $\Fr C_G$ contained in the closure of $E$.
As was mentioned in Introduction, we call such a compact core {\sl nice}.

Following Thurston, we call a pair $(C,P)$ of a compact irreducible $3$-manifold
and a subsurface of its boundary  a {\sl pared manifold} when
\begin{enumerate}
\item
$P$ consists of disjoint incompressible tori and annuli,
\item
every incompressible (\ie $\pi_1$-injective) map from a torus to $C$ is homotopic into $P$,
\item 
and every incompressible map $(S^1 \times I, S^1 \times \partial I)
\rightarrow (C,P)$ is relatively homotopic to a map into $P$.
\end{enumerate}

The subsurface $P$ above is called the {\sl paring locus}.
When we consider a pared manifold which is a relative compact core of the non-cuspidal part of a hyperbolic $3$-manifold, we call its paring locus the {\sl parabolic locus}.

A {\sl compression body} $W$ is a connected $3$-manifold obtained from finitely
many product
$I$-bundles $S_1 \times I , \dots , S_m \times I$ over closed surfaces which
are not spheres
by attaching $1$-handles to $\cup S_k \times \{1\}$.
We assume that there is at least one $1$-handle;
hence we do not regard a product $I$-bundle as a compression body.
Exceptionally handlebodies are also regarded as compression bodies.
The union of the $S_k \times \{0\}$, called the interior boundary, is denoted by
$
\partial_i W$, and
the remaining boundary component coming from $S_k \times \{1\}$, which is called
the exterior boundary,
is denoted by $\partial_e W$.
We use the same symbols $\partial_e$ and $\partial_i$ to denote the unions of the exterior boundaries and the interior boundaries respectively for a disjoint union of compression bodies.
For a compact irreducible $3$-manifold $C$, there exists a submanifold $V$ which
is a
 disjoint union of compression bodies such that $V \cap 
\partial C= 
\partial_e
V$ is the union of compressible boundary components of $C$,
and each component of  $\partial_i V$ (unless it is empty) is either an incompressible surface in
$\mathrm{Int}C$ which is not parallel into $\partial C$ or an incompressible boundary component of $C$.
Such a manifold is unique up to isotopy
and is called the {\sl characteristic compression body} of $C$.
If $V$ is a characteristic compression body of $C$, then every incompressible (\ie $\pi_1$-injective)
map from a closed surface to $C$ is homotoped into $(C \setminus V) \cup
\partial_i V$.
The closure of $C\setminus V$ is a (possibly empty) boundary-irreducible manifold and none of
its components are product $I$-bundles unless $C$ itself is a product $I$-bundle.
These facts are proved by Bonahon \cite{BoA}, to which we refer the reader for further details.
%We need to use the following fact later in \S 5.

%\begin{lemma}
%\label{homotopy in C}
%Let $C$ be a compact irreducible $3$-manifold and $V$ its characteristic compression body.
%Let $s$ be an essential  simple closed curve on a component $S$ of $\partial_i V$.
%Suppose that $s$ is homotopic in $C$ to another simple closed curve $s'$ lying on $\partial V$ disjoint from $s$.
%Then they are homotopic in the manifold obtained by cutting $C$ along $S$.
%\end{lemma}
%\begin{proof}
%We first observe that there is no essential embedded annulus in $V$ both of whose boundary components lie on $\partial_i V$.
%This can be easily seen by using the fact that an essential annulus cannot intersect a compressing disc essentially at a simple closed curve, and hence an essential annulus whose boundary components lie on $\partial_i V$ cannot intersect a compressing disc essentially.

%Now, consider an annulus $A$ cobounded by $s$ and $s'$.
%This annulus must be essential since $s$ is essential on $\partial_i V$ hence also in $C$.
%We see that $\Int A$ can intersect $\partial_i
%\end{proof}

Let $S$ be a hyperbolic surface, possibly with geodesic boundary.
A geodesic lamination on $S$ is a closed set disjoint from $\partial S$
consisting of disjoint simple geodesics.
A measured lamination is a geodesic lamination with a holonomy-invariant
transverse measure.
We always assume that the support of the measure is the entire lamination.
For a measured lamination $\mu$, its support is denoted by $|\mu|$.
The space of measured laminations with the weak topology on transverse arcs is
denoted by $\mathcal{ML}(S)$.
The projective lamination space is $(\mathcal{ML}(S) \setminus
\{\emptyset\})/\reals_+$, where the action of $\reals_+$ is that of scalar
multiplication of the transverse measures, and is denoted by $\mathbb P\mathcal{ML}(S)$.
Thurston defined a compactification of the Teichm\"{u}ller space
$\mathcal{T}(S)$ whose
boundary is identified with $\PML(S)$.
When $S$ has boundary, we define its Teichm\"{u}ller space to be the space 
of marked hyperbolic structures with respect to which the lengths of the
boundary components are constant.
We sometimes say that $m_i\in \mathcal{T}(S)$ {\sl diverges towards} $\lambda \in
\PML(S)$ when $\{m_i\}$ converges to $\lambda$ in the Thurston compactification.

A geodesic lamination is said to be  {\sl arational} when every component of its
complement is
either simply connected, or an annulus around a cusp or a boundary component.
An arational measured lamination is always minimal, i.e., it does not have a proper sublamination.
We say that a measured lamination is {\sl maximal} when it is not a proper sublamination of another measured lamination.
Arational measured laminations are always maximal.

Consider a compact $3$-manifold $C$ and  an essential  subsurface $S$ contained in a component of $\partial C$.
We assume that no boundary components of $S$ are meridians (i.e.\ boundaries of compressing discs).
In the measured lamination space $\mathcal{ML}(S)$, we define the following subsets.
First, we set $\mathcal{WC}(S)$ to be the subset of $\mathcal{ML}(S)$ consisting
of disjoint weighted union of meridians lying on $S$.
Except for the case when $S$ has only one isotopy class of compressing discs,
we
define the {\sl Masur domain} of $S$ by $$
\mathcal{M}(S) =
\{\lambda \in \mathcal{ML}(S)\mid i(\lambda, \nu) >0 \  \text{for any } \nu \in
\overline{\mathcal{WC}}(S)\},$$
where $\overline{\mathcal{WC}}(S)$ denotes the closure of $\mathcal{WC}(S)$ in $\mathcal{ML}(S)$.
When $S$ has only one isotopy class of compressing discs, we define the Masur
domain by
\begin{align*}\mathcal{M}(S) =
\{\lambda \in \mathcal{ML}(S) \mid i(\lambda, \nu)>0, & \text{ for any } \nu  \text{
that is}
\\
&\text{ disjoint from a meridian } \}.
\end{align*}
We note that $\mathcal{M}(S)$ coincides with the entire measured lamination space if $S$ is incompressible.

We need to define another domain in $\mathcal{ML}(\partial V)=\mathcal{ML}(\partial_e V) \cup \mathcal{ML}(\partial_i V)$ larger than
the Masur domain when $V$ is a compression body,
following  Lecuire \cite{Lecuire}.
For a compression body $V$, we set
% Verfier la r?f?rence
\begin{align*}
\mathcal{D}(V)=\{\lambda \in \mathcal{ML}(\partial V)\mid
&\text{ there exists }
\eta >0 \text{ such that }
 i(\lambda, \partial D) > \eta\\&  \text{ for any compression disc
and any essential annulus } D\}.
\end{align*}

%We define $\mathcal D_e(V)$ to be $\mathcal D(V) \cap \mathcal{ML}(\partial_e V)$.
%How about $\mathcal{D}$?

The subspaces $\mathcal{WC}(S)$, $\mathcal{M}(S)$ and $\mathcal{D}(V)$ are all
invariant under scalar multiplication.
We put $\mathbb{P}$ to denote their images in the projective lamination space
$\PML(S)$ or $\PML(\partial V)$ .

The subgroup of the mapping class group of $S$ consisting of classes of
diffeomorphisms homotopic to the identity in $V$ is denoted by
$\mathrm{Mod}_0(V)$.
This group $\mathrm{Mod}_0(V)$ acts on $\PM(S)$ properly
discontinuously, and the limit set of $\mathrm{Mod}_0(V)$ in $\PML(S)$
is equal to $\overline{\mathbb P\mathcal{WC}}(S)$.
(Refer to Otal \cite{OtT}.)

We let  $\tilde{V}$ be the universal cover of $V$ and  $\tilde{S}$
the preimage of $S$ lying on  $\partial \tilde{V}$.
We fix a hyperbolic metric on $S$, which induces one on $\tilde{S}$.
Let $l$ be a leaf of a geodesic
lamination $\mu$ on $S$.
Consider a lift $\tilde{l}: \reals \rightarrow \tilde{S}$ of $l$.
We say that $l$ is {\sl homoclinic} if there are  sequences of points
$\{s_i\}, \{t_i\}$ on $\reals$ such that $|s_i - t_i| \rightarrow \infty$
whereas
$d_{\tilde{S}}(\tilde{l}(s_i), \tilde{l}(t_i))$ is bounded above.

We need to use the following lemma of Otal \cite{OtT}.
(This is contained in the proof of Proposition 2.10 of \cite{OtT}.
Kleineidam-Souto stated this as Lemma 4 in  \cite{KS}.)

\begin{lemma}
\label{no homoclinic}
The support of a measured lamination contained in $\mathcal{M}(S)$ cannot be
extended to a geodesic lamination with a
homoclinic leaf.
\end{lemma}
%\begin{proof}
%Let $\lambda$ be a measured lamination contained in $\mathcal{M}(S)$.
%First we note that there is a meridian $m$ on $S$ with respect to which
%$\lambda$ is in tight position, that is there is no arc on a leaf of $\lambda$
%which can be homotoped in $V$ into $m$ fixing the endpoints.

For a hyperbolic $3$-manifold $M$, there exists a unique minimal convex
submanifold that is a deformation retract.
Such a submanifold is called the {\sl convex core} of $M$.
When the convex core of $M_G$ is compact, $G$ is said to be {\sl
convex
cocompact}, and when the convex core has finite volume, $G$ is said to be
{\sl geometrically finite}.

More generally, for a hyperbolic 3-manifold $M$, an end of $M_0$ is said to be geometrically finite when it has
a neighbourhood intersecting no closed geodesics.
%Let $C$ be a relative compact core of $M_0$ and 
 %$e$  a geometrically infinite end of $M_0$ facing a component $\Sigma$ of the
%frontier of $C$.
By the resolution of Marden's conjecture, we can take a nice compact core $C$ of $M_0$ so that each end $e$ has a neighbourhood homeomorphic
to $\Sigma \times \reals$ for a component $\Sigma$ of $\Fr C$ facing $e$.
Suppose that $e$ is geometrically infinite.
%We can enlarge a compact core so that the component of $M_0 \setminus C$
%containing $e$ is $\Sigma \times (0,\infty)$, where $\Sigma \times \{0\}$ is
%identified with $\Sigma$,
%by the relative compact core theorem of McCullough.
%Citation
By Proposition 10.1  in Canary \cite{CaJ}, there is a sequence of closed geodesics
$c_k^*$ tending to $e$ such that $c_k^*$ is homotopic in $\Sigma \times [0,
\infty)$ to a simple closed curve $c_k$ such that $\{r_k c_k\}$ converges to a
measured lamination $\mu$ in $\mathcal{M}(\Sigma)$ for some $r_k \in (0,\infty)$.
In this situation, we say that $\mu$ {\sl represents the ending lamination} of
$e$.
Actually the transverse measure is irrelevant for the ending lamination.
The geodesic lamination which is the support of a measured lamination representing the ending lamination is called the ending lamination of $e$.

A {\sl pleated surface} is a map $f: S \rightarrow M$ from a hyperbolic surface to a
hyperbolic $3$-manifold  taking a puncture or a boundary component to a cusp
such that for any point $x \in S$, there is a geodesic
on $S$ passing through $x$ which is mapped to a geodesic in $M$, and the length
metric on $S$ induced by $f$ from $M$ coincides with that induced from the hyperbolic metric on $S$.
When we consider a component $S$ of $\partial C \setminus P$ for some pared manifold $(C,P)$, we always assume that {\em pleated surfaces map frontier components of $S$ to cusps of $M$}. 
%We call such a pleated surface {\sl proper}.
In some situation, we need to consider a pleated surface taking each boundary
component to a closed geodesic.
We call such a pleated surface {\sl a pleated surface with boundary}.
A geodesic lamination $\lambda$ on $S$ is said to be {\sl realised} by a  pleated
surface $f$ when $f|\lambda$ is totally geodesic.

The following dichotomy for measured laminations in the Masur domain was proved
by Otal \cite{OtT}.
See Section 2, above all, Th\'{e}or\`{e}me 2.2 of \cite{OtT}.

\begin{lemma}
Let $(C,P)$ be a pared manifold and $S$ a component of $\partial C \setminus  P$.
Consider a map $g: S \rightarrow M$ to a hyperbolic $3$-manifold $M$ sending
the frontier of $S$ to cusps of $M$.
Let $\lambda$ be a measured lamination contained in the Masur domain of $S$.
Then one and only one of the following holds:
\begin{enumerate}
\item either $\lambda$ is realised by a pleated surface homotopic to $g$ keeping
the frontier mapped
to cusps, 
\item or for any sequence of weighted simple closed curves $\{r_k c_k\}$ on $S$
converging
to $\lambda$, the closed geodesics $c_k^*$ freely homotopic to $g(c_k)$ tend to
an end of $M_0$ after passing to a subsequence.
\end{enumerate}
\end{lemma}

For a geodesic or measured lamination $\lambda$ on $\partial C$, a subsurface of $
\partial C$
 containing $\lambda$ without boundary components null-homotopic on $\partial C$ which is minimal up to isotopies is called the {\sl minimal
supporting surface} of $\lambda$ and is denoted by $T(\lambda)$.
In a special case when $\lambda$ is a (weighted) simple closed curve, we
define its minimal supporting surface to be its annular neighbourhood.
The minimal supporting surface is unique up to isotopy.

We shall next define for two geodesic laminations on the boundary of a 3-manifold to be isotopic.
In contrast to the case of simple closed curves, we need to take more care since without specifying surfaces on which lamination lie, the meaning of isotopy for laminations is not so clear.
Two disjoint minimal geodesic laminations $\lambda_1, \lambda_2$ on $\partial C$ are said to be
{\sl isotopic} when there is a sequence of properly embedded essential annuli
$A_j$ such that
$\partial A_j$ converges in the Hausdorff topology to a geodesic lamination on $T(\lambda_1) \cup 
T(\lambda_2)$ containing $\lambda_1
\cup
\lambda_2$  as $j \rightarrow \infty$.
In the case when no boundary component of $T(\lambda_1)$ and $T(\lambda_2)$ are meridians, this definition of isotopic laminations implies that they lie on homotopic minimal supporting surfaces as we shall see below.

\begin{lemma}
\label{homotopic laminations}
Let $\lambda_1$ and $\lambda_2$ be two disjoint minimal geodesic laminations on $\partial C$ which are isotopic.
Suppose that no boundary component of $T(\lambda_1)$ and $T(\lambda_2)$ is a
meridian.
Suppose moreover that $\lambda_i$ is contained the Masur domain of
$T(\lambda_i)$ for $i=1,2$.
Then the following hold.
\begin{enumerate}
\item The minimal supporting surfaces $T(\lambda_1), T(\lambda_2)$ of
$\lambda_1$ and $\lambda_2$ are both incompressible in $C$.
\item 
 $T(\lambda_1)$  is homotopic to $T(\lambda_2)$ in $C$.
%(This is a consequence of the Jaco-Shalen-Johannson theory. We can apply it to
%$(C, T(\lambda_1) \cup T(\lambda_2))$ since the minimal supporting surfaces are
%incompressible if they have no boundary components which are meridians.)
%Peut-?tre il est mieux d`ajouter le cas o? les deux surfaces co?ncident.
\end{enumerate}
\end{lemma}

\begin{proof}
(i) Suppose that $T(\lambda_1)$ is compressible seeking a contradiction. Let
$\{A_j\}$ be a sequence of essential annuli whose intersection with $
\partial C$ converge to $\lambda_1 \cup \lambda_2$ in the Hausdorff topology. 
As was shown in  Morgan-Shalen \cite{MS2}, there is an incompressible branched surface $B$ properly embedded in $C$ which carries all the annuli $A_i$.
%%pr?ciser la r?f?rence Floyd-Oertel?
Moreover, it was shown there that such a branched surface can be chosen so that $B \cap 
\partial C \subset T(\lambda_1) \cup T(\lambda_2)$, and $B=\tau \times I$ for some train track $\tau$ such that $\tau \times \{0\}$ carries $\lambda_1$ and $\tau \times \{1\}$ does $\lambda_2$.
 By choosing a
train track approximating  $\lambda_1$ closely enough and removing 
branches redundant for carrying $\lambda_1$, we can make $B$  carry no compressing
discs for $T(\lambda_1)$ since $\lambda_1$ is contained in the Masur domain.
Since $B$ is
both incompressible and boundary-incompressible, by the standard cut-and-paste
argument (which uses discs and semi-discs) starting from some compressing disc
for $T(\lambda_1)$, we see that we can reduce the intersection with $B$, and
finally there must be a compressing disc $D$ contained in $T(\lambda_1)$
disjoint from $B$. 
Since $\lambda_1$ is arational in $T(\lambda_1)$ and carried by $B$, this is
possible only when such a disc is isotopic to a boundary component of
$T(\lambda_1)$.
This contradicts our assumption, and we have thus shown that $T(\lambda_1)$ is incompressible.
The same argument works for $T(\lambda_2)$.

(ii)
%Suppose that no boundary components of $T(\lambda_1)$ and $T(\lambda_2)$ are
%meridians.
By (i), we know that both $T(\lambda_1)$ and $T(\lambda_2)$ are
incompressible.
We can assume that they are disjoint  by moving their boundaries slightly by
an isotopy if they intersect each other at their boundaries since $\lambda_1$ and $\lambda_2$ are assumed to be disjoint.
%Since our claim is evident when they coincide, we assume that they are disjoint.
We can then apply the Jaco-Shalen-Johannson theory (\cite{JS}, \cite{Jo}) to $(C, T(\lambda_1) \cup
T(\lambda_2))$ to obtain a characteristic $I$-pair $(\Sigma, \Phi)$
properly embedded in $(C, T(\lambda_1) \cup T(\lambda_2))$ such that any
essential annulus properly embedded in $(C, T(\lambda_1) \cup T(\lambda_2))$ is
properly homotopic into $(\Sigma, \Phi)$.
Since $\lambda_1$ is arational in $T(\lambda_1)$ and there is an essential
annulus whose intersection with $T(\lambda_1)$ approximates $\lambda_1$
arbitrarily closely, we see that $\Phi \cap T(\lambda_1)$ must be isotopic to
$T(\lambda_1)$.
The same holds for $T(\lambda_2)$.
Since $\Sigma$ is an $I$-bundle, this means that  $T(\lambda_1)$ and
$T(\lambda_2)$ cobound a product $I$-bundle  in $C$.
\end{proof}
%Refer also to Claim 6.4 of Lecuire \ref{Lecuire}.
%If there are boundary components bounding  compressing discs, paste compressing
%discs to these components

%We can generalise this lemma to the case when $\lambda$ lies on $\partial_i V$ of the characteristic compression body $V$ of $C$ but not on $\partial C$ as follows.
%%(Note that $\lambda$ may not lie on $\partial C$ in this setting.)
%We say a lamination $\lambda_1$ on $\partial_i V \setminus \partial C$ is homotopic to another lamination $\lambda_2$ on $\partial C$ when there is a sequence of (not necessarily properly) embedded annuli $A_j$ in $C$ whose boundary lies on $\partial_i V \cup \partial C$ such that $\lambda_1 \cup \lambda_2$ is the Hausdorff limit of $\partial A_i$.

%\begin{corollary}
%\label{homotopic lamination 2}
%Let $\lambda_1$ be a lamination on $\partial_i V$ homotopic to another lamination $\lambda_2$ on $\partial C$ which is contained in the Masur domain of $T(\lambda_2)$.
%Then $T(\lambda_1)$ is homotopic to $T(\lambda_2)$ in $C$.
%\end{corollary}
%\begin{proof}
%We consider the annuli $A_i \cap W$ and $A_i \cap V$ for each $i$, whose number of components can be assumed to be constant by passing to a subsequence.
%By applying the argument of the previous lemma for each component, each component of the limit of $A_i \cap W$ or $A_i \cap W$ gives rise to a homotopy between two surfaces lying on $\partial W$ or $\partial V$.
%Connecting them, we get a homotopy between $T(\lambda_1)$ and $T(\lambda_2)$.
%\end{proof}

An {\sl $\reals$-tree} is a geodesic metric space in which two points are
connected by a unique simple arc.
An isometric action of a group $G$ on an $\reals$-tree $T$ is said to have {\sl
small edge-stabilisers}
when for any non-trivial segment $c$ of $T$, the
stabiliser of $c$ is either a finite group or a finite extension of
$\integers$.
Let $S$ be a hyperbolic surface of finite area, and suppose that there is an
action of $\pi_1(S)$ on an $\reals$-tree $T$ by isometries.
A geodesic lamination $\lambda$ on $S$ is said to be {\sl realised} in $T$ when
there is an equivariant map from the universal cover $\hyperbolic^2$ of $S$ to
$T$
which maps each component of the preimages of the leaves of $\lambda$
injectively.
It was proved by Otal \cite{Otal} that if a measured lamination $\lambda$  is
realised, then there is a train track $\tau$ carrying $\lambda$ which is
realised by the same equivariant map: that is, there is an equivariant map
from the universal cover $\tilde{\tau}$ of $\tau$ such that every
branch of $\tilde{\tau}$ is mapped to a non-degenerate
segment and  any train route is mapped locally injectively at every switch of $\tilde \tau$.

For a Kleinian group $G$, we define its deformation space to be the space of
faithful discrete representations of $G$ into $\PSL_2 \complexes$ modulo
conjugacy.
We endow this space with the induced topology as a quotient space of the
representation space with  the topology of point-wise convergence, and denote
it by $AH(G)$.
We denote an element of $AH(G)$, which is a conjugacy class of representations,
by a representation $\phi : G
\rightarrow
\PSL_2 \complexes$ representing the conjugacy class or by a pair $(\Gamma,
\phi)$, where $\Gamma$ is the Kleinian group $\phi(G) \subset \PSL_2
\complexes$.
The quasi-conformal deformations of $G$ modulo conjugacy form a
subspace of $AH(G)$, which is denoted by $QH(G)$.

A sequence of Kleinian groups $\{G_i\}$ is said to converge {\sl
geometrically} to a Kleinian group $H$ if every element 
of $H$ is the limit of a sequence $\{g_i\}$ for $g_i \in 
G_i$, and 
the limit of any convergent sequence $\{g_{i_j} \in 
G_{i_j}\}$ for  a subsequence $\{G_{i_j}\} \subset 
\{G_i\}$ is contained in $H$.
If $\phi_i : G \rightarrow \PSL_2 \complexes$ converges to $\phi$ as representations,
then its images $\phi_i(G)$ converge geometrically (up to extracting a subsequence) to a
Kleinian group containing $\phi(G)$.
When the geometric limit and the algebraic limit coincide, we say that the
convergence is {\sl strong}.

Suppose that a sequence of Kleinian groups $\{G_i\}$ converges to $H$ geometrically.
Let $p_i: \hyperbolic^3 \rightarrow M_{G_i}$ and $q: \hyperbolic^3 \rightarrow M_H$ be universal coverings.
Fix a point $x$ in $\hyperbolic^3$.
Then, $(M_{G_i}, p_i(x))$ converges to $(M_H, q(x))$ in the pointed Gromov-Hausdorff topology:
that is, there exists a $K_i$-bi-Lipschitz map to its image $\rho_i : B_{r_i}(M_{G_i}, p_i(x)) \rightarrow B_{K_ir_i}(M_H, q(x))$  with $K_i \rightarrow 1$ and $r_i \rightarrow \infty$.
In this situation, we say a map $f_i: N_i \rightarrow M_{G_i}$ from some Riemannian manifold  $N_i$ converges geometrically to a map $g: N' \rightarrow M_H$ if $(N_i,z_i)$ converges to $(N',z')$ geometrically with some basepoint $z_i$ which has bounded $d_{M_{G_i}}(p_i(x), z_i)$ and with approximate isometries $\bar \rho_i$, and $\rho_i \circ f_i \circ \bar \rho_i^{-1}$ converges to $g$ uniformly on every compact subset of $N'$.

%Add some facts about geometric convergence plus geometric convergence of maps.
\section{Construction of sequences
and the statement of the convergence 
theorem}
Consider a minimally parabolic, geometrically finite Kleinian group $G$ as was given in Theorem
\ref{main}.
Let $S_1, \dots, S_p$ be the boundary components of $C$ that are not 
tori.
We shall define a sequence of marked conformal structures $\{n_i^k\}$ on
each component $S_k$.
First consider core curves of $P \cap S_k$,  one on each component of $P \cap S_k$, and denote
them by $c^k_1, \dots, c_q^k$ and their union by $C^k$.
On each component of $S_k \setminus C^k$, either a marked conformal
structure or a measured lamination is given in Theorem \ref{main}.
We denote the components of $S_k \setminus C^k$ on which marked
conformal structures are given by $\Sigma^k_1, \ldots , \Sigma^k_s$, and the
given marked conformal structures by $m^k_1, \ldots , m^k_s$.
Let $\mu^k_{s+1}, \ldots , \mu^k_r$ be the measured laminations given on
the rest of the components $\Sigma^k_{s+1}, \ldots , \Sigma^k_r$, and $M^k$ their
union.
Sometimes it is more convenient to consider a compact surface  which is
obtained by deleting a collar neighbourhood of the frontier from
 $\Sigma^k_j$ than $\Sigma^k_j$ itself.
Slightly abusing notation, we use the same symbol $\Sigma_j^k$ to denote
such a compact surface.

\begin{definition}
\label{n}
We define  hyperbolic structures $n^k_i$ on $S_k$ 
in such a way that
\begin{enumerate}
\item $\mathrm{length}_{n^k_i}(c^k_j)= 1/i$ for every $j=1, \dots , q$,
\item $\mathrm{length}_{n^k_i}(\mu^k_j)$ is bounded above independently of $i$
for every $j=s+1,
\dots , r$,
%Change this to "monotone non-increasing, and add another condition that it%%
%diverges to a projective lamination [\mu_j^k]$ lying on the boundary.
\item $\{n^k_i|\Sigma^k_j\}$ diverges in the Teichm\"{u}ller space of $\Sigma_j^k$ towards
the projective lamination $[\mu_j^k]$ on the Thurston boundary for $j=s+1, \dots , r$,
and 
\item $\{n^k_i|\Sigma^k_j\}$ converges to $m^k_j$ as $i \rightarrow \infty$ for every
$j=1,
\dots , s$.
\end{enumerate}
\end{definition}

\begin{lemma}
\label{existence}
Hyperbolic structures $n^k_i$ satisfying the four conditions above exist.
\end{lemma}

\begin{proof}
We shall construct hyperbolic structures on each $\Sigma_j^k$, and get
structures
as
we want by pasting them together.
First we consider $\Sigma^k_j$ for $j=1, \dots , s$.
The conformal structure $m_j^k$ is realised by a hyperbolic structure making
the frontier of $\Sigma^k_j$  cusps.
Such a structure is approximated by  hyperbolic structures with respect to
which  every component of $\partial \Sigma^k_j$ is a closed
geodesic of length
$\delta$ with $\delta \rightarrow 0$.
Therefore we can choose a hyperbolic structure $n^k_i(j)$ converging to
$m^k_j$ on
$\Sigma^k_j$ so that each component of $\partial \Sigma^k_j$ is a closed
geodesic of length
$1/i$.

Next we consider $\Sigma^k_j$ for $j=s+1, \dots , r$.
For each $j$, take a sequence of hyperbolic structures $\{{n'}_i^j\}$
converging 
to a complete hyperbolic structure ${n'}_0^j$ making $\partial \Sigma^k_j$ cusps
such that the length of each component of $\partial \Sigma^k_j$ with respect
to ${n'}_i^j$ is $1/i$.
We should note that since $\{{n_i'}^j\}$ converges to ${n'}_0^j$, the lengths of
$\mu^k_j$ with respect to the ${n'}_i^j$ are bounded above.
We let $n_i^k(j)$ be the hyperbolic structure obtained from ${n'}_i^j$ by the
earthquake along $i\mu^k_j$, which is a measured lamination with the
transverse measure of $\mu^k_j$
multiplied by $i$.
Then we have
$\mathrm{length}_{n_i^k(j)}(\mu^k_j)=\mathrm{length}_{{n'}_i^j}(\mu^k_j)$,
which is bounded above independently of $i$.
Also, since $\{{n'}_i^j\}$ converges and the earthquake is performed along
$i\mu_j^k$, we see that $\{n_i^k(j)\}$ diverges towards the projective
lamination represented by $\mu^k_j$ on the Thurston boundary of
$\mathcal{T}(\Sigma^k_j)$.

%Let $\{w_l \gamma_l\}$ be a sequence of weighted simple closed curves
%converging to $\mu_j^k$ as $l \rightarrow \infty$.
%By considering a pants decomposition,  we can construct a hyperbolic structure
%$n_l(i)$ on $\Sigma^k_j$ such that $\mathrm{length}(\gamma_l)=1/iw_l$ and  each component
%of  $\partial \Sigma^k_j$ is a closed geodesic of length $1/i$.
%Taking the limit of $n_l(i)$ as $l \rightarrow \infty$, we get a hyperbolic
%structure $n_i^k(j)$ with $\mathrm{length}_{n_i^k(j)}(\mu_j)=1/i$.

Pasting these hyperbolic structures $n_i^k(j)$ along the $c^k_j$, we get a
hyperbolic structure as we wanted.
\end{proof}

For a sequence satisfying the four conditions, we have the following.

\begin{lemma}
\label{going to 0}
Let $\{n_i^k\}$ be a sequence of hyperbolic structures on $S_k$ satisfying the
four conditions above.
Then for each $j=s+1, \dots ,r$, there exists a sequence of weighted simple
closed curves
$\{r_ic_i\}$ on $\Sigma^k_j$ with the following two properties.
\begin{enumerate}
\item $\mathrm{length}_{n_i^k}(r_i c_i)$ goes to $0$ as $i \rightarrow
\infty$.
\item $\{r_ic_i\}$ converges to a measured lamination having the same
support as $\mu^k_j$.
\end{enumerate}
\end{lemma}

\begin{proof}
Let $c_i$ be the shortest non-peripheral closed geodesic on $\Sigma^k_j$ with
respect to the metric $n_i^k$.
By Bers' theorem (see \cite{Bers2}), there is a constant $K$ bounding the
$\mathrm{length}_{n_i^k}(c_i)$
from above.
Taking a subsequence, we can assume that there is a sequence of positive
numbers $r_i$ such that $\{r_i c_i\}$ converges to a non-empty measured
lamination $\nu_j^k$ on $\Sigma^k_j$.
Note that $r_i$ goes to $0$ if $c_i$ is not eventually constant.
Even in the latter exceptional case, $r_i$ is bounded above.

We shall show that $\nu_j^k$ and $\mu^k_j$ have the same support.
Suppose that they do not, seeking a contradiction.
Since $\mu^k_j$ was assumed to be arational, then we have $i(\mu_j^k,
\nu_j^k)>0$.
Since $n_i^k|\Sigma_j^k$ diverges towards $[\mu_j^k]$, this implies, by Lemma 8.II.1 in Fathi-Laudenbach-Po\'{e}naru \cite{FLP},
that $\mathrm{length}_{n_i^k}(r_i c_i) \rightarrow \infty$.
Since $r_i$ is bounded above, it follows that $\mathrm{length}_{n_i^k}(c_i)$
also goes to $\infty$.
This contradicts the fact that the lengths are bounded above by $K$.

Since $\mu^k_j$ is not a closed geodesic, neither is $\nu_j^k$, and we see that $r_i \rightarrow 0$.
Hence $\mathrm{length}_{n_i^k}(r_ic_i) \leq r_i K \rightarrow 0$.
Thus we have completed the proof.
\end{proof}

%We first pick up some marked conformal structure $n_0^k$ on each $S_k$
%for $k=1, \ldots , p$.
%% so that its restrictions to $\Sigma^k_1, \ldots , 
%%\Sigma^k_s$ coincide with $m^k_1, \ldots , m^k_s$.
%We define $n_i^k$ to be the marked conformal structure obtained by 
%composing the
%earthquake deformation starting from
% $n_0^k$ with respect to $iM^k$, where the scalar multiplication
%means that of the transverse measure,
% and a  deformation of 
%the structure in a small neighbourhood  of $\Sigma^k_1 \cup \dots \cup
%\Sigma^k_s$,
%which makes the lengths of the components of $C^k$ the
% $i$-th of the original  without changing  the lengths of the
%components of $M^k$
%in such a way that 
%the marked conformal structures on $\Sigma_1, \ldots , \Sigma_s$ converge
%to $m^k_1, \dots , m^k_s$.
%This means that the length of $\mu^k_j$ with respect $n_i^k$ is
%constantly equal to that with respect to $n_0^k$, whereas every simple closed
%curve on $\Sigma_1, \dots , \Sigma_s$ has length with respect to $n_i^k$
%bounded as $i \rightarrow \infty$.
%The latter deformation can be constructed  by considering a pants 
%decomposition of $\Sigma_k$ and shortening only the geodesics on the 
%boundary 
%corresponding to components of $C^k$ without changing the lengths and 
%pasiting maps for other simple closed curves consistuting the pants 
%decomposition.

Let $q : \mathcal{T}(\Omega_G/G) \rightarrow QH(G)$ be the 
Ahlfors-Bers map, which is a (possibly ramified) covering map.
This map is obtained as follows.
For any point $m \in \mathcal{T}(\Omega_G/G)$, there is a Beltrami
differential $\mu$ on $\Omega_G$ giving rise to a quasi-conformal homeomorphism from the original conformal structure on $\Omega_G/G$ to $\mu$, which is automorphic under the action of $G$.
We extend $\mu$ to the entire Riemann sphere by setting the value outside
$\Omega_G$ to be $0$, which is obviously also automorphic.
This defines a quasi-conformal deformation of $G$, which is defined to be
$q(m)$.

Recall that $\Omega_G/G$ is topologically identified with the union of 
non-torus components of $\partial C$.

\begin{definition}
\label{nu_i}
For $n_i^1, \dots , n_i^p$ defined in Definition \ref{n}, we let $\{\nu_i= (n_i^1, \ldots , 
n_i^p)\}$ be  a sequence regarded as lying in $\mathcal{T}(\Omega_G/G)$.
\end{definition}

One of the key ingredients to show our main theorem is the following.

\begin{theorem}
\label{convergence}
The sequence $\{q(\nu_i)\}$ converges in the deformation space 
$AH(G)$, after passing to a subsequence.
\end{theorem}

We shall prove this theorem in the following two sections.
The proof is based on the work of
Kleineidam-Souto \cite{KS}, Lecuire \cite{Lecuire} and Ohshika \cite{Oh1}.
Although there is another way to use  more general result by Kim-Lecuire-Ohshika
\cite{KLO} which is still unpublished, we have chosen to use only \cite{KS},  \cite{Lecuire}  and \cite{Oh1} since the argument in \cite{KLO} is much more complicated than these three.
We note that even if we invoke the main theorem of \cite{KLO}, we cannot obtain Theorem \ref{convergence} immediately and the general line of our argument does not change, although some part can be shorten.

Let $\phi_i : G \rightarrow \PSL_2(\complexes)$ be a geometrically finite
representation corresponding to $q(\nu_i)$, and denote its image by $G_i$,
which is a geometrically finite group.
Since $\phi_i$ is induced by a quasi-conformal deformation, there is a
homeomorphism $\Phi_i : M_G \rightarrow M_{\phi_i(G)}$
inducing the isomorphism $\phi_i$ between the fundamental groups.
%The proof of this theorem consists of two steps. In the first step, we shall
%show that for each component of the 
%characteristic compression body of $C$, the corresponding subgroups of
%$q(\nu_i)$ converge. In the second step, we shall show that subgroups
%corresponding to the 
%complement of the characteristic compression body also converge. We can easily
%show that the entire groups converge once we prove these two facts.

\section{Limit laminations of the boundary-irreducible part}
In \cite{Oh1}, we first analysed the hyperbolic structures on the characteristic compression body.
In contrast, in the present argument, we begin by analysing the behaviour
of the
hyperbolic structures on the complement of the characteristic compression body.

Consider the characteristic compression body $V$ of $C$.
We call the closure of the complement of $V$
the {\sl boundary-irreducible part} of $C$.
Let $W$ be
 a component of the boundary-irreducible part.
Then $\pi_1(W)$ injects to $\pi_1(C) \cong G$ by the homomorphism 
induced from the inclusion.
We denote its image by $H^W$.

We recall the following lemma essentially 
due to Thurston.
Its detailed proof can be found in Theorem 3.1 of Ohshika \cite{Oh1} except for the last sentence concerning the Hausdorff convergence.
It should be also noted that although it was assumed that the Kleinian group does not have parabolic elements in \cite{Oh1}, the existence of rank-$2$ parabolic subgroups does not affect its proof at all.

We note that the codimension-one measured lamination  which we get in this lemma describes only the behaviour of the lengths of the simple closed curves on $\partial W$.
It may not capture the behaviour of the lengths of the closed curves in $W$ in general, hence may not be a lamination dual to the action of $\pi_1(W)$ on an $\reals$-tree which is a rescaled Gromov limit of $(\Gamma_i,\psi_i)$.
Remark 3.2 in \cite{Oh1} gave such an example.

%In the statement, we use the following notation.
%We consider a freely indecomposable, minimally parabolic Kleinian group $H$.
%For a sequence $\{(\Gamma_i, \psi_i)\}$ in $AH(H)$ for a freely indecomposable Kleinian
%group $H$, we denote a homotopy equivalence from $\hyperbolic^3/H$ to
%$\hyperbolic^3/\Gamma_i$
%inducing $\psi_i$ between the fundamental groups by $\Psi_i$.
%Let $W$ be a compact core of $\hyperbolic^3/H$.
%For a measured lamination $\lambda$ on $W$, we denote by
%$\mathrm{length}(\Psi_i(\lambda))$ the sum of the lengths of realisations of
%the components of $\lambda$ by  pleated surfaces in $\hyperbolic^3/\Gamma_i$
%homotopic to $\Psi_i$.
%If some component is not realisable, we define the length of the component to
%be
%$0$.

\begin{lemma}
\label{Thurston}
Let $H$ be a minimally parabolic, freely indecomposable Kleinian group.
Let $W$ be a compact core of $M_H$.
%Let $\mathcal{C}$ be the set consisting of the conjugacy classes of $H$ 
%that are represented by 
%closed curves on $\partial W$.
Suppose that we have a sequence $\{(\Gamma_i, \psi_i) \in AH(H)\}$ 
which does not have a convergent subsequence.
We denote a homotopy equivalence from $M_H$ to
$M_{\Gamma_i}$ inducing $\psi_i$ between the fundamental groups by $\Psi_i$.
Then, after passing to a subsequence of
$\{(\Gamma_i,\psi_i)\}$, there is a sequence of  disjoint, non-parallel unions of 
essential annuli
 $\{A_i^1\sqcup \ldots \sqcup A_i^\kappa\}$ 
properly embedded in $W$, whose number $\kappa$ is independent of 
$i$,  with positive weights $w_i^1, \ldots , w_i^\kappa$ such that $w_i^1
A_i^1 \sqcup \dots \sqcup w_i^\kappa A_i^\kappa$ converges in the measure topology to some
codimension-$1$ measured lamination in $W$ and the following hold.
For any two convergent sequences of measured laminations $\alpha_i,
\beta_i$ lying on  $\partial W$, we have
$$\lim_{i \rightarrow 
 \infty}\frac{\mathrm{length}(\Psi_i(\alpha_i))}{\mathrm{length}(\Psi_i(\beta_i))}=
 \lim_{i \rightarrow \infty}\frac{i(\alpha_i, w_i^1 A_i^1) + \dots +
i(\alpha_i, w_i^\kappa A_i^\kappa)}{i(\beta_i, w_i^1 A_i^1) + \dots +
i(\beta_i, w_i^\kappa A_i^\kappa)},$$
provided that either the numerator or the denominator of the second term tends to
a positive number.
Here, for a measured lamination $\lambda$ on $W$, we denote by
$\mathrm{length}(\Psi_i(\lambda))$ the sum of the lengths of realisations of
the components of $\lambda$ by  pleated surfaces in $M_{\Gamma_i}$
homotopic to $\Psi_i$.
%$$\lim_{i \rightarrow \infty}
%\frac{\mathrm{length}(\psi_i(c))}{\mathrm{length}(\psi_i(c'))}=\lim_{i 
%\rightarrow \infty}\frac{i(w_i^1A^1_i,c )+ \dots 
%+i(w_i^kA^\kappa_i, c )}{
%i(w_i^1A^1_i,c' )+ \dots 
%+i(w_i^kA^\kappa_i, c' )}$$ for every $c, c' \in {\mathcal{C}}$ provided that either the numerator or the denominator of the second term tends to
%a positive number, where
%$\mathrm{length}$ denotes the translation length.
Moreover,  if $i(\alpha_i, w_i^1A^1_i )+ \dots 
+i(\alpha_i, w_i^k A^\kappa_i )$ converges to a positive number, then
$\mathrm{length}(\Psi_i(\alpha_i))$ goes to $\infty$.
The annuli above can be chosen so that $(w_i^1 A_i^1 \sqcup \dots \sqcup
w_i^\kappa A_i^\kappa)$
converges  to a codimension-1 measured lamination whose support is  equal to the Hausdorff limit
of $A_i^1 \sqcup \dots \sqcup  A_i^\kappa$.
%The annuli above can be chosen so that $(w_i^1 A_i^1 \sqcup \dots \sqcup
%w_i^\kappa A_i^\kappa) \cap \partial W$
%converges  to a measured lamination with support equal to the Hausdorff limit
%of $(A_i^1 \sqcup \dots \sqcup  A_i^\kappa) \cap \partial W$ if we regard the
%latter multi-curve as a geodesic lamination.
%(If there are two components in $(A_i^1 \sqcup \dots \sqcup  A_i^\kappa) \cap \partial W$ which are homotopic on $\partial C$, we merge them into one.
%For $(w_i^1 A_i^1 \sqcup \dots \sqcup
%w_i^\kappa A_i^\kappa) \cap \partial W$, we add up all the weights on homotopic curves on $\partial W$ and define it to be the weight of the merged curve.)
\end{lemma}

\begin{proof}
By Theorem 3.1 of \cite{Oh1}, there is a sequence of essential annuli $A_i^1, \dots , A_i^\kappa$ with weights $w_i^1, \dots , w_i^\kappa$ satisfying the conditions of our lemma except for the last sentence.
We only need to show that such annuli can be chosen so that they satisfy the last condition.

Let $L$ be a codimension-1 measured lamination in $W$ which is the limit in the measure topology of $w_i^1 A_i^1 \sqcup \dots \sqcup w_i^\kappa A_i^\kappa$.
As was shown in the proof of Theorem 3.1 in \cite{Oh1} using the theorem of Morgan-Shalen \cite{MS}, there is an incompressible branched surface $B$ carrying both $L$ and $A_i^1 \sqcup \dots \sqcup A_i^\kappa$ for large $i$, which carries only surface with null Euler characteristic.
Consider a weight system $\omega$ on $B$ such that $(B,\omega)$ carries $L$.
We approximate $\omega$ by  rational weight systems $\omega_i$ such that no coordinate of $\omega_i$ goes to $0$ as $i \rightarrow \infty$.
Recall that if $B$ is given a rational weight system, then there is a weighted disjoint union of  essential annuli carried by it.
(If $B$ with the weight system  carries a M\"{o}bius band, we consider an annulus lying on the frontier of a twisted regular neighbourhood of the band which doubly covers it.
We then give the half of the weight on the M\"{o}bius band to the annulus. 
This does not change the effect of  the transverse measure to the simple closed curves on $\partial W$.)
If we let $w_i^1 A_i^1 \sqcup \dots \sqcup w_i^\kappa A_i^\kappa$ be weighted annuli carried by $(B, \omega_i)$, then the Hausdorff limit of $A_i^1 \sqcup \dots \sqcup A_i^\kappa$ coincide the support of $L$ provided that we let the annulus covering a M\"{o}bius band closer and closer to the band as $i \rightarrow \infty$ when we need to consider the Hausdorff limit of M\"{o}bius bands.
\end{proof}

This implies the following corollary. 
%proved in \cite{OhI} using the fact that he
%lengths of $\mu_j$ and of core curves of $P$ with
%respect to $\nu_i$ goes to $0$ and applying Sullivan's theorem proved in
%\cite{EM}.

\begin{corollary}
\label{no V}
Let $\{(G_i, \phi_i)\}$ be a sequence as in Theorem \ref{main}.
If the characteristic compression body $V$ of $C$ is empty (i.e. $W=C$), then $\{(G_i,\phi_i)\}$ converges in $AH(G)$ after
passing to a subsequence.
\end{corollary}

Before starting the proof of Corollary \ref{no V}, we shall consider the general situation when $V$ may not be empty.

Suppose that every subsequence of $\{ \phi_i|H^W\}$ diverges in $AH(H^W)$.
Then by the lemma above, we get a sequence of weighted disjoint, non-parallel union of annuli $a_i=w_i^1 A_i^1 \sqcup \dots \sqcup w_i^\kappa A_i^\kappa$ describing the divergence.
These annuli are disjoint from the torus boundary components of $W$ (passing to a
subsequence if necessary) since the translation length of every parabolic element is $0$.
By taking a subsequence, we can assume that if we regard $a_i$ as a codimension-$1$ measured lamination in $W$,
then it converges to some non-empty codimension-$1$ measured lamination
$\Lambda^W$ in $W$ as $i \rightarrow \infty$.
We call this measured lamination $\Lambda^W$ the {\sl limit lamination} of
$W$.
The limit lamination may depend on the choice of a subsequence which was taken in Theorem \ref{Thurston}.
We fix a subsequence in the following argument.
%The limit lamination is unique up to a scalar multiplication of the transverse
%measure.

%Note that $a_i$ admits a natural $I$-bundle structure since it consists
%of annuli.
By the Jaco-Shalen-Johannson theory (\cite{JS}, \cite{Jo}), all the annuli $A_i^k$ are properly
isotopic into the union of 
the characteristic pairs $(X_j, Z_j)$ of $(W,
\partial W)$ each of which is  either an $I$-pair or a solid torus.
Let $\Lambda_0$ be a component of $\Lambda^W$.
Then there is a characteristic pair $(X_0, Z_0)$ containing $\Lambda_0$
since $\Lambda_0$ is approximated by weighted unions of annuli.
If $X_0$ is an $I$-bundle, we can assume that   the annuli approximating
$\Lambda_0$
are all vertical
 with respect to the $I$-bundle structure of  $X_0$.
Then $\Lambda_0$ is also vertical, and 
admits an $I$-bundle structure itself whose associated $\partial I$-bundle is $\Lambda_0 \cap \partial W$.
(Refer to Morgan-Shalen \cite{MS2} for a detailed account of this.)
If $X_0$ is a solid torus, then $\Lambda_0$ itself is either a weighted annulus or
a weighted M\"{o}bius band.
In either case, $\Lambda_0$ admits an $I$-bundle structure.

Consider an involution $\iota^W$ on $\Lambda^W \cap \partial W$ such that for a
point
$x \in \Lambda^W \cap \partial W$, its image $\iota^W(x)$ is the other endpoint of the
fibre containing $x$ with respect to the
$I$-bundle structure obtained as above.
Let $\lambda_0$ be a component of $\Lambda^W \cap \partial W$.
Then $\iota^W(\lambda_0)$ either coincides with $\lambda_0$ or is disjoint from
$\lambda_0$.
When we talk about $\iota^W(\lambda_o)$, we regard it as having the transverse measure induced from that of  $\lambda_0$.
Moreover, we have the following lemma.

\begin{lemma}
\label{twisted}
The following two hold for the involution $\iota^W$ defined above.
\begin{enumerate}
\item
If $\iota^W(\lambda_0)$ is disjoint from $\lambda_0$, then
their supports are isotopic in $W$.
\item  In the case when $\iota^W(\lambda_0)=\lambda_0$,  there are two
possibilities:
(a) There is a twisted characteristic $I$-pair
over a non-orientable surface $S'$
in $W$ such that $\lambda_0$ is homotopic in $\partial W$ to a double cover of  a measured
lamination on the zero-section of $S'$ (where we regard $I$ as $[-1,1]$).
If $\lambda_0$ is not a simple closed curve, this is the only possibility. 
(b) There is a solid torus component of the characteristic pair
in $W$ such that $\lambda_0$ is homotopic to a double covering of its core
curve.
\end{enumerate}
\end{lemma}

\begin{proof}
Let $\Lambda_0$ be the component of $\Lambda^W$ containing $\lambda_0$.
As was remarked above, $\Lambda_0$ can be assumed to be vertical with respect to
the $I$-bundle structure of the characteristic $I$-pair $(X_0, Z_0)$ containing $\Lambda_0$.
Since the union of annuli $A_i^1 \sqcup \dots \sqcup A_i^\kappa$ converges to
the support of $\Lambda^W$ in the Hausdorff topology by the  last
sentence of Lemma \ref{Thurston},
there is a subset
$a_i'$ of $A_i^1 \sqcup \dots \sqcup A_i$, which consists of  annuli, converging to the support of $\Lambda_0$ in the Hausdorff topology.
If $\lambda_0$ is disjoint from $\iota^W(\lambda_0)$, the $I$-bundle structure of
$\Lambda_0$ must be trivial.
This means that the supports of $\lambda_0$ and $\iota^W(\lambda_0)$ are isotopic by our definition of isotopy between two laminations, and we are done in this case.

We should also note that if $\Lambda_0$ is contained in a characteristic pair
$(X_0, Z_0)$  that is a product $I$-bundle, then the $I$-bundle
structure of $\Lambda_0$ must be trivial.
Therefore, to consider the remaining case when $\lambda_0=\iota^W(\lambda_0)$, we can assume that $(X_0,Z_0)$ is twisted.

Now, suppose that $\iota^W(\lambda_0)=\lambda_0$.
Then the characteristic pair $(X_0,Z_0)$  is either
a twisted $I$-bundle or a solid torus.
First consider the case when $X_0$ is a twisted $I$-bundle.
We consider the base surface $S_0$ of the $I$-bundle, which is a
non-orientable surface, and by identifying $S_0$ with the image of its
section, we regard $S_0$ as embedded in $X_0$ horizontally.
We consider the multi-curve $a_i' \cap X_0$, which is regarded
as a geodesic lamination and is denoted by $\beta_i$.
Then $\beta_i$ converges to a geodesic lamination $\beta_\infty$ over which $|\Lambda_0|$ is a
twisted $I$-bundle.
The bundle must be twisted since we assumed that $\iota^W(\lambda_0)=\lambda_0$.
We can give a transverse measure to $\beta_\infty$ by identifying it with $\Lambda_0 \cap S_0$.
This is the case corresponding to (a) of (ii).

Next suppose that $X_0$ is a solid torus.
Then,  $\Lambda_0$ is either an annulus or a M\"{o}bius
band.
Since $\iota^W(\lambda_0) = \lambda_0$, the latter is the case.
Then obviously $\lambda_0$ is homotopic to a double cover over a core curve of the
M\"{o}bius band.
\end{proof}

%These can be shown by Jaco-Shalen-Johannson theory, and isotoping $\Lambda^W$
%to a horizontal lamination with respect to the $I$-bundle structure.
%(Refer to Morgan-Shalen \cite{MS2}.)
%In the former case, either the minimal supporting surface  of $\lambda$ is
%homotopic to that of  $\lambda'$ in $W$.
%In the latter case, the minimal supporting surface of $\lambda$ is homotopic to a double
%cover of  a subsurface of the base surface $S'$.

Now we shall start a proof of Corollary \ref{no V}.
Although this has already proved in \cite{OhI}, 
the proof here will serve as a perspective for our general argument.

\begin{proof}[Proof of Corollary \ref{no V}]
Suppose, seeking a contradiction, that $\{\phi_i=\phi_i|H^W\}$ diverges and set
$a_i=w_i^1 A_i^1 \sqcup \dots
\sqcup w_i^\kappa A_i^\kappa$ for weighted annuli given in Lemma \ref{Thurston}.
Then $\{a_i\}$ converges to the limit lamination $\Lambda^W$ in the measure topology.
Let $\lambda$ be $\Lambda^W \cap \partial C$, which is disjoint from $T$ 
as was seen above.
Suppose first that one of the laminations  $\mu_j$ given in the
assumption
 of Theorem \ref{main} intersects $\lambda$ essentially.
%By Lemma \ref{going to 0}, there is a sequence of measured laminations
%$\{\lambda_i\}$ on $\Sigma_j$ converging to a measured lamination $\nu_j$ with the same support as $\mu_j$ such that 
%$\mathrm{length}_{\nu_i}(\lambda_i) \rightarrow 0$. 
% 
% Since we assumed that $\mu_j$ intersects $\lambda$ transversely, 
Then, $i(\mu_j, a_i)$ converges to a positive number as $i \rightarrow \infty$ passing to a subsequence.
It follows that  the length of $\Phi_i(\mu_j)$ in $M_{\phi_i(G)}$ goes to
infinity by  Lemma \ref{Thurston}.
On the other hand,  since the length of $\mu_j$ with
respect to $\nu_i$ is bounded,  by Sullivan's theorem (see
Epstein-Marden \cite{EM}) or a generalised version of Bers' inequality (see
Lemma 2.1 in Ohshika
\cite{OhI}), we see that the length of $\Phi_i(\mu_j)$ (on a pleated surface
realising it)
is also bounded. This is a contradiction.

The same argument applies when $\lambda$ intersects a  core curve of a component
of
$P$ essentially.
Therefore we can assume that each component of $\lambda$ is contained in
one of $\Sigma_1, \dots, \Sigma_m$, say $\Sigma_j$.
Suppose that $j \leq n$ first.
There is a simple
closed curve $c$ on $\Sigma_j$ intersecting $\lambda$ essentially.
It follows that the length of $\Phi_i(c)$ goes to infinity by Lemma
\ref{Thurston}.
This contradicts (using Sullivan's theorem again) the assumption that the marked
conformal structures on
$\Sigma_j$ converge to $m_j$, which implies that the length of every
simple closed curve is bounded as $i \rightarrow \infty$.

Therefore $\lambda$ is contained in $\Sigma_{n+1} \sqcup \dots \sqcup
\Sigma_m$.
Since $\lambda$ cannot intersect $\mu_j$ essentially and $\mu_j$ is arational
in $\Sigma_j$, the support of a component of $\lambda$ coincides with that of
some $\mu_j$.
This implies that each component $\lambda_0$ of $\lambda$ is arational in some
$\Sigma_j$;
hence in particular, $\Sigma_j$ is the minimal supporting surface of
$\lambda_0$.
Suppose that $\lambda_0$ is a component of $\lambda$ such that
$\iota^W(\lambda_0)$ is disjoint from $\lambda_0$.
Then $\lambda_0$ is isotopic to $\iota^W(\lambda_0)$.
This implies that the minimal supporting surfaces  $T(\lambda_0)$ and
$T(\iota^W(\lambda_0))$, both of which are among $\Sigma_{n+1}, \dots ,
\Sigma_m$, are homotopic (by (ii) of Lemma \ref{homotopic laminations}), and 
cobound relative
to $P$ a product $I$-bundle since no two distinct components of $P$ are
homotopic.
This is an excluded case in Theorem \ref{main}.

Suppose next that $\iota^W(\lambda_0)=\lambda_0$.
Again the minimal supporting surface of $\lambda_0$ is some $\Sigma_j$.
Then, there is a twisted $I$-bundle bounded relative to
$P$ by $\Sigma_j$, and  $\lambda_0$ doubly covers a lamination on the
zero-section by Lemma \ref{twisted}, as was observed above. 
This is also an excluded case in Theorem \ref{main}.

Thus we have shown that $\{(G_i,\phi_i)\}$ converges in $AH(G)$ in this situation.
\end{proof}

From now on until the end of the proof of Theorem \ref{convergence}, we assume
that $V$ is not empty.
Consider two distinct components $W, W'$ of the irreducible part for which
both $\{\phi_i|H^W\}$ and $\{\phi_i|H^{W'}\}$ diverge in $AH(H^W)$ and
$AH(H^{W'})$ respectively, supposing that there are such components.
%Let $a_i$ and $a_i'$ be the weighted union of essential annuli
%in $W$ and $W'$ respectively obtained by Lemma \ref{Thurston}.
Let $\Lambda^W$ and $\Lambda^{W'}$ be the limit laminations of $W$ and $W'$
respectively.
We say that $W$ {\sl dominates} $W'$  if for closed curves $\gamma$ on $\partial
W$ intersecting $\Lambda^W$ essentially and $\gamma'$ on $\partial W'$
intersecting $\Lambda^{W'}$ essentially, we have $$\lim_{i \rightarrow 
\infty} \frac{\mathrm{length}\phi_i(
\gamma')}{\mathrm{length}\phi_i(\gamma)} = 0,$$
where $\mathrm{length}$ denotes the translation length.
Obviously this definition does not depend on the choice of $\gamma$ and
$\gamma'$.
By taking a subsequence, we can assume that for any two components $W, W'$ of the irreducible parts unless one of them dominates the other, $\displaystyle \lim_{i \rightarrow \infty}\frac{\mathrm{length}\phi_i(
\gamma')}{\mathrm{length}\phi_i(\gamma)} $ exists and is a positive number.
It also follows from the definition that the relation of domination is transitive: if $W$ dominates $W'$ and $W'$ dominates $W''$, then $W$ dominates $W''$.
Therefore we can regard this relation as ordering, and denote $W' \prec_\mathrm{dom} W$ if $W$ dominates $W'$.

Let  $\mathcal{B}$ be the  set of all components of the irreducible part that are maximal with respect to the ordering $\prec_\mathrm{dom}$.
Then by definition, no component contained in $\mathcal{B}$ is dominated by another
component.
Also, by a remark in the previous paragraph, if both $W$ and $W'$ are contained in $\mathcal B$, then $\displaystyle \lim_{i \rightarrow \infty}\frac{\mathrm{length}\phi_i(\gamma')}{\mathrm{length}\phi_i(\gamma)} $ is a positive number,
and  if $W
\in
\mathcal{B}$ and $W'
\not\in \mathcal{B}$ then $W$ dominates $W'$.
We see that $\mathcal{B}$ is  non-empty if $\{\phi_i|H^W\}$
diverges for some component $W$ of the boundary-irreducible part.
We call components contained in $\mathcal{B}$  the {\sl dominating
components}.
We define  $\mathcal{B}$ to be empty if there
is no component $W$ of the boundary-irreducible part such that $\{\phi_i|H^W\}$
diverges.

\section{Convergence of function groups}
%To show the convergence of subgroups of the $G_i$ corresponding to a component of the characteristic compression body of $C$, we use a technique analogous to the one used by Kleineidam-Souto \cite{KS}.
Take a compressible boundary component $S_k$ of $C$.
We denote by $V_k$ the component of the characteristic compression 
body $V$ of $C$ which contains $S_k$ as its exterior boundary.
We consider a subgroup $H^k$ of $G$ associated to the image of 
$\pi_1(S_k)$ in $\pi_1(C) \cong G$ by a homomorphism induced by the
inclusion.
The image of $\pi_1(S_k)$ can be identified with
$\pi_1(V_k)$.

Recall that we denote the union of core curves of $P$ on $S_k$ by $C^k$
and the union of measured laminations $\mu_j^k$ on $S_k$ by $M^k$.
The compression body $V_k$ may have interior boundary components contained in
$\partial C$.
In this case we add core curves of $P$ and the $\mu_j$ on these components to
$C^k$ and $M^k$ respectively.
{\em We denote the unions by $\bar{C}^k$ and $\bar{M}^k$.}
If $\partial_i V_k \cap \partial C = \emptyset$, then we set $\bar{C}^k=C^k$
and $\bar{M}^k=M^k$.

\begin{lemma}
\label{contradicting component}
Suppose that $\{\phi_i\}$ diverges in $AH(G)$ even after passing to a
subsequence.
Then one of the following two conditions holds:
\begin{enumerate}
\item  $\mathcal{B}$ is empty, or
\item  there are a dominating component $W$ with its limit lamination $\Lambda^W$, a component $V_k$ of $V$ meeting $W$,
and a component of $\Lambda^W \cap \partial_i V_k$ whose support is not isotopic in $V_k$ to the support of any component of $C^k\cup M^k$.
\end{enumerate}
\end{lemma}

\begin{proof}
We assume that there is at least one dominating component, and shall show that
the second alternative holds.
If $\Lambda^W$ is disjoint from $V$ for some dominating component $W$ (i.e.\ the
boundary of $\Lambda^W$ lies in $(
\partial C \setminus V)$),
then we can argue entirely
in the same way as the proof of Corollary \ref{no V} and get a contradiction.
Therefore, we can assume that $\Lambda^W \cap V \neq \emptyset$ for every
dominating component $W$.

Let $V_j$ be a component of $V$ intersecting some dominating component.
Suppose that for every dominating component $W$ intersecting $V_j$, each 
component $\lambda$ of the
lamination $\Lambda^W \cap \partial_i V_j$ has support isotopic in $V_j$ to the support
of a component $\lambda'$ of $C^j \cup M^j$.
(Otherwise, we have only to set $V_k$ to be $V_j$ and get the second
alternative.)

Take a component $\lambda$ of $\Lambda^W \cap \partial_i V_j$ for a dominating
component $W$ intersecting $V_j$.
We consider the lamination $\iota^W(\lambda)$.
As was observed in the previous section, either $\iota^W(\lambda)=\lambda$ or
$\iota^W(\lambda)$  is disjoint from $\lambda$.

(a) Let us first consider the case when $\iota^W(\lambda)=\lambda$.\\
Recall that we assumed that there is a component $\lambda'$ of $C^j \cup
M^j$ lying on $S_j$ whose support is 
isotopic  in $V_j$ to that of $\lambda$.
Suppose first that $\lambda'$ is contained in $C^j$.
Then $\lambda$ is a simple closed curve, hence is homotopic to double
covering of either a simple closed curve on the base surface of a twisted
$I$-bundle or
a core curve of a solid torus, either of which is embedded in $W$,
 as was observed in Lemma \ref{twisted}.
This is a contradiction because every component of $C^j$ is a
core curve of
 a paring locus $P$, hence represents a primitive class of $\pi_1(C)$.

Suppose next that $\lambda'$ is contained in $M^j$.
Again by Lemma \ref{twisted}, there is an embedded twisted $I$-bundle $X$ in $W$
and $\lambda$ is homotopic to a double cover of a measured lamination on the
zero-section of the base surface.
By assumption, $\lambda'$ is arational on a component $\Sigma_u$ of $S_j
\setminus C^j$ and is contained in the Masur domain of $\Sigma_u$.
Let $T(\lambda)$ be the minimal supporting surface of $\lambda$ on $\partial_i V_j$.
Then $T(\lambda)$ doubly covers a subsurface $T'$ on the zero-section of the base
surface of $X$.
Since $\lambda$ is isotopic to $\lambda'$, and neither $\partial T(\lambda)$
lying on $
\partial_i V_j$ nor $\partial \Sigma_u$ contains a meridian,
the surface $T(\lambda)$ is homotopic to $\Sigma_u$ in $V_j$ by Lemma
\ref{homotopic laminations}.

Since 
$T(\lambda)$ doubly covers non-orientable $T'$ and every boundary component of $T'$ is an orientation-preserving curve on $T'$,
 each boundary component
of $T(\lambda)$ is one of the two components of the preimage of a boundary component of $T'$, hence is homotopic in $W$ to another boundary component of $T(\lambda)$.
Therefore,  each component of $c$ of $\partial \Sigma_u$ is homotopic in $V_j \cup W$ to  another boundary component of $\partial
\Sigma_u$.
(The boundary of $\Sigma_u$ cannot be empty since  $S_j$ cannot be homotopic to a
surface covered by a subsurface of a component of $
\partial_i V_j$.)
Since the boundary components of $\Sigma_u$ are all in $C^j$, which are core
curves
of  the paring locus $P$, two of them are homotopic in $C$ only
when they are parallel on $\partial V_j$.
If this happens, then the entire $S_j$ is obtained by pasting annuli in $P$ to
$\Sigma_u$, and $V_j \cup X$ itself is a twisted $I$-bundle.
This is impossible since $V_j$ is a compression body.

(b)
Next we consider the case when $\iota^W(\lambda)$ is disjoint from $\lambda$.\\
If $\iota^W(\lambda)$ is also contained in $\partial_i V_j$, then $\iota^W(\lambda)$ is also isotopic in $V_j$ to a component of $C^j \cup M^j$,  and by an argument similar to the above, we are lead to a contradiction with either the assumption that $C^j$
consists of core curves of the paring locus $P$ and the fact that $V_j$ cannot
be contained in a product $I$-bundle if $W$ is not empty.

Suppose now that $\iota^W(\lambda)$ lies on $\partial W \setminus V$, hence on an
incompressible component of $\partial C$.
If either one of the $\mu_j$ or a core curve in $P$ intersects $\iota^W(\lambda)$
essentially, then we get a contradiction as in the proof of Corollary \ref{no
V}.
Therefore, $\iota^W(\lambda)$ is contained in a component $\Sigma$ of $\partial C
\setminus  P$, and by the same argument as in Corollary \ref{no V}, we see that it is arational
there (i.e.\ $T(\iota^W(\lambda))=\Sigma$).
Moreover, since $\Sigma$ lies on $\partial W$, it is incompressible and
$\iota^W(\lambda)$ is contained in the Masur domain of $\Sigma$, which is the
entire measured lamination space.
Then $\Sigma$ and $\Sigma_u$ are homotopic in $C$ by Lemma \ref{homotopic laminations}.
If $\Sigma$ and $\Sigma_u$ have non-empty boundaries, their boundaries, which lie
in distinct components of $\partial C$,
consist of (simple closed curves homotopic to)
core curves of $P$.
This is a contradiction since no two distinct components of $P$ are homotopic in
$C$.
If $\Sigma$ and $\Sigma_u$ are closed surfaces, then $W$ itself is a product $I$-bundle, which contradicts
the definition of the characteristic compression body.

Suppose next that $\iota^W(\lambda)$ is contained in another component $V_k$ of
$V$ with $k \neq j$.
If $\iota^W(\lambda)$ is not isotopic in $V_k$ to a component of $C^k \cup
M^k$, then this $V_k$ is what we were looking for to get the condition (ii).
Suppose, on the contrary, that $\iota^W(\lambda)$ is isotopic in $V_k$ to a component
$\lambda''$ of $C^k \cup M^k$.
Then $\lambda'$ and $\lambda''$ are isotopic in $C$.
If $\lambda''$ is simple closed curve, i.e., contained in $C^k$, then this
means that there
are core curves of distinct two  components of $P$ which are homotopic in $C$.
This is a contradiction.
Suppose that $\lambda''$ is contained in $M^k$.
Then $\lambda''$ is arational in some component $\Sigma_v$ of $S_k
\setminus
C^k$ and $\Sigma_v$ is homotopic to the minimal supporting surface of
$\iota^W(\lambda)$.
(Since no boundary component of $\Sigma_v$ or $T(\iota^W(\lambda))$ is a
meridian, we can use Lemma \ref{homotopic laminations}.)
On the other hand, $\lambda'$ is arational in a component $\Sigma_u$ of $S_j
\setminus
C^j$ which is homotopic to the minimal supporting surface $T(\lambda)$.
Since $T(\lambda)$ is homotopic to 
$T(\iota^W(\lambda))$, this implies that $\Sigma_u$ is homotopic to
$\Sigma_v$.
In particular,  a boundary component of $\Sigma_u$, which is contained in
$C^j$,
is
homotopic to a boundary component of $\Sigma_v$, which is contained in $C^k$.
(If $\Sigma_u$ is closed, we get a contradiction as before.)
This contradicts the fact that no two core curves of distinct  components of $P$
are homotopic.
Thus we have completed the proof.
\end{proof}

Now we start the proof of Theorem \ref{convergence}.
We shall first show that for each component $V_k$ of $V$, the 
sequence $\{(\phi_i(H^{k}), \phi_i|H^{k}) \in AH(H^{k})\}$ converges if
$\mathcal{B}$
is empty, using Lecuire's result on measured laminations in $\mathcal{D}$.
For  that, we need to prove that  we can extend $\bar{C}^k \cup
\bar{M}^k$
 to a measured lamination contained in $\mathcal{D}(V_k)$.

\begin{lemma}
\label{extension}
Take an arational measured lamination on each component of   $(V_k \cap \partial C)
\setminus P$ among $\Sigma_1, \dots , \Sigma_n$, which is contained in its Masur domain 
and does not have  support isotopic in $V_k$ to the support of a  component of $M^k$ in such a way that no two of them have
supports isotopic to each other in $V_k$.
(Recall that $\Sigma_1, \dots , \Sigma_n$ are components of $\partial C \setminus P$ where the marked conformal structures converge in the Teichm\"{u}ller spaces.)
Denote the union of all of these arational measured laminations by $L^k$.
Take also an arational measured lamination on each interior boundary component of
$V_k$ that is not contained in $\partial C$,
and let $Q^k$ be their union. 
Then $\bar{C}^k \cup \bar{M}^k \cup L^k \cup Q^k$ is contained in 
$\mathcal{D}(V_k)$.
\end{lemma}

\begin{proof}
By the definition of $\mathcal{D}$, we have only to show that there 
is $\eta>0$ such that  every meridian and the boundary of every essential annulus
has intersection
 number greater than $\eta$ with $\bar{C}^k \cup \bar{M}^k \cup L^k \cup Q^k$.
Suppose, seeking a contradiction, that there is a sequence of
 meridians $\{\partial D_l\}$  or of boundaries of essential annuli $\{
\partial A_l\}$ such that
$i(\partial D_l,
\bar{C}^k
\cup \bar{M}^k \cup L^k \cup Q^k)$ or $i(\partial A_l, \bar{C}^k
\cup \bar{M}^k \cup L^k \cup Q^k)$ goes to $0$ as $l \rightarrow
\infty$. 
Then after taking a subsequence, we can assume that $\partial D_l$ or
$\partial A_l$
is
 disjoint from $\bar{C}^k$ for every $l$, for otherwise it has intersection
number at least $1$ with $\bar{C}^k$.
For meridians, we set $c_l = \partial D_l$, and for annuli, we set $c_l$ to be
one of the components of $
\partial A_l$.
Then, passing to a subsequence, all of the $c_l$ can be assumed to be contained
in
 a component $\Sigma$ of $
\partial V_k \setminus \bar{C}^k$ since they
are disjoint from $\bar{C}^k$ and there are only finitely many components of
$\partial V_k \setminus \bar{C}^k$.
We can also assume that the same holds even if we choose the other boundary component of $A_l$ for each $l$.
Since $\mu_j$ or a component of $L^k \cup Q^k$ given on $\Sigma$ is  
arational, it follows
 that its support coincides with a minimal component 
$\ell$ of 
the Hausdorff limit of the $c_l$, for both choices of the boundary components of $A_l$.
If the $c_l$ are meridians (of $\Sigma$),  their
Hausdorff
 limit contains a homoclinic leaf.
(Refer to Th\'{e}or\`{e}me 1.8 of Otal \cite{OtT}.)
By Lemma \ref{no homoclinic}, this contradicts the assumptions that any $\mu_j$
or
 any component of $L^k$
is contained in the Masur domain, and that $Q^k$ lies on the interior boundary of
$V_k$.

Next suppose that the $c_l$ are boundary components of the $A_l$.
We consider first the case when the two boundary components of $A_l$ are
contained
in distinct components $\Sigma$ and $\Sigma'$ of $\partial V_k
\setminus \bar{C}^k$.
Let $\ell, \ell'$ be the Hausdorff limits of $A_l \cap \Sigma$ and
$A_l \cap \Sigma'$, which are isotopic in $V_k$ by definition.
Since  the Hausdorff limit of $\partial A_l$ does not intersect
$\bar{M}^k \cup L^k \cup Q^k$ essentially, both $\ell$ and $\ell'$ consist of a
unique arational minimal component and isolated non-compact leaves spiralling around them.
Their minimal supporting surfaces are $\Sigma$ and
$\Sigma'$ respectively.
Also, the minimal components of both $\ell$ and $\ell'$ carry transverse measures
which give rise to
measured laminations in the Masur domains of $\Sigma$ and $\Sigma'$ since they
coincide with the supports of components of $\bar{M}^k \cup Q^k \cup L^k$.
Since their boundaries do not contain meridians, by Lemma \ref{homotopic
laminations},
$\Sigma$ and $\Sigma'$ are homotopic.
Moreover, since no two distinct components of $\bar{C}^k$ are homotopic and no
interior boundary component is homotopic to another boundary component,
this is possible only
when  $\Sigma$ and $\Sigma'$ lie on $\partial V_k \cap \partial C$, all the
boundary components of $\Sigma$
 are homotopic in $
\partial V_k$ to
boundary components
of $\Sigma'$, and the same holds interchanging $\Sigma$ with $\Sigma'$.
This happens only when $V_k$ is a handlebody and homeomorphic to a product
$I$-bundle over a compact surface as a pared manifold.
This in particular implies that the boundary-irreducible part is empty; hence
$V_k = C$, and $\Sigma,  \Sigma'$ can be assumed to be
$\Sigma_1, \Sigma_2$.

By the assumption of Theorem \ref{main}, if the $\mu_j$ are given on both of
$\Sigma_1$ and $\Sigma_2$, i.e. $n=0$, then the supports of $\mu_1$ and $\mu_2$
are not isotopic.
This contradicts the fact that $\ell$ and $\ell'$ are isotopic and the
supports of  $\mu_1, \mu_2$ are their minimal components.
If $\mu_j$ is not given on  $\Sigma_1$, i.e., $n\geq 1$, then by our assumption,
the component $\lambda_1$ of $L^k$ on $\Sigma_1$ has support which is not
isotopic to that of 
$\mu_2$ or $\lambda_2$. 
This again contradicts the fact that $\ell$ and $\ell'$
are isotopic.
 
Next we consider the case when both of the boundary components of $A_l$ lie in
the same component $\Sigma$ of $\partial V_k  \setminus
\bar{C}^k$.
We need to divide our argument into two sub-cases depending on whether
$\Sigma$ is compressible or not.
Suppose first that $\Sigma$ is compressible.
Let $D$ be a compressing disc of $\Sigma$.
Then, since $\{\partial A_l\}$ converges to an arational lamination, $A_l$
intersects $D$ essentially for large $l$.
Therefore, we can either boundary-compress $A_l$ along an outermost semi-disc
on $D$ bounded by $A_l \cap D$, or get a compressing disc intersecting $A_l$
along fewer arcs than $D$ which is obtained by cutting $D$ along an arc cobounding with
$\partial A_l$ an
outermost semi-disc on $A_l$.
In the former case, we get a meridian $d_l$ of $\Sigma$ which is disjoint from
$A_l$ by
boundary-compressing $A_l$.
Also in the latter case, by repeating the same operation, we eventually get a
meridian $d_l$ of $\Sigma$ which is disjoint from $A_l$.
Then the Hausdorff limit $d'$ of the $d_l$ does not intersect the Hausdorff limit
of $\partial A_l$ transversely, hence neither does it a component of $\bar{M}^k \cup
L^k \cup Q^k$ on $\Sigma$, which is
an arational lamination.
Then by the same argument using the homoclinicity as for the Hausdorff limit of
$\{D_l\}$ above, we get a contradiction.

Suppose next that $\Sigma$ is incompressible.
Then by applying the Jaco-Shalen-Johannson theory to $(V_k, \Sigma)$, we
see that there is a characteristic pair $X$ which is a union of  $I$-bundles and
solid tori intersecting $
\partial V_k$ in $\Sigma$, to which every essential annulus with boundary on
$\Sigma$ can be
properly homotoped.
Since $\{\partial A_l\}$ converges to an arational lamination, 
there is an $I$-bundle component $X_0$ of $X$ to which $A_l$ can be properly
homotoped for large $l$, and $ \Sigma$ is isotopic in $V_k$ to the
associated $\partial I$-bundle  of $X_0$.
This is possible only when $X_0$ is a twisted $I$-bundle.
We also see that $V_k = X_0$   since no two
 components of $\partial \Sigma$ not homotopic on $\partial V_k$ are homotopic
in $V_k$.
This implies that $V_k$ is a twisted $I$-bundle over a non-orientable
surface, and contradicts the fact that $V_k$ is a
compression body.
\end{proof}

We need to use the following lemma, which was proved by Lecuire in
\cite{Lecuire}.

\begin{lemma}
\label{Lecuire}
Let $V$ be a compression body and $S$ its exterior boundary.
Suppose that $\pi_1(V)$ acts on an $\reals$-tree $\mathcal{T}$ by isometries 
with small edge-stabilisers.
Let $\mu$ be a measured lamination contained in $\mathcal{D}(V)$.
Then there exists a $\pi_1(S)$-equivariant map $F: \hyperbolic^2 
\rightarrow \mathcal{T}$ which realises at least one component of $\mu$.
(Here we regard $\pi_1(S)$ as acting on $\mathcal{T}$ by pre-composing the
epimorphism from $\pi_1(S)$ to $\pi_1(V)$ induced by the inclusion.)
\end{lemma}

%Now, recall that we have a sequence of quasi-conformal deformations 
%of $G$, which was given as $q(\nu_i)$.
%Let $\phi_i : G \rightarrow \PSL_2(\complexes)$ be a geometrically finite
%representation corresponding to $q(\nu_i)$.

As a first step of the proof of Theorem \ref{convergence}, we shall show the 
following proposition.
This is a case which Lecuire's Theorem 6.6 in \cite{Lecuire} already covers.
Still, we shall give an outline of proof here based on Lecuire's lemma above so that we
can refer to it in the argument for the next case when $\mathcal{B}$ is not empty.

\begin{proposition}
\label{convergence of function groups}
Suppose that $\mathcal{B}$ is empty.
Then for every component $V_k$ of $V$, the quasi-conformal deformations of $H^k$
given by the $\phi_i|H^k$ converge after taking a subsequence and  
conjugates.
\end{proposition}

\begin{proof}
Suppose, seeking a contradiction, that $\{\phi_i|H^k\}$ does not 
have a convergent subsequence in $AH(H^k)$.
Then, by the Morgan-Shalen-Bestvina-Paulin theory (\cite{MS}, \cite{Be} and
\cite{Pa}), there is an isometric
 action with small edge-stabilisers of $H^k$ on an $\reals$-tree $\mathcal{T}$
which is a Gromov limit
of the rescaled action of $\phi_i(H^k)$ on $\epsilon_i \hyperbolic^3$ with
$\epsilon_i \rightarrow 0$.
This can be regarded as an action of $\pi_1(V_k)$ on $\mathcal{T}$.
By Lemmata \ref{extension} and \ref{Lecuire}, one of the components 
of $\bar{C}^k \cup \bar{M}^k \cup L^k \cup Q^k$, which we shall denote by $\chi$,
is realised
 in $\mathcal{T}$.

Suppose first that $\chi$ is a component of $\bar{C}^k$.
Then there is an element $\gamma \in H^k$ whose conjugacy class is 
represented by $\chi$.
Since $\chi$ is realised in $\mathcal{T}$, the translation length of 
$\phi_i(\gamma)$ goes to infinity as $i \rightarrow \infty$.
On the other hand, by our construction, the length of $\chi$ with 
respect to the conformal structure at infinity $\nu_i$ goes to $0$.
This implies, as was shown in Theorem 6.2 Sugawa \cite{Su} (or by the main theorem of Canary \cite{CaD}), that
the translation length of $\phi_i(\gamma)$ also goes to $0$.
Thus we are lead to a contradiction.

Next consider the case when $\chi$ is a component $\mu$ of $\bar{M}^k$.
By Lemma \ref{going to 0}, there is a sequence of weighted simple closed
curves $\{r_nc_n\}$ on a component of $S_k \setminus C^k$
such that $r_n c_n$ converges to
 $\mu$ and 
$\mathrm{length}_{\nu_i}(r_i c_i)$ goes to $0$ as $i \rightarrow 0$.
We should also note that for any meridian $m$ in $\Sigma$, its length with
respect to $\nu_i$ goes to $\infty$ since $i(m,\mu)>0$.
Therefore, by the main theorem of Canary \cite{Ca},  
$r_i\mathrm{length}(\phi_i(c_i))$ in $M_{\phi_i(H^k)}$  goes to $0$,
where $\mathrm{length}$ denotes the translation length.
Since $\chi$ is arational, the Hausdorff limit of $c_n$, which we denote by
$\xi$,
is the union of $\chi$ and finitely many isolated non-compact leaves.
Since $\chi$ is realised in $\mathcal{T}$, by the same argument as
Th\'{e}or\`{e}me 3.1.4 of Otal \cite{Otal}, we see that $\xi$ is also realised
in $\mathcal{T}$.
Furthermore, Otal's argument implies that there is a train track carrying $\xi$
which is realised in $\mathcal{T}$.
This 
implies that  $r_i\mathrm{length}(\phi_i(c_i))$ in $M_{\phi_i(H^k)}$
must go to
 infinity as $i \rightarrow \infty$ 
since $\tau$ is mapped to a train
track with geodesic branches and small exterior angles at switches in
$M_{\phi_i(H^k)}$
(see Chapitre 3 of
\cite{Otal}).
This is  a contradiction.

Next we consider the case when $\chi$ is a component of $L^k$.
Then $\chi$ lies on $\Sigma_j$ on which the conformal structure converges to
$m_j$.
Therefore, the length with respect to $n_i^k$ of $r_ic_i$ taken as before is
bounded 
as $i \rightarrow \infty$.
Thus, by the same argument as the previous case, we get a 
contradiction.

Finally, we consider the case when $\chi$ is a component of $Q^k$.
 Let $\Sigma$ be the interior boundary component not contained in $\partial C$ on
which $\chi$ lies.
 There is a component $W$ of the boundary-irreducible part containing $\Sigma$
as a boundary component.
 Since we assumed that $\mathcal{B}$ is empty, $\phi_i|H^W$ converges, hence
in particular, the length of  $\Phi_i(\chi)$ in
$M_{\phi_i(G)}$ (hence also that in $M_{\phi_i(H^k)}$)
is also bounded.
 This contradicts the fact that $\chi$ is realised in $\mathcal{T}$ by the
same argument as above.
\end{proof}

Now we consider the case when $\mathcal{B}$ is not empty.
By Lemma \ref{contradicting component}, there are a dominating component $W$
and a
component $V_k$ of
$V$ such that there is a
component of
 $\Lambda^W
\cap
\partial_i V_k$ whose support is not isotopic in $V_k$ to the support of a component of
$\bar{C}^k
\cup \bar{M}^k$.
We shall show that this will lead to a contradiction.

\begin{proposition}
\label{converge if contradicting}
In the settings of Theorem \ref{convergence}, if $V$ is not empty, then $\mathcal B$ must be empty.
\end{proposition}

\begin{proof}
To prove this proposition, we need to analyse an action on an $\reals$-tree in the
same way as the proof of Lemma \ref{Lecuire} in \cite{Lecuire}.

Suppose, seeking a contradiction, that neither $V$ nor $\mathcal B$ is  empty.
Then there is  a component $V_k$ of $V$ intersecting $\Lambda^W$ as we mentioned just before the proposition.
Then $\{\phi_i|H^k\}$ must diverge in $AH(H^k)$ since $\Lambda^W$ intersects $\partial_i V_k$.
By the same argument as Proposition \ref{convergence of function groups},
there is a  limit isometric action $\rho$ of $H^k\cong \pi_1(V_k)$ on an
$\reals$-tree $\mathcal{T}$ having small edge-stabilisers.
We first consider a special case when the restriction of $\rho$ to every
interior boundary component of $V_k$ has a global fixed point in $\mathcal{T}$.
In this case, we can argue as in the proof of Proposition \ref{convergence of
function groups}.
We extend $\bar{C}^k \cup \bar{M}^k$ by adding $L^k$ and $Q^k$ using Lemma
\ref{extension}.
By Lemma \ref{Lecuire}, one of the components of $\bar{C}^k \cup \bar{M}^k \cup
L^k \cup Q^k$ is realised by $\rho$.
By the same argument as the proof of Proposition \ref{convergence of function
groups}, we see that it is impossible that a component of $\bar{C}^k \cup \bar{M}^k \cup
L^k$ is realised by $\rho$.
If a component $\chi$ of $Q^k$ is realised by $\rho$, then we consider the
component $\Sigma$ of $\partial_i V_k$ on which $\chi$ lies.
Then the restriction of $\rho$ to $\pi_1(\Sigma)$ is non-trivial.
This contradicts the assumption of our special case here.

Now, we assume until the end of the proof that there is at least one component
of $\partial_i V_k$ on which the restriction of $\rho$ is non-trivial, i.e.,
does not have a global fixed point.
Regard  $\rho$ as an action of $\pi_1(S_k)$ (recall that $S_k=\partial_e V_k$) by pre-composing the 
epimorphism from $\pi_1(S_k)$ to $\pi_1(V_k)$ induced by the inclusion.

Take a sequence of weighted multi-curves $\{c_n\}$ decomposing $S_k$ into pairs of
pants
 which converges to $C^k \cup M^k \cup (L^k \cap S_k)$ in $\mathcal{ML}(S_k)$ as $n \rightarrow
\infty$.
Since $C^k \cup M^k \cup (L^k \cap S_k)$ is maximal, by approximating it by a
train track with complementary regions whose vertices correspond one-to-one to
ideal vertices of the complementary regions of $C^k \cup M^k \cup (L^k \cap
S_k)$ and giving rational weights, we can assume that a union of some
components of
$|c_n|$ (the support of $c_n$) converges to the support of $C^k \cup M^k \cup (L^k \cap S_k)$ with
respect to the Hausdorff topology.
As was shown in  Morgan-Otal \cite{MO}, there are an action $\rho_n$ of
$\pi_1(S_k)$
on an $\reals$-tree $\mathcal{T}_n$, which is dual to a measured lamination
$\zeta_n$ on $S_k$, and a morphism $\pi_n$ from $\mathcal{T}_n$ to $\mathcal{T}$ 
such that the action obtained by pushing $\rho$ forward by $\pi_n$, which we denote by $(\pi_n)_* \rho_n$, coincides
with $\rho$, and the translation length of $c_n$, \ie the weighted sum of the translation lengths of the components of $c_n$, with respect
to $\rho_n$ is equal to that with respect to  $\rho$.
Now we consider the Hausdorff limit $\zeta_\infty$ of $\{|\zeta_n|\}$.
As was shown in the proof of Theorem 3 in Kleineidam-Souto \cite{KS}, if a component of $C^k\cup M^k \cup
(L^k
\cap S_k)$ intersects
$\zeta_\infty$ essentially, then it is realised in
$\mathcal{T}$.
This is a contradiction as in the proof of Proposition \ref{convergence of
function groups}.

Therefore, $\zeta_\infty$ does not intersect $C^k \cup M^k \cup (L^k \cap S_k)$
essentially.
%We shall first show that $\mu_\infty$ is disjoint from $(L^k \cap S_k)$.
%If not, there is a minimal component $\mu'$ of  $\mu_\infty$ having the same
%support as a component $\lambda_j$ of 
%$(L^k
%\cap S_k)$.
%Recall that $\lambda_j$ was given on a component $\Sigma_j$ of $S_k \setminus P$,
%which is among $\Sigma_1, \dots , \Sigma_n$.
%Since a component of the support of $c_n$ converges to the support of $\lambda_j$
%with respect to the Hausdorff topology,
%for sufficiently large $n$, there is a component $\gamma_n$ of the support 
%of $c_n$ which is contained in $\Sigma_j$.
%Since $\gamma_n$ intersects $\mu_n$ transversely for large $n$, it must have a
%positive translation number in 
%On the other hand, by our definition of $\nu_i$ and Sugawa's inequality, the
%translation length of $\phi_i(c)$ is bounded as $i \rightarrow \infty$.
%This is a contradiction.
%
%%Let $\partial_i V_k$ denote the interior boundary of $V_k$.
We consider the measured lamination $\Lambda^W \cap \partial_i V_k$ for all the
components $W \in \mathcal{B}$, and let $\hat{\upsilon}$ be the union of all such laminations on $\partial_i V_k$.
By Lemma \ref{contradicting component}, there are a component $F$ of
$\partial_i V_k$  and a component $\upsilon_0$ of $\hat{\upsilon}$
having support which is not
isotopic in $V_k$ to the support of a component of $C^k \cup M^k$.
Moreover, we see that no component of $\hat{\upsilon}$ is isotopic in $V_k$ to a
component of $L^k \cap S_k$ as follows.

Let $\lambda_j$ be a component of $L^k \cap S_k$, and suppose that $\lambda_j$
lies on 
$\Sigma_j$ which is among $\Sigma_1, \dots , \Sigma_n$.
Then we can take a simple closed curve $c$ on $\Sigma_j$ intersecting
$\lambda_j$ essentially.
Suppose, seeking a contradiction, that  $\lambda_j$ is isotopic in $V_k$ to a
component $\upsilon_0$
of $\hat{\upsilon}$.
Then the minimal
supporting surface $T(\upsilon_0)$ is homotopic to $\Sigma_j$ in $V_k$ by Lemma
\ref{homotopic laminations}.
This implies that there is a simple closed curve $c'$ on $T(\mu_0)$ homotopic
to $c$ intersecting $\upsilon_0$ essentially.
In particular, the translation length of $\phi_i(c')$ goes to
infinity as $i \rightarrow \infty$.
By using the fact  that the length of $c$ with respect to
$\nu_j$ is bounded, and either Theorem 6.2 in Sugawa \cite{Su} or the main theorem of Canary \cite{CaD}, we get a contradiction.

Recall that we assumed that for some component of $\partial_i V_k$, the
restriction of $\rho$ to the component is non-trivial.
Consider the component $F$ of $\partial_i V_k$ containing $\upsilon_0$ as defined above.
Since all the $W$ in $\mathcal{B}$ are  dominating and some $\Lambda^W$ with $W \in \mathcal B$
intersects $F$, if the
restriction of $\rho$ to $F$ is trivial, then the restriction of $\rho$
to every component of $\partial_i V_k$ is trivial.
Therefore, under the present assumption, the restriction of $\rho$ to $\pi_1(F)$ is
a non-trivial action having small edge-stabilisers since $F$ is incompressible.
By
Skora's theorem \cite{Sk}, it is dual to some measured lamination $\nu$
on $F$.
By the definition of $\Lambda^W$, if we consider the restriction of $\phi_i$
to the subgroup corresponding to $\pi_1(F)$ and its rescaled Gromov limit, we see that $\nu$ must
coincide with $\hat{\upsilon} \cap F$ up to a scalar multiple.
%On the other hand, since there is a morphism from $\mathcal{T}_n$ to
%$\mathcal{T}$, the support of  $\mu_n$ contains that of $\hat{\mu}$; hence so
%does $\mu_\infty$.
%
%Recall that by Lemma \ref{contradicting component}, we can choose a component $F$ of $\partial_i V_k$ such that
%$\hat{\upsilon} \cap F$ has a component $\upsilon_0$ whose support is not
%isotopic in $V_k$ to the support of  a component of $C^k
%\cup M^k$, and that such a component cannot be isotopic to $L^k \cap S_k$ either as was shown above.
Since $V_k$ is a compression body, there is a surface $\hat F$ on $S_k$ each of whose boundary component is a meridian such that the surface obtained by attaching disjoint compressing discs to $\partial \hat F$ is isotopic to $F$ in $V_k$.
Let $\upsilon_1$ be a measured lamination on $\hat F$ which is isotopic to $\upsilon_0$.

Recall that we have actions $\rho_n$ of $\pi_1(S_k)$ on $\reals$-trees $\mathcal T_n$, which are dual to measured laminations $\zeta_n$ on $S_k$.
As was shown in the proof of Proposition 6.1 in Lecuire \cite{Lecuire}, we see
that $|\zeta_n|$
is constant for large $n$, which implies that $\zeta_\infty$ is the support of a
measured lamination, which we denote by $\zeta$.
As was seen above, $\zeta_\infty$ cannot intersect $C^k \cup M^k \cup (L^k \cap S_K)$ essentially; hence we have $i(C^k \cup M^k \cup (L^k \cap S_K), \zeta) =0$.
%If $i(C^k \cup M^k \cup (L^k \cap S_K), \zeta) >0$, then there is a component of $C^k \cup M^k \cup (L^k \cap S_K)$ which is realised by $\rho$ in $\mathcal T$ by the proof of Theorem 3 in Kleineidam-Souto \cite{KS}.
%$i(c_n, \zeta_n)$, which is equal to the translation length  of $\rho(c_n)$ in $\mathcal T$, is bounded below by a positive constant, and we have a subsequence $\{\phi_{i_n}\}$ of $\{\phi_i\}$ such that $\length(\phi_{i_n}(c_n))$ goes to $\infty$.
%This is a contradiction as before.
Since $C^k\cup M^k \cup (L^k \cap S_k)$ is maximal, we see that the support of $\zeta$ is contained in that of $C^k\cup M^k \cup (L^k \cap S_k)$.
In particular, $|\zeta|$ cannot contain $|\upsilon_1|$ as a component since $\upsilon_1$ intersects $C^k\cup M^k \cup (L^k \cap S_k)$ essentially.

%We shall show that there exists a component of $\zeta$ whose support coincides with that of $\upsilon_1$.
We shall now show that $\zeta$ cannot intersect $\upsilon_1$ essentially.
Suppose, on the contrary, that there is a component $\zeta_0$ of $\zeta$ intersecting $\upsilon_1$ essentially.
We shall first consider the case when $\zeta_0$ is contained in $\hat F$.
Since $\rho_n$ is dual to $\zeta_n$ whose support is equal to $|\zeta|$ for large $n$, there is a $\pi_1(S_k)$-equivariant map $q_n: \hyperbolic^2 \rightarrow \mathcal T_n$ which maps each component of the preimage of $\zeta$ to a point in $\mathcal T_n$.
Then there is a train track $\tau$ carrying $\zeta_0$ each component of whose preimage is mapped to a point by $q_n$, as was shown in the proof of Theorem 4 in Kleineidam-Souto \cite{KS}.
Let $s_n k_n$ be a sequence of weighted simple closed curves converging to $\zeta_0$ which is carried by $\tau$.
(Such a sequence can be taken by approximating the weight system for $\zeta_0$ by rational numbers.)
Then each lift of $k_n$ in $\hyperbolic^2$ is mapped to a point by $q_n$; hence the translation length of $\rho(k_n)=(\pi_n)_* \rho_n(k_n)$ is $0$.
On the other hand, since $\zeta_0$ intersects $\upsilon_1$, by homotoping $\zeta_0$ and $\upsilon_1$ in $V_k$ into $F$, we see that the translation length of $\rho(k_n)$ must go to $\infty$ if $\zeta_0$ is not a simple closed curve, and is positive if $\zeta_0$ is a simple closed curve.
This is a contradiction.

Next we consider the case when $\zeta_0$ is not contained in $\hat F$.
Since $\zeta$ is a measured lamination, the leaves of $\zeta$ intersecting $\partial \hat F$ cannot accumulate inside $\hat F$, and consist of finitely many parallel families of geodesic arcs with endpoints on $\partial \hat F$.
By joining outermost geodesic arcs in parallel families and arcs on $\partial F$ disjoint from $\zeta$, we can construct essential simple closed curves on $\hat F$ which are disjoint from $\zeta$.
Since $\zeta_0$ intersects $\upsilon_1$ essentially, so does one of the simple closed curves, which we denote by $\gamma$.
This curve $\gamma$ is homotopic to a simple closed curve $\gamma'$ on $F$ which intersects $\upsilon_0$ essentially.
Since the restriction of $\rho$ to $\pi_1(F)$ is dual to $\nu$ whose support contains that of $\upsilon_0$, we see that $\rho(\gamma)$ has a positive translation length in $\mathcal T$.
This contradicts the fact that $\rho_n(\gamma)$ has null-translation length in $\mathcal T_n$ (since $\gamma$ is disjoint from $\zeta$) and $\rho=(\pi_n)_* \rho_n$.

Thus we have shown that $\zeta$ cannot intersect $\upsilon_1$ essentially.
The only remaining possibility is that $\upsilon_1$ is disjoint from $\zeta$.
In this case, we consider a complementary region $U$ of $\zeta$ containing $\upsilon_1$.
We note that $U \cap \hat F$ contains $T(\upsilon_1)$.
If $\upsilon_1$ is not a weighted simple closed curve, then there is a simple closed curve $\delta$ in $T(\upsilon_1)\subset U \cap \hat F$ intersecting $\upsilon_1$ essentially.
Since $\delta$ is disjoint from $\zeta$, the translation length of $\rho_n(\delta)$, hence also that of $\rho(\delta)$, is $0$.
On the other hand, $\delta$ is homotopic to a simple closed curve $\delta'$ on $F$ intersecting $\upsilon_0$ essentially.
Therefore $\rho(\delta)$ must have positive translation length in $\mathcal T$ as before.
This is a contradiction.

Next suppose that $\upsilon_1$ is a weighted simple closed curve.
By our choice of $\upsilon_0$, there is a component $\xi$ of $C^k\cup M^k \cup (L^k \cap S_k)$ intersecting $\upsilon_1$ essentially.
%Since we are assuming that $\zeta$ is disjoint from $\upsilon_1$ and $|\zeta|$ is contained in the support of $C^k\cup M^k \cup (L^k \cap S_k)$, 
Since $\zeta$ is assumed to be disjoint from $\upsilon_1$, there is a complementary region $U$ of $\zeta$ such that $U \cap \hat F$ contains $\upsilon_1 \cup \xi$.
If $U\cap \hat F$ contains a simple closed curve intersecting $\upsilon_1$ essentially, we are done by arguing as in the previous paragraph.
Otherwise, we take a simple geodesic arc $\beta$ in $U \cap \hat F$ with endpoints at $\partial \hat F$ intersecting $\upsilon_1$ essentially.
Starting from $\beta$, and using arcs on $\partial \hat F$ disjoint from $\zeta$ and outermost geodesic arcs in parallel families of $\zeta \cap \hat F$ intersecting $\partial \hat F$, as in the previous case when $\zeta_0$ intersects $\upsilon_1$ but is not contained in $\hat F$, we can construct a simple closed curve $\gamma$ in $\hat F$ which intersects $\upsilon_1$ essentially and is disjoint from $\zeta$.
Then we are lead to a contradiction in the same way as before.
Thus, we have reached a contradiction in every case, and  completed the proof.
\end{proof}

%Thus we have shown that if neither $\mathcal{B}$ nor $V$ is  empty,  we reach
%a contradiction since $\{\phi_i|H^k\}$ cannot diverge as shown above, whereas
%$V_k$ intersects $\Lambda^W$ for some $W \in \mathcal{B}$.
Thus we have shown that either $V$ or $\mathcal B$ must be empty.
If $V$ is empty, by Corollary \ref{no V}, we are done.
If $\mathcal{B}$ is empty, then by Proposition \ref{convergence of
function groups}, for every component $V_k$ of $V$, the restriction $\phi_i|H^k$
converges in $AH(H^k)$.
Since $\mathcal{B}$ is empty, the restriction of $\phi_i$ to each
component $W$ of the boundary-irreducible part also converges (after taking a
subsequence and conjugates).
Then,  the same argument as 
the proof of Lemmata 4.5 and 4.6 in \cite{Oh1} shows Theorem
\ref{convergence}.

\section{Unrealisable laminations and Ending laminations}

Having proved Theorem \ref{convergence}, we now know that $\{q(\nu_i)\}$ converges in
$AH(G)$ after passing to  a subsequence.
Let $\phi : G \rightarrow \PSL_2\complexes$
be a  representation with  image $\Gamma=\phi(G)$ such that $(\Gamma, \phi)$ is the limit of (a subsequence of )  $\{q(\nu_i)\}$ in $AH(G)$.
We  consider the 
hyperbolic $3$-manifold $M_\Gamma$.
We use the symbol $\Phi$ to
denote a homotopy equivalence from 
$M_G$ to $M_\Gamma$ induced by the 
isomorphism $\phi$.
Let $C'$ be a relative compact core of $(M_\Gamma)_0$.
If we denote $C' \cap \partial(M_\Gamma)_0$ by $P'$, 
then the pair $(C', P')$ is a pared manifold.

In this section, we shall prove that for the given laminations $\mu_j$ in
Theorem \ref{main}, their images $\Phi(\mu_j)$ actually represent ending laminations of ends of $(M_\Gamma)_0$ (Proposition \ref{ending lamination}).
A similar result on the equivalence of being unrealisable and representing 
an ending
lamination was given independently by Namazi-Souto \cite{NS} as we mentioned in Introduction.

Before stating the main proposition, we shall show that $\Phi$ can be homotoped to take $(C, P \cup T)$ to $(C',P')$ and that $\Phi(\mu_j)$ is unrealisable in $M_\Gamma$, where $P$ and $T$ are as given in Theorem \ref{main}.

\begin{lemma}
\label{parabolic}
The homotopy equivalence $\Phi$ can be homotoped so that $\Phi(C,P \cup T) \subset
(C',P')$ as pairs and $\Phi|(P\cup T)$ is an embedding into $P'$.
\end{lemma}

\begin{proof}
Since $C'$ is a compact core, we can homotope $\Phi$ so that $\Phi(C) \subset C'$.
Since any immersed incompressible torus in $C'$ is homotopic into (a component of ) the union  of 
torus components of $P'$, which we denote by $T'$, we can make  $\Phi(T) \subset T'$.
Let $T_0$ be a component of $T$.
Then, since $\pi_1(T_0)$ is a maximal abelian subgroup in $\pi_1(C)$, we see that $\Phi|T_0$ induces an isomorphism from $\pi_1(T_0)$ to $\pi_1(T'_0)$ for some component $T_0'$ of $T'$,
hence is homotopic to a homeomorphism
to $T'_0$. 
Thus we have shown that $\Phi$ can be homotoped so that $\Phi|T$ is an
embedding into $T'$.

Let $c$ be a core curve of a component of $P$.
By the definition of $\nu_i$, we have $\mathrm{length}_{\nu_i}(c) \rightarrow 0$.
By Theorem 6.2 in Sugawa \cite{Su} or the main theorem of Canary \cite{CaD}, this implies that $\mathrm{length}(\phi_i(c)) \rightarrow
0$, hence also that $\phi(c)$ is parabolic.
Therefore, $\Phi$ can be homotoped so that $\Phi(c) \subset P'$.
Since $c$ represents a generator of a maximal abelian group, (for $(C, P\cup T)$ is a pared manifold,) $\Phi|c$ is
homotopic to a homeomorphism to a core curve of an annulus component of $P'$.
This completes the proof.
\end{proof}

\subsection{Unrealisability of $\mu_j$}
Recall that the aim of this section is to show that $\Phi|\mu_j$ represents an
ending lamination for some end of $(M_\Gamma)_0$.
We shall first see that $\Phi|\mu_j$ cannot be realised by a pleated surface.
Recall that $\mu_j$ is contained in a component $\Sigma_j$ of $\partial C
\setminus (P \cup T)$.

\begin{lemma}
\label{unrealised}
There is no pleated surface homotopic to $\Phi|\Sigma_j$ realising $\mu_j$.
%which
%takes a neighbourhood of the frontier into the cuspidal part.
\end{lemma}

To prove this lemma, we need to invoke the following lemma which appeared as
Lemma 4.10 in \cite{OhM}.
Although we allow $\Gamma$ to have parabolic elements here, the proof in
\cite{OhM} works with only a slight refinement as we shall see below.

In the following, we say that a map $f: \Sigma \rightarrow M_\Gamma$ is adapted to a tied neighbourhood  $N$ of a train track $\tau$ on $\Sigma$  if it maps each branch of $\tau$ to a geodesic segment and each tie of $N$ to a point.

\begin{lemma}
\label{stay near}
Let $\Sigma$ be a component of $\partial C \setminus P$.
Suppose that $\mu$ is an arational measured lamination in $\mathcal{M}(\Sigma)$ 
which can be realised by a pleated surface $f: \Sigma \rightarrow M_\Gamma$.
%sending the frontier of $\Sigma$ to cusps of $M_\Gamma$.
Let $\{w_kc_k\}$ be a sequence of weighted essential simple closed 
curves which converges to $\mu$.
Then for any $\delta >0$ and $t <1$, there exist a continuous map $h 
: \Sigma \rightarrow M_\Gamma$ homotopic to $f$ and 
a subsequence $\{w_{k(l)} c_{k(l)}\}$ of $\{w_kc_k\}$ with the following
two properties.
\begin{enumerate} 
\item The map $h$ is adapted to a tied neighbourhood $N$ of a
train track 
$\tau$ which carries  $\mu$ and the $w_{k(l)}c_{k(l)}$ for 
sufficiently large $l$.
%(A map is said to be adapted to a tied neighbourhood  $N$ of a train track $\tau$ if it maps each branch of $\tau$ to a geodesic segment and each tie of $N$ to a point.)
Moreover, $N$ can be taken to contain $\mu$ and $c_{k(l)}$ in
such a way that their leaves are transverse to the ties of $N$ (without moving
$\mu$ and $c_{k(l)}$ by a homotopy).
\item For sufficiently large $l$,  the simple closed curve $h(c_{k(l)})$
represents a loxodromic class.
Moreover the closed geodesic $c_{k(l)}^*$ in
 $M_\Gamma$ homotopic to  $h(c_{k(l)})$ has a 
part with length at least $t\mathrm{length}h(c_{k(l)})$ which is 
contained in the $\delta$-neighbourhood of the closed curve 
$h(c_{k(l)})$.
\end{enumerate}
\end{lemma}

\begin{proof}
We need to show that the argument in \cite{OhM} works even if we allow
parabolic elements to exist.
It is easy to see that parabolic elements corresponding to punctures of
$\Sigma$ do not affect the argument.

We shall discuss the case when there may be a closed curve not homotopic to a
puncture of $\Sigma$ which represents a parabolic class of $\Gamma$.
(Such an element is called an accidental parabolic element in some
literature.)
We need to show that for every $\{w_kc_k\}$ converging to $\lambda$, after passing
to a subsequence, every $h(c_k)$ represents a loxodromic class of $\Gamma$.
Once we prove this, the argument involving the area estimate of piecewise
geodesic annulus cobounded by $h(c_{k(l)})$ and $c_{k(l)}^*$, which is originally due to
Bonahon \cite{Bo},
works in the same way as in \cite{OhM}.

Now suppose that this is not the case.
Then, by extracting a subsequence, we can assume that all the $h(c_k)$
represent parabolic elements.
(It is impossible that infinitely many $h(c_k)$ are null-homotopic since $\mu$
is contained in the Masur domain.)
Then we can construct a piecewise geodesic ideal annulus $A_k$ cobounded by
$h(c_k)$ and a cusp.
The area of this annulus is equal to the total exterior angle $e_k$ of $h(c_k)$.
We can make each of the exterior angles and $w_k e_k$ arbitrarily small and
the length of the image of each branch of $\tau$ arbitrarily long by
approximating a pleated surface realising $\lambda$ by $h$ closely.
(See Lemma 4.8 of \cite{OhM}.)
For each point of $\partial A_k= h(c_k)$, we consider a geodesic on $A_k$ (with 
respect to the two-dimensional hyperbolic metric induced on $A_k$) perpendicular
to $\partial A_k$.
By the Gauss-Bonnet formula, these geodesic arcs are disjoint and can be
extended indefinitely if they start outside the $\eta$-neighbourhoods of
the vertices, where $\eta$ goes to $0$ as the exterior angle goes to $0$.
Since $w_k\length h(c_k)$ converges to a positive constant, the area of $A_k$ multiplied by $w_k$ goes to infinity as $h$ approximates
the pleated surface realising $\lambda$ closer and closer.
This contradicts the fact shown above that $w_k e_k$, which is equal to $w_k \Area(A_k)$, goes to $0$.
%To show that the argument in \cite{OhM} works even in this case, we need to
%prove that we can take cusp neighbourhoods of $\hyperbolic^3/\Gamma$ which are
%disjoint from all the $c^*_k$ and the $h(\Sigma)$.
%Recall first that since $c^*_k$ is contained in $\mathcal{M}(\Sigma)$ (for
%sufficiently large $k$), there is a pleated surface $f_k$ homotopic to $f$
%realising $c_k$ as $c^*_k$.
%Suppose that $c^*_k$ goes deeper and deeper into a cusp neighbourhood as $k
%\rightarrow \infty$.
%Then, by the same argument as \S6.3 of Bonahon \cite{Bo}, we can bound by a
%number going to $0$
%from above the
%intersection number with a curve on $\Sigma$ representing an accidental
%parabolic element and $c_k/\mathrm{length}_{m_0}(c_k)$ for a fixed hyperbolic
%metric $m_0$ using the fact that there exists a positive $K$ such 
%that $\mathrm{length}(c_k^*) \leq K\mathrm{length}_{m_0}(c_k)$.
%This contradicts our assumption that $c_k$ is arational.
\end{proof}

\begin{proof}[Proof of Lemma \ref{unrealised}]
Suppose that $\mu_j$ is realised by a pleated surface $g_\infty:(\Sigma_j, n_\infty) \linebreak\rightarrow M_\Gamma$ homotopic to
$\Phi|\Sigma_j$. 
 %taking a neighbourhood of the frontier into the cuspidal part.
Then every measured lamination on $\Sigma_j$ with the same support as $\mu_j$ is
also realised by this pleated surface.
Let $\{\gamma_k\}$ be a sequence of simple closed curves on $\Sigma_j$ such that $\{w_k \gamma_k\}$ converges to a measured lamination with the same support as $\mu_j$, and $\gamma_k^*$ the closed geodesic homotopic to $\Phi(\gamma_k)$.
Since $w_k \length_{n_\infty}(\gamma_k)$ converges to $\length_{n_\infty}(\mu_j)$, and $\length(\gamma_k^*) \leq \length_{n_\infty}(\gamma_k)$, we see that $\length(\gamma_k^*)$ grows at most in the order of $w_k^{-1}$ as $k \rightarrow \infty$.

%Then by Lemma \ref{stay near},  for every sequence of simple closed curves
%$\gamma_k$ such that $\{w_k \gamma_k\}$ converges to a measured lamination with
%the same support as $\mu_j$, the closed geodesic $\gamma_k^*$ in
%$\hyperbolic^3/\Gamma$ freely homotopic to $\phi(\gamma_k)$ has length
%growing at least in the same order as $w_k^{-1}$ as $k \rightarrow \infty$.

Let $g_k: (\Sigma_j, n_k) \rightarrow M_\Gamma$ be a pleated surface
homotopic to $\Phi|\Sigma_j$ which realises $\gamma_k$ as $\gamma_k^*$.
We shall show that we can assume that  the $g_k$ intersect cusp neighbourhoods of
$M_\Gamma$ only at a thin neighbourhood of the frontier of $\Sigma_j$, i.e.,
that $g_k(\Sigma_j)$ is disjoint from a sufficiently small cusp neighbourhood corresponding to accidental parabolic
elements for every $k$ if we take a subsequence.
Suppose that  this is not the case.
Then passing to a subsequence, the image of $g_k$ outside a
neighbourhood of the frontier
goes deeper and deeper into a cusp neighbourhood as $k \rightarrow \infty$, for there are only finitely many cusps in $M_\Gamma$.
Since a pleated surface can intersect a small cusp-neighbourhood only at its thin annulus (for the diameters of the thick parts of the pleated surfaces are uniformly bounded),  there is a non-peripheral simple closed curve $d_k$ such that
$\phi(d_k)$  either represents an accidental parabolic element or  is
null-homotopic, and $g_k(d_k)$ is contained in  the  $\epsilon_k$-cusp neighbourhood $U_{\epsilon_k}$  of the same cusp with
$\epsilon_k \searrow 0$ ($\epsilon_k$ goes monotonically to $0$)
as $k
\rightarrow \infty$ after passing to a subsequence.
Furthermore, any point $x$ such that $g_k(x)$ is contained in $U_{\epsilon_k}$ is contained in the $\eta_k$-thin part of $\Sigma_j$ with respect to $n_k$ where $\eta_k \rightarrow 0$ as $k \rightarrow \infty$.

Take positive numbers $v_k$ so that $\{v_k d_k\}$ converges to a non-empty
measured lamination $\xi$ after passing to a subsequence.
Suppose that $i(\xi, \mu_j)= 0$.
Since $\mu_j$ is arational, this implies that the supports of $\xi$ and $\mu_j$
coincide.
%If $\phi(d_k)$ is null-homotopic, $\xi$ is contained in $\overline{\mathcal{C}}(\Sigma_j)$ and contradicts the assumption that $\mu_j$ lies in $\mathcal{M}(\Sigma_j)$.
Since $\xi$, which has the same support as $\mu_j$, is  realised by a
pleated surface homotopic to $\Phi|\Sigma_j$ realising $\mu_j$, by Lemma \ref{stay near}, $\phi(d_k)$ must represent a loxodromic element.
Since $\phi(d_k)$ is either parabolic or trivial, this is a contradiction.
Thus we have $i(\xi, \mu_j)> 0$.
Then $i(d_k, \gamma_k)$ grows in the order of $w_k^{-1}v_k^{-1}$.

Now, let $D_k$ be the distance between $\partial U_{\epsilon_1}$ and $U_{\epsilon_k}$.
Since $\epsilon_k \searrow 0$, we have $D_k \rightarrow \infty$ as $k \rightarrow \infty$.
For each intersection $p$ of $\gamma_k$ with $d_k$, there is an arc  $a_p$ on $\gamma_k$ containing $p$ such that $g_k(a_p)$ starts from  $\partial U_{\epsilon_1}$, goes into $U_{\epsilon_k}$ and comes back to $\partial U_{\epsilon_1}$.
Then  $g_k(a_p)$ has length at least $2D_k$.
(See Figure \ref{fig:deep cusp}.)
Since each intersection corresponds to an arc traversing the $\epsilon_k$-thin annulus around $d_k$, we see that $a_p \cap a_{p'}=\emptyset$ if $p\neq p'$.
It follows that  the geodesic $\gamma_k^*$ has part with length of order $2w_k^{-1}v_k^{-1} D_k$.
%On the other hand, as was seen above, $\length(\gamma_k^*)$ grows at most in the same order as $w_k^{-1}$.
Since $v_kd_k$ converges in $\mathcal{ML}(\Sigma_j)$, we see that $v_k$ is bounded above as $k \rightarrow \infty$.
Therefore $v_k^{-1} D_k$ goes to $\infty$ as $k \rightarrow \infty$, which implies $w_k^{-1}v_k^{-1}D_k$, hence also $\length(\gamma_k^*)$, grows in a higher order than $w_k^{-1}$.
This contradicts the fact that $\length(\gamma_k^*)$ grows at most in the order of $w_k^{-1}$, which was proved above.

\begin{figure}[t]
\includegraphics[height=6cm]{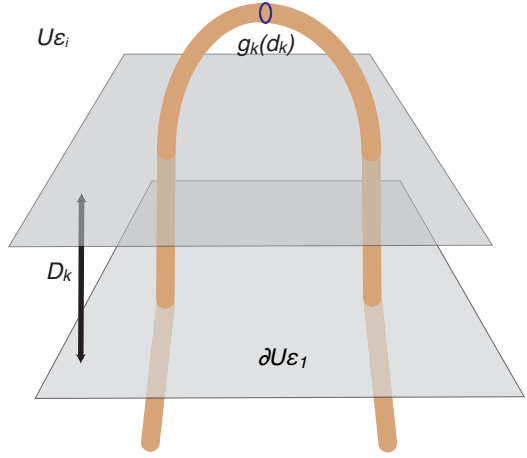}
\caption{Each intersection with $d_k$ contributes $2D_k$ to the length of $\gamma_k^*$.}
\label{fig:deep cusp}
\end{figure}
%implies that for any $\epsilon$, there
%is $\delta >0$ such that
% the realisation $\gamma_k^*$ of $\gamma_k$ by $g_k$ has a part of length
%$\delta
%\mathrm{length}(\gamma_k^*)=\delta \mathrm{length} g_k(\gamma_k)$ which lies in the $\epsilon$-cusp neighbourhood containing the $\phi(d_k)$.
%This again contradicts Lemma \ref{stay near}.

%As was shown in the proof of Lemma \ref{stay near}, we can assume that outside
%a thin neighbourhood of the frontier of $\Sigma_j$, the
%image of $g_k$ is disjoint from  the cusp neighbourhoods of
%$\hyperbolic^3/\Gamma$.

Thus, we have shown that $g_k(\Sigma_j)$ does not intersect cusps outside a
thin neighbourhood of the frontier; hence
there is a uniform upper bound for the diameters of the non-cuspidal parts of the pleated surfaces $g_k$.
Moreover, by Lemma \ref{stay near}, $g_k(\Sigma_j)$ is within a uniformly bounded distance from $g_\infty(\Sigma_j)$.
Then, passing to a subsequence, $\{g_k\}$ converges to a pleated surface realising $\mu_j$ uniformly on every compact set by the compactness of marked pleated surfaces.
(The proof of the compactness is similar to the incompressible case proved in Theorem 5.2.18 of Canary-Epstein-Green \cite{CEG}.
It can be found in the proof of Th\'{e}or\`{e}me 2.3 in Otal \cite{OtT} how to generalise this to the compressible case.)
%Num?ro
In particular, this implies that $\length(\gamma_k^*)$ grows  in the same order as $w_k^{-1}$.
The closed geodesic
$\gamma_k^*$ can be projected down to the geometric limit 
of $\{\phi_i(G)\}$ and then can be pulled back
by an approximate isometry to a $(K_k, \delta_k)$-quasi-geodesic in $M_{\phi_i(G)}$ which is homotopic to $\phi_i(\gamma_k)$ for sufficiently large $i$, with $K_k \rightarrow 1$ and $\delta_k \rightarrow 0$.
Then we can see that the closed geodesic
$\gamma_{i}^{+}$ representing the conjugacy class of $\phi_{i}([\gamma_{i}])$
also
has length growing  in the same order as $w_{i}^{-1}$ as $i \rightarrow
\infty$.
(See also the argument in the proof of Lemma 6.12 in \cite{OhM}.)
 On the other hand, as we showed in the proof of Proposition \ref{convergence
of function groups}, there is a sequence of weighted simple closed curves $r_i
c_i$
converging to a measured lamination $\mu_j'$ having the same support as $\mu_j$
such that $r_i\mathrm{length}(\phi_i([c_i]))$ goes to $0$ as $i
\rightarrow \infty$.
%Take positive real numbers $r_i$ so that $r_i c_i$ converges to a non-empty
%measured lamination $\nu$.
%If $i(\nu, \mu)\neq 0$, then, by the definition of the Thurston
%compactification of the Teichm\"{u}ller space, the length of $r_ic_i$ goes to
%infinity as $i \rightarrow \infty$.
%This contradicts our choice of $c_i$.
%Hence we have $i(\nu, \mu)=0$.
%Since we took $\mu$ to be  arational, this means that $\nu$ has the
%same support as $\mu$.
Since in the argument above, $w_i\gamma_i$ could be taken to be any sequence of weighted simple closed curves converging to a measured lamination with the same support as $\mu_j$, we can let it be $r_i c_i$.
Then we reach a contradiction since $r_i\mathrm{length}(\phi_i([c_i]))$ goes to $0$  whereas it was proved above that $\mathrm{length}(\phi_i([c_i])$ grows in the order of $r_i^{-1}$.
% since
%if the length of $\gamma_{i}$ on the boundary at infinity is bounded, so is the
%length of $\gamma_{i}^{+}$ as
%was shown by Sugawa \cite{Su}
Thus, we have shown   that $\mu_j$ cannot be realised.
\end{proof}

As a consequence of Lemma \ref{unrealised}, the following holds.
(See Lemma 4.4 and Proposition 4.14 in \cite{OhM}, the former of which is based on
the argument of
 Otal
\cite{Otal}.)
\begin{corollary}
\label{pleated surface}
For any sequence of weighted simple closed curves $\{w_k \gamma_k\}$ converging
to a measured lamination with the same support as $\mu_j$, there is a sequence of pleated 
surfaces $f_k$ homotopic to $\Phi|\Sigma_j$  which realise the $\gamma_k$ and tend to an
end of $(M_\Gamma)_0$ as $k \rightarrow \infty$.
Moreover, if $\Sigma_j$ is not null-homologous in $C$  relative to $P$, then  the end
to which the $f_k$ tend is topologically tame and has a neighbourhood
homeomorphic to $\Sigma_j \times \reals$ such that $f_k$ is homotopic  within $\Sigma_j \times \reals$ to
a homeomorphism onto $\Sigma_j
\times \pt$.

 Even in the case when $\Sigma_j$ is null-homologous, the result of Brock-Souto
\cite{BS} implies that $\Gamma$ is topologically tame.
\end{corollary}
In this latter case, we do not know a priori if pleated surfaces $f_k$ realising
$w_k \gamma_k$ give a product structure near the end.
This makes the argument for the latter case more complicated.

\subsection{Main proposition}
In the remaining of this section, we shall show that $\Phi(\mu_j)$ indeed
represents an ending lamination.
\begin{proposition}
\label{ending lamination}
In the setting of Theorem \ref{main},
we can homotope $\Phi$ so that the following holds.
There is a nice compact core $C'$ of $(M_\Gamma)_0$ with $P'=C' \cap \partial (M_\Gamma)_0$ such that
$\Phi|\Sigma_j$ is a homeomorphism to a component of $\partial C' \setminus P'$ and 
 $\Phi(\mu_j)$ represents the ending lamination for the end facing
$\Phi(\Sigma_j)$ for every $j=n+1, \dots , m$.
\end{proposition}

Since the proof of this proposition is rather long, taking up most of the rest of the paper, we shall first summarise our argument before really starting it.

\medskip
\noindent
{\bf Summary of Proof.}
We shall first show that we can reduce the general case to a special case when  the compact core $C$ is a compression body.
(Lemma \ref{compression body}.)
Suppose that $C$ is a compression body.
Then, we have only to consider the cases when $\Sigma_j$ represents a non-trivial second homology class relative to $P$, and when $C$ is a handlebody and $\partial C \setminus P$ is connected.
In the former case, which is easier to deal with, a generalised version of Bonahon's intersection lemma (Lemma \ref{intersection}) implies that $\Phi$ can be homotoped to a homeomorphism for which $\Phi(\mu_j)$ represents the ending lamination.
Most of our discussion will be devoted to the latter case.

In the latter case, we shall first show that $\partial C' \setminus P'$ is also connected.
(Lemma \ref{connected}.)
As a corollary of this lemma, we can show that the convergence of $\{(G_i, \phi_i)\}$ to $(\Gamma, \phi)$ is strong.
Let $\rho_i$ be an approximate isometry between $M_{G_i}$ and $M_\Gamma$ associated to the geometric convergence of $\{G_i\}$ to $\Gamma$.
In the final step, we shall consider a pleated surface $g_i: \partial C \rightarrow M_{G_i}$ homotopic to $\Phi_i$ realising $C^1 \cup \mu_1$.
 (Note that since $\partial C \setminus P$ was assumed to be connected, there is only one $\mu_j$.)
We shall then show that $g_i$ is incompressible in the complement of $\rho_i^{-1}(C')$.
 (Lemma \ref{incompressible}.)
Using another generalised version of Bonahon's intersection lemma (Lemma \ref{uniform intersection}), we can finally show that $\Phi$ can be homotoped so that $\Phi(\mu_1)$ represents the ending lamination.
 \qed
 \medskip

Now, we start the proof.
First we should remark the following.

\begin{lemma}
\label{compression body}
To prove Proposition \ref{ending lamination}, it is sufficient to deal with the case when $C$ is a compression body.
\end{lemma}
\begin{proof}
Let $\mu_j$ be one of  $\mu_{n+1}, \dots , \mu_m$.
Take a nice compact core $(C',P')$ of $(M_\Gamma)_0$, and homotope $\Phi$ so that $\Phi(C)$ is contained in $C'$.
Consider the boundary component of $C$ on which $\mu_j$, hence also $\Sigma_j$ lies, and denote it by $S$.
Let $G^S$ be a subgroup of $G \cong\pi_1(C)$ corresponding to $\iota_\# \pi_1(S)$, where $\iota_\#$ denotes the homomorphism between the fundamental groups induced by the inclusion $\iota : S \rightarrow C$, and set $\Gamma^S$ to be $\phi(G^S)$.
Then $M_{G^S}$ has a nice compact core $\bar C$ to whose boundary $S$ is lifted homeomorphically as $\bar S$.
The core $\bar C$ is homeomorphic to $S  \times I$ if $S $ is incompressible, and is a compression body if $S$ is compressible.
Let $\bar \Sigma_j$ and $\bar \mu_j$ be the lifts of $\Sigma_j$ and $\mu_j$ to $\bar S$.
We denote by $\bar \Phi$ a homotopy equivalence from $M_{G^S}$ to $M_{\Gamma^S}$ induced by $\phi|G^S$.

Now, suppose that Proposition \ref{ending lamination} is valid under the assumption that a compact core of $M_G$ is a compression body.
We shall then apply this for $M_{G^S}$.
Using Lemma \ref{parabolic}, we can assume that $\Phi|P$ is an embedding into $P'$.
Since we also know that Proposition \ref{ending lamination} is true for the case when every boundary component of $\partial C \setminus P$ is incompressible by the main theorem
%num?ro
of \cite{OhL}, this means that there is a relative compact core $(\bar C', \bar P')$ of $(M_{\Gamma^S})_0$ such that $\bar \Phi|\bar \Sigma_j$ is a homeomorphism to a component of $\partial \bar C' \setminus \bar P'$ facing an end $\bar e_j$ for which $\bar \Phi(\bar \mu_j)$ represents the ending lamination.
Now, Canary's covering theorem proved in \cite{CaT} shows that there is an end $e_j$ of $(M_\Gamma)_0$ whose neighbourhood is finitely covered by a neighbourhood of $\bar e_j$. This also implies that there is a component $\Sigma'_j$ of $\partial C' \setminus P'$ such that $\Phi|\Sigma_j$ is homotopic relative to $P$ to a finite-sheeted covering to $\Sigma'_j$ for every $j=n+1, \dots , m$.

For $\Sigma_j$ with $j=1, \dots , n$, we can show that $\Phi|\Sigma_j$ is a homeomorphism to a boundary component by using a result of Abikoff as follows.
Recall that in Definition \ref{n}, the conformal structure $n_i$ on $S$ was taken to converge to $m_j$ if restricted to $\Sigma_j$.
Let $\Omega^S$ be a component of the domain of discontinuity of $G$ such that $\Omega^S/G^S$ is the surface at infinity corresponding to $S$, and $G^{\Sigma_j}$ a subgroup of $G^S$ corresponding to $\iota_\#\pi_1(\Sigma_j)$.
The quasi-conformal deformation $(G_i, \phi_i)$ of $G$ is realised by a quasi-conformal homeomorphism $h_i : S^2_\infty \rightarrow S^2_\infty$.
By Lemma 3 of Abikoff \cite{Ab}, there is an  open subset $\Omega$ of $\Omega^S$ invariant under $G^{\Sigma_j}$ such that $h_i|\Omega$ converges to an equivariant homeomorphism to a component $\Omega'$ of $\Omega_\Gamma$ whose stabiliser is $\phi(G^{\Sigma_j})$.
Since every frontier component of $\Sigma_j$ is mapped by $\Phi$ to a closed curve representing a parabolic element, the surface $\Omega'/\phi(G^{\Sigma_j})$ is of finite type.
Therefore, there is a frontier component of $(M_\Gamma)_0$ facing $\Omega'/\phi(G^{\Sigma_j})$.
Since $C'$ is a nice compact core, there is a component of $\partial C' \setminus P'$ isotopic to this frontier component, which we denote by $\Sigma'_j$.
Since $\Sigma'_j$ carries a group corresponding to $\phi(G^{\Sigma_j})$ in $\pi_1(M_\Gamma) \cong \Gamma$, we can homotope $\Phi|\Sigma_j$ to a homeomorphism to $\Sigma'_j$.

Therefore, $\Phi$ is homotopic to a map whose restriction to every component of $\partial C \setminus P$ is a  finite-sheeted covering.
Since $\Phi|P$ is an embedding, this implies that $\Phi|\partial C$ is a covering into $\partial C'$.
Since $\Phi|C$ is a homotopy equivalence, by Waldhausen's theorem \cite{Wa}, $\Phi|C$ must be homotopic to a homeomorphism to $C'$, and in particular $\Phi|\Sigma_j$ is homotopic to a homeomorphism to $\Sigma'_j$ even for $j=n+1, \dots , m$.
Furthermore, since $\bar \Phi(\bar \mu_j)$ represents the ending lamination of $\bar e_j$, its projection $\Phi(\mu_j)$ represents the ending lamination of the end of $(M_\Gamma)_0$ facing $\Sigma'_j$.
This shows Proposition \ref{ending lamination} for the original $G$ and $\Gamma$.
\end{proof}

{\em We assume from now on until the end of the proof of Proposition \ref{ending lamination} that $C$ is a compression body.}
Recall that we are considering a component $\Sigma_j$ of $\partial C \setminus P$ on which $\mu_j$ lies.
We need to divide the proof into two parts:
\begin{enumerate}[(1)]
\label{two cases}
\item
 the case when $\Sigma_j$ represents a non-trivial second homology class of
$H_2(C,P)$, and 
\item
the case when $\Sigma_j$ is null-homologous relative to $P$.
\end{enumerate}
Since we are assuming that $C$ is a compression body,  (2) corresponds exactly to the case when $C$ is a
handlebody and $
\partial C \setminus P$ is connected.
In both of these cases, we need to use the following Lemma \ref{intersection}, which is a
generalisation of Proposition 3.4 in 
Bonahon \cite{Bo}.
Before stating the lemma, we shall define some condition for  a compact core and a compressible surface, which we need to use in the statement of the lemma.

Let $M$ be a topologically tame hyperbolic $3$-manifold, and $C'$ a nice compact core of $M_0$.
 If $C'$ is boundary-reducible,  Proposition 5.1 in Canary \cite{CaJ} implies that we can isotope $C'$ so that the following holds:
There is a double branched covering $p: \tilde{M} \rightarrow M$ such that
$\tilde{M}$ admits a pinched negatively curved metric with respect to which
 $p$ is an isometry to its image if restricted to the complement of
$p^{-1}(\mathrm{Int} C')$, and such that 
$p^{-1}(C')$ is a boundary-irreducible
compact core of $\tilde{M}$ whose interior contains the tubular neighbourhood of the branching locus where the metric is deformed.

We call such a compact core $C'$ {\sl adequate} with respect to the covering $p$.
We
also say that an embedded surface $F$ in $M$ is {\sl liftable} (with respect to the
covering $p$) when
  $F$
lies outside $\mathrm{Int} C'$ for some adequate $C'$, hence in particular is
lifted isometrically to $\tilde{M}$.

\begin{lemma}
\label{intersection}
Set $M=M_\Gamma$, and consider an adequate core $C'$ of $M$.
Let $M_{0'}$ be the complement of some of the cusp neighbourhoods of $M$.
(This may be equal to $M$ or $M_0$ or something in between.) 
Let $F$ be a properly embedded
compact surface in $M_{0'}$ separating $U \cong F \times \reals$ from
$M_{0'}$.
In the case when $F$ is compressible, we further assume that $F$ is liftable.
We fix a hyperbolic metric $m_F$ on $F$.

Let $\Sigma$ be a component of $\partial C \setminus P$ and $f: \Sigma \rightarrow
M$ a map homotopic to $\Phi|\Sigma$ taking the components of the frontier of $\Sigma$ to cusps of $M$.
Let $\{c_k\}$ be a sequence of simple closed curves contained in the Masur domain
of $\Sigma$
whose projective classes converge to the class represented by a measured
lamination $L_c$
in the 
Masur domain of $\Sigma$, and $\{d_k\}$ either another such sequence or a sequence
of
simple closed curves in the
Masur domain of $F$, such that projective classes $[d_k]$ converge to the
projective
class represented by a measured lamination $L_d$ in  the Masur domain of either $\Sigma$ or $F$.
Suppose that the closed geodesics $c_k^*$ and $d_k^*$ which are respectively homotopic
to either $f(c_k), f(d_k)$ or $f(c_k),
d_k$, depending on the definition of $d_k$, are contained in $U$ for all $k$.
In the case when the $d_k$ lie on $F$, we further assume that a homotopy between
$d_k$ and $d_k^*$ can be taken to be contained in $U \cup F$.
Then, taking subsequences of $\{c_k\}$ and $\{d_k\}$ and denoting them by the same symbols, there exist sequences of simple closed curves $\{C_k\}$ on $\Sigma$ and $\{D_k\}$ on either $\Sigma$ or  $F$ with the following properties.

\begin{enumerate}
\item Let $[\lambda_c], [\lambda_d]$ be  projective laminations to which the projective classes $\{[C_k]\}$ and $\{[D_k]\}$ converge respectively.
Then we have $i(\lambda_c, L_c)=i(\lambda_d, L_d)=0$.
\item
The closed geodesics $C_k^*, D_k^*$ homotopic in $M$ to $f(C_k), f(D_k)$ or $f(C_k), 
D_k$ lie in $U$.
\item
In the case when the $d_k$ lie on $F$, a homotopy between $ D_k^*$ and $D_k$ can be taken to lie in $U \cup F$.
\item
Let $\bar{C}_k$ be a closed curve on $F$ homotopic to $C^*_k$ in $U \cup F$.
We set $\bar{D}_k=D_k$ if the $d_k$ lie on $F$, and define $\bar{D}_k$ to be
a closed curve on $F$ homotopic to $D^*_k$  in $U \cup F$ otherwise.
Then we have an inequality:
$$i(\bar{C}_k, \bar{D}_k) \leq Ke^{-D} \mathrm{length}(\bar{C}_k)
\mathrm{length}(\bar{D}_k),$$
where $D$ is $\min\{d(\bar{C}_k, C^*_k), d(\bar{D}_k, D^*_k)\}$, the
$\mathrm{length}$ denotes the geodesic length with respect to $m_F$,
and $K$ is independent of $k$.
To be more precise, $K$ depends only the pinching
constant for
the branched covering $\tilde{M}$
 and a positive constant bounding from below the
lengths of essential simple closed curves on $F$ with respect to the metric induced from $M$.
\item
If $c_k=d_k$, we can take $C_k, D_k$ to be equal.
\end{enumerate}

\end{lemma}
\begin{proof}
This lemma was shown by Bonahon in \cite{Bo} under the
assumptions that $F$ is
incompressible and that $c_k^*$ and $d_k^*$ can intersect $\epsilon$-Margulis tubes
only at their core curves for some fixed positive constant $\epsilon$ by setting $C_k=c_k, D_k=d_k$. (In this case, the
constant $K$ depends only on a constant bounding the lengths of the essential simple closed curves on $F$ from below. 
See the argument in the proof of Proposition 3.4 in \cite{Bo}.) 

To remove the assumption of incompressibility, still under the assumption that
$c_k^*$ and $d_k^*$  intersect $\epsilon$-Margulis tubes only at their axes,
we apply Canary's construction of a branched covering.
Since we assumed that $F$ is liftable, there is a branched covering  $p:
\tilde{M} \rightarrow M$,
and $F$ is lifted isometrically to an incompressible surface $\tilde{F}$
separating $\tilde{F} \times \reals$ from $\tilde{M}_{0'}$ (the lift of $M_{0'}$ to $\tilde M$), which is contained in a
component of the complement of the boundary-irreducible compact core $p^{-1}(C')$ of
$\tilde{M}$.
Then we can apply Bonahon's argument for the case of incompressible surface.
The only difference is that the metric is not hyperbolic in $\tilde{M}$.
(When we realise homotopies as piecewise totally geodesic ones, they may go out of $U$ and intersect the part where the curvature is not constant.)
Still the argument works since the sectional curvature is pinched between  a
negative constant and $-1$.
This affects the constant $K$, but in such a way that it  only depends on
the pinching constant of $\tilde{M}$.

Now we shall see how to deal with the case when the closed geodesics $c_k^*$ and $d_k^*$ intersect
Margulis tubes outside the core curves.
What we are going to show is  that we can replace $c_k$ with $C_k$ and $d_k$ with $D_k$ in such
a way that neither $C_k^*$ nor $D_k^*$ intersects thin Margulis tubes outside the core curves.
The argument is the same for $c_k$ and $d_k$ except for the case when $d_k$ lies on $F$.
In the latter case, we need to show additionally that the condition (iii) is satisfied.
We shall deal with only $c_k$ from now on until the  last paragraph 
of the proof, and shall explain how the condition (iii) is satisfied in the latter case at the end of the proof.

Since $c_k$ is a simple closed curve in the Masur domain, there is a pleated surface $f_k$ properly homotopic to $f$ which realises $c_k$ as $c^*_k$.
Since we have only to consider the case when $c^*_k$ is far enough from
$\bar{C}_k$ and the diameters of the thick parts of pleated surfaces are uniformly bounded, we can assume that the image of  $f_k$ is entirely contained in $U$.

Suppose that for any small $\epsilon >0$, there is some
 $f_k(S)$ which intersects an $\epsilon$-Margulis tube although $c_k^*$ is not the core curve of the tube.
Then, passing to a subsequence, we can assume that $f_k$
intersects an
 $\epsilon_k$-Margulis tube $T_k$ with core curve different from $c^*_k$ such that
$\epsilon_k \rightarrow 0$.
Let $m_k$ denote the hyperbolic metric on $\Sigma$ induced by $f_k$ from $M$.
Since the diameter of the thick part of $(\Sigma, m_k)$ is uniformly bounded, $f_k(\Sigma)$ can intersect $T_k$ only its thin part.
%By the simple area computation using a polar coordinate, we see that for any
%$A \in \reals$, there exists $k_0 \in \naturals$ such that for $k \geq k_0$,
% if
%$f_k^{-1}(T_k)$ contains a component with an inessential boundary component, then
%it bounds a disc on $S$ with  area with respect to $m_k$ greater than
%$A$. 
%We let such $A$  be greater than $-2\pi \chi(F)=\mathrm{Area}(S)$.
%Then all the components of  $f_k^{-1}(T_k)$ are incompressible on $S$ for $k \geq
%k_0$.
Replacing $\epsilon_k$ if necessary,  with the condition that $\epsilon_k \searrow 0$ preserved, this implies that $\Sigma$ contains an essential simple closed curve $\gamma_k$ with length with respect to $m_k$ going to $0$ such that $f_k(\gamma_k)$ is contained in $T_k$.

Let $L_\gamma$ be a measured lamination to whose projective class $\{[\gamma_k]\}$ converges after passing to a subsequence.
By our definition of $\gamma_k$, the closed curve $f_k(\gamma_k)$ is either null-homotopic or homotopic to an iteration of the core curve of $T_k$.
%First consider the case when $f_k(\gamma_k)$ is homotopic in $T_k$ to (an
%iteration of) the axis of
%$T_k$ for every $k$ after extracting a subsequence. We should note that this is
%always the case when $f$ is incompressible.
If $i(L_\gamma, L_c)=0$, then $L_\gamma$ cannot be contained in $\overline{\mathcal{WC}}(\Sigma)$ since $L_c$ is contained in the Masur domain.
Therefore $f_k(\gamma_k)$ is homotopic to an iteration of the core curve of $T_k$.
By setting $C_k$ to be
$\gamma_k$, this case is reduced to the case when the $c_k^*$  intersect Margulis tubes only at their core curves.
Note also that $C_k^*$ lies on the core curve of $T_k$, hence is homotopic to $f_k(C_k)$ in
$T_k$, which is contained in $U$ for large $k$.
Thus we are done in this case.

Suppose next that $i(L_\gamma,L_c)>0$.
It follows that $c_k$ intersects $\gamma_k$ essentially for large $k$.
Now, we can show the following claim.
\begin{claim}
\label{arcs}
In this situation, we can find a piecewise geodesic simple closed curve $\delta_k$ on $(\Sigma,m_k)$ as follows, after taking a subsequence with respect to $k$.
Here, we assume that $c_k$ and $\gamma_k$ have been isotoped to closed geodesics with respect to $m_k$.
\renewcommand{\labelenumi}{(\alph{enumi})}
\begin{enumerate}
\item $\delta_k$ consists of two arcs, $a_k$ on $c_k$ and $b_k$ on $\gamma_k$.
%,where $c_k$ and $\gamma_k$ are assumed to be isotoped to closed geodesics with respect to $m_k$.
\item The free homotopy class of $\delta_k$ is constant with respect to $k$.
\item The length of $a_k$ goes to infinity whereas that of $b_k$ goes to $0$
as $k \rightarrow \infty$.
\item At the two endpoints, $a_k$ comes to $b_k$ from the opposite sides.
\end{enumerate}
\renewcommand{\labelenumi}{(\roman{enumi})}
\end{claim}
\begin{proof}
The proof of this claim is similar to the argument which can be found in Affirmation 2.3.4
of Otal \cite{OtT}.
 If we fix a hyperbolic metric $m_0$ on $\Sigma$, the simple closed curves $c_k$ and $\gamma_k$, realised as closed geodesics with respect to $m_0$, converge in the Hausdorff topology, after passing to a subsequence, to geodesic laminations $\Lambda_c$ and $\Lambda_\gamma$ which contain the supports of $L_c$ and $L_\gamma$ respectively.
Since we assumed that $i(L_c, L_\gamma)>0$, there is a minimal component $\lambda_c$ of $\Lambda_c$ which is not a simple closed curve and  intersects a minimal component $\lambda_\gamma$ of $\Lambda_\gamma$ transversely.
Now, take an arc $a$ on a leaf of $\lambda_c$ and $b$ on a leaf of $\lambda_\gamma$ such that 
$a$ meets $\lambda_\gamma$ only at its endpoints, $\partial a=\partial b$, and at the endpoints, $a$ comes to $b$ from the opposite sides of the leaf containing $b$.

Since $c_k$ converges to $\Lambda_c \supset \lambda_c$ and $\gamma_k$ converges to $\Lambda_\gamma \supset \lambda_\gamma$, there are arcs $a_k$ on $c_k$ and $b_k$ on  $\gamma_k$ sharing their endpoints, which converge to $a$ and $b$ on $(\Sigma,m_0)$ with respect to the Hausdorff topology.
We set $\delta_k=a_k \cup b_k$.
Realise $c_k$ and $\gamma_k$ as closed geodesics with respect to $m_k$ and we denote the arcs on them corresponding to $a_k$ and $b_k$ by the same symbols.
We shall show that these $a_k$, $b_k$ and $\delta_k$ have the required properties.

The properties (a)  and (d) follow from our definition of $a_k$ and $b_k$.
Since $a_k$ and $b_k$ on $(\Sigma,m_0)$ converge to $a$ and $b$ in the Hausdorff topology, if we take a subsequence, we can assume that the homotopy class of $\delta_k$ is constant, which shows (b).
Since the length of $\gamma_k$ with respect to $m_k$ goes to $0$, the length of $b_k$ goes to $0$, and the arc $a_k$ which traverses a thick annulus around $\gamma_k$ must have length going to $\infty$ with respect to $m_k$.
This completes the proof of our claim.
\end{proof}

Now we return to the proof of Lemma \ref{intersection}.
Recall that $a_k$ is mapped geodesically into $c_k^*$ by $f_k$.
Since the length of $b_k$ goes to $0$, we see that $f_k(\delta_k)$ must represent a non-trivial free homotopy class and the closed geodesic $\delta_k^*$
homotopic to $f_k(\delta_k)$ stays within a distance going to $0$ as $k
\rightarrow \infty$ from $f_k(\delta_k)$, and its length goes to infinity as $k
\rightarrow \infty$.
This is a contradiction since $\delta_k$ represents a constant free homotopy
class.
Thus we have shown that $i(L_\gamma, L_c)=0$ always holds and this completes the proof except for the case when the $d_k$ lie on $F$.
%Suppose next that $\gamma_k$ is not homotopic to an iteration of the axis.
%Then $f_k(\gamma_k)$ must be null-homotopic for every $k$ after passing to a
%subsequence.
%Since $L_c$ is contained in the  Masur domain and $L_\gamma$ is the
%limit of weighted meridians, we have $i(L_c, L_\gamma)>0$.
%Therefore, we can apply the same argument as above for the case when $L_c$ and
%$L_\gamma$ have distinct supports, and
%get a contradiction.

Now finally, we consider the case when the $d_k$ lie on $F$.
We consider the inclusion $\iota : F \rightarrow M_0$ in place of $f$ in the
argument above.
All the argument works without any modification except for the proof of the
condition (iii).
We shall now show that  the condition (iii) holds.
There is a pleated surface $\iota_k : F \rightarrow M_{0'}$ realising $d_k$ which is homotopic to $\iota$.
In general, a homotopy between $\iota$ and $\iota_k$ may not lie in $U \cup F$.
Still, since $F$ is liftable, there is a simple closed curve $D_k$ on $F$ homotopic to   $d_k^*$ in $U \cup F$.
The closed curves $D_k$ and $d_k^*$ are lifted to a simple closed curve $\tilde{d}_k$ on the lift
$\tilde{F}$
of $F$ and a closed geodesic $\tilde{d}_k^{*}$ which is homotopic to $\tilde d_k$.
We can see that for large $k$ there is a pleated surface in
$\tilde{M}_{0'}$ realising $\tilde{d}_k$ as $\tilde{d}_k^{*}$ homotopic to the
inclusion of $\tilde{F}$.
Such a pleated surface is projected to a pleated
surface in $U$ which is homotopic to the inclusion of $F$ in $U \cup F$.
By redefining $\iota_k$ to be a pleated surface obtained as this and apply the argument as above considering $\iota_k$ in place of $f_k$, we get the
condition (iii).
%$U \cup F$
%between
%$\iota_k$ and a homeomorphism $\iota' : F \rightarrow F$ homotopic to $\iota$
%in $M_0$.
%After getting a simple closed curve $\gamma_k$ (such that $\iota_k(\gamma_k)$ is
%homotopic to an iteration of the axis of $T_k$)
%as in the argument above, we let $D_k$ be $\iota'(\gamma_k)$.
%Then the closed geodesic $D_k^*$, which lies on the axis of $T_k$, is
%homotopic to $D_k$ in $U \cup F$.
%This is the condition (iii), which we wanted to prove.
%
%Next we shall see how to generalise the argument to deal with the case when $c^*,
%d^*$ intersect Margulis tubes outside their axes.
%In $\tilde{M}$, any simple closed curve can be realised by a negatively curved
%pleated surface.
%If $c^*$ is far enough from $c$, then a pleated surface $h_c$ realising $c$
%lies in $U$, hence is disjoint from $\tilde{C}$, and is a pleated surface in
%usual sense.
%Therefore the argument above for the incompressible case works.
\end{proof}

\subsection{The case when $\Sigma_j$ is homologically non-trivial.}
%Revoir si on a besoins d'argument plus d?taill?.
Now we reformulate Proposition \ref{ending lamination} in the case (1) in the alternatives in p.\ \pageref{two cases}, and shall prove it.
\begin{lemma}
\label{homologically non-trivial}
Suppose that either  $C$ is not a handlebody or  $\partial C \setminus P$ is
disconnected.
Let $(C', P')$ be a nice compact core of $(M_\Gamma)_0$ which is adequate
with respect to a branched covering.
Then,  $\Phi$ can be homotoped so that  $\Phi|\Sigma_j$ is a homeomorphism 
to a component of $\partial C'\setminus P'$ and $\Phi(\mu_j)$ represents the ending
lamination of the end facing that component for every $j=n+1, \dots , m$.
\end{lemma}

\begin{proof}
By Corollary \ref{pleated surface}, there is a sequence $\{f_k\}$ of pleated surfaces
which realise $w_k \gamma_k$ converging to $\mu_j$ and tend to a
topologically tame end $e$ of $(M_\Gamma)_0$ with a neighbourhood
homeomorphic to $\Sigma_j
\times
\reals$.
Let $\gamma_k^*$ be the closed geodesic homotopic to $\Phi(\gamma_k)$ which is contained
in the image of $f_k$.
By properly isotoping $C'$, we can assume that this
neighbourhood $\Sigma_j \times (0, \infty)$ is a component of the complement of
$C'$ in $(M_\Gamma)_0$ with $\Sigma_j \times \{0\}$ identified with a
component of $\partial C' \setminus P$.
Let $\pi$ be the projection of $\Sigma_j \times [0,\infty)$ to $\Sigma_j \times
\{0\}$.
Since $f_k(\Sigma_j)$ is homotopic to $\Sigma_j \times \pt$ within
$\Sigma_j \times \reals$ by Corollary \ref{pleated surface},
we see that $\pi
\circ f_k$ is homotopic to a homeomorphism, which we shall denote by $h_k$.
Note moreover that $h_k$ can be extended to a map from $C$ to $C'$ inducing an
isomorphism between the fundamental groups which is conjugate to $\phi$.

Now, since $e$ is topologically tame and geometrically infinite (for $f_k$
tends to $e$), there is a measured lamination $\lambda$ on $\Sigma_j \times
\{0\}$ representing the ending lamination for $e$, which is contained in the Masur domain 
$\mathcal{M}(\Sigma_j
\times
\{0\})$.
This means that there is a sequence of weighted simple closed curves $\{s_k
d_k\}$ converging to $\lambda$ in $\mathcal{M}(\Sigma_j \times \{0\})$ such that
there are closed geodesics $d_k^*$ in $\Sigma_j \times (0,\infty)$ tending to $e$ which are
homotopic to the $d_k$ in $\Sigma_j
\times [0, \infty)$.
On the other hand, as was seen above, $h_k(\gamma_k)$ is homotopic to $\gamma_k^*$ in $\Sigma_j \times [0,\infty)$.

Since both $\gamma_k^*$ and $d_k^*$ tend to the end of $(M_\Gamma)_0$ contained in $\Sigma_j
\times [0, \infty)$ as $k \rightarrow \infty$, by Lemma \ref{intersection}, 
there are sequence of simple closed curves $C_k$ on $\Sigma_j$ and $D_k$ on
$\Sigma_j \times \{0\}$ such that $C_k/\mathrm{length}(C_k)$ and $D_k/\mathrm{length}(D_k)$ converge to 
measured laminations $L_C, L_D$ on $\Sigma_j$ and $\Sigma_j \times\{0\}$
respectively with $i(L_C, \mu_j)=0, i(L_D, \lambda)=0$ and 
 $\displaystyle i(\frac{D_k}{\mathrm{length}(D_k)} ,
\frac{h_k(C_k)}{\mathrm{length}(h_k(C_k))})
\rightarrow 0$.
This implies that $\{h_k(C_k)/\mathrm{length}(h_k(C_k))\}$ converges to a
measured
lamination with the same support as $\lambda$ since $\lambda$ is arational, hence
in particular, is contained in $\mathcal{M}(\Sigma_j \times \{0\})$ for large
$k$.
Since $\mu_j$ is also arational and  contained in the Masur domain, $L_C$ has the
same support as $\mu_j$, hence in particular is also contained in the Masur domain.
Therefore $C_k$ lies in the Masur domain of
$\Sigma_j$ for
large $k$.

Since both $\{[C_k]\}$ and $\{[h_k(C_k)]\}$ converge in the projectivised Masur
domains, and
the group of
homeomorphisms of $\Sigma_j \times \{0\}$ which are homotopic to the identity 
in $C'$ acts on the projectivised Masur domain
properly discontinuously, this implies that for sufficiently large $k$, the
homeomorphism $h_k$ does not depend on $k$.
We denote $h_k$ for such large $k$ by $h$.
Since $\lambda$, hence also $L_D$, represents the ending lamination for $e$, so does $h(L_C)$,
hence also $h(\mu_j)$.
Since $\Phi|\Sigma_j$ is homotopic to $h$ in $(M_\Gamma)_0$, this shows that we can homotope $\Phi$
near $\Sigma_j$ so that $\Phi|\Sigma_j$ is
a homeomorphism to $\Sigma_j \times \{0\}$ and $\Phi(\mu_j)$ represents the
ending lamination for the end facing it.
Since we can achieve this only changing $\Phi$ near $\Sigma_j$, we can repeat  the
same operation for each $\Sigma_j$, one by one for all $j=n+1, \dots , m$, and
complete the proof.
\end{proof}

\subsection{Uniform version of Lemma \ref{intersection} for the case (2)}

In Lemma \ref{intersection}, the hyperbolic $3$-manifold in which the closed geodesics
lie is fixed.
We need to consider a sequence of hyperbolic $3$-manifolds and a pair of
sequences of closed
geodesics contained in them, one pair in each manifold, for the proof of the case
(2) in the alternatives in p.\ \pageref{two cases}.
We start with clarifying the setting of our situation where we need to use the lemma.

Suppose that the nice compact core $C$ of $M_G$ is a handlebody and
$\partial C \setminus P$ is connected.
(This is equivalent to the assumption of (2) since we are assuming that $C$ is a compression body.)
Let $C'$ be a nice compact core of $(M_\Gamma)_0$ as before.
Since $\Gamma \cong G$ is a free group, $C'$ is also a handlebody.
{\em We further assume that $
\partial C'
\setminus P'$ is connected.}
(It will be proved in Lemma \ref{connected} that this  always holds in the
case (2).) 
{\em Suppose that $\{\phi_i\}$ converges to $\phi$ strongly, \ie $\{\phi_i(G) \in AH(G)\}$ converges to $\Gamma=\phi(G)$ also geometrically.}
(This will be proved to be the case in Lemma \ref{strong}.)
Set $M_i=M_{\phi_i(G)}$ and $M=M_\Gamma$.
%Let $M_{0'}$ be a
%manifold obtained by deleting some of the cusp neighbourhoods from $M$ as in
%Lemma \ref{intersection}.
Fix a branched covering $p_\infty: \tilde{M} \rightarrow M$ obtained by Canary's construction as was explained just before Lemma \ref{intersection} .
Since $\{G_i\}$ converges strongly to $\Gamma$, we can take branched coverings $p_i :
\tilde{M}_i \rightarrow M_i$  so that
$\tilde{M}_i$ converges geometrically to $\tilde{M}$ and approximate
isometries are equivariant with respect to the covering translations.
(Refer to Chapter 5  of \cite{OhM} for a proof of this fact.)
We assume that $(C',P')$ is a nice compact core of $M$ which is adequate with respect
to $p_\infty$.
Let $F_\infty$ be $\partial C' \setminus P'$ for this adequate $C'$ separating
$U\cong F_\infty \times \reals$ from $M_0$.
We denote  by $\rho_i$ an approximate isometry
between $M_i$ and $M$, which we can assume to be $C^\infty$, and set $F_i=\rho_i^{-1}(F_\infty
\cup P')=\rho_i^{-1}(\partial C')$.
We fix a  hyperbolic metric $m_{\partial C'}$ on $F_\infty$, and let $m_{F_i}$ be its pull-back $\rho_i^*(m_{\partial C'})$.
Suppose that $F_i$ is an embedded surface in $(M_i)_0$ which is liftable with
respect to $p_i$ and separates $U_i \cong F_i
\times
\reals$ from  $(M_i)_0$.
(We shall show that we can take $C'$ so that this is the case in Lemma
\ref{adequate}.)

Let $\Sigma$ be  $\partial C \setminus P$, which we assumed to be connected.
Also, let $f_i : \Sigma \rightarrow M_i$ be a map taking each component of the
frontier to a closed geodesic
such that $\rho_i \circ f_i$ converges uniformly on every compact set in
$\mathrm{Int} \Sigma$ to
$f: \Sigma
\rightarrow
M$ homotopic to $\Phi|\Sigma$ taking each component of the frontier
to a cusp.
%D?finir la convergence forte pour les surfaces pliss?es.
\begin{lemma}
\label{uniform intersection}
In this situation, suppose that $\{c_k\}$ is a sequence of simple closed curves on $\Sigma$ such that $\{r_k c_k\}$ converges 
in the Masur domain of $\Sigma$  to a measured lamination $L_c$ and that
$\{d_k\}$ is one on either $\Sigma$ or $F_\infty$ such that $\{s_k d_k\}$ converges  in the Masur domain of either $\Sigma$ or $F_\infty$ to a measured lamination $L_d$.
Suppose that the closed geodesics $c_i^*$ and $d_i^*$ which are homotopic to
$f_i(c_i)$ and either $f_i(d_i)$ or
$\rho_i^{-1}(d_i)$ respectively are contained in $U_i$ for each $i$, and that
both
$d_{M_i}(c_i^*, F_i)$ and $d_{M_i}(d_i^*, F_i)$ go to $\infty$ as $i \rightarrow \infty$.
Furthermore, we assume that $d_i^*$ is homotopic to $\rho_i^{-1}(d_i)$ in $U_i
\cup F_i$ in the case when the $d_k$ lie on $F_\infty$.
Then, after taking a subsequence of $\{\phi_i\}$, there exist sequences of simple closed curves $\{C_i\}$ on $\Sigma$ and
$\{D_i\}$ on either $\Sigma$ or  $F_\infty$ depending on where the $d_k$ lie, with the following properties.

%%V?rifier les conditions suivantes.
\begin{enumerate}
\item In the projective lamination space of $\Sigma$  and that of either $\Sigma$ or $F_\infty$
respectively, $\{[C_i]\}$ converges to a projective lamination $[\mu_c]$ and $\{[D_i]\}$ converges to $[\mu_d]$ such that $i(\mu_c, L_c)=0$ and $i(\mu_d, L_d)=0$.
\item
The closed geodesics $C_i^*, D_i^*$, which are homotopic in $M_i$ to $f_i(C_i)$ and either
$f_i(D_i)$ or $\rho_i^{-1}(D_i)$ respectively, lie in $U_i$.
\item The closed geodesic $ D_i^*$ is homotopic to $\rho_i^{-1}(D_i)$
in $U_i \cup F_i$ in the case when the $d_k$ lie on $F_\infty$.
\item
Let $\bar{C}_i, \bar{D}_i$ be the closed curves on $F_i$ homotopic to $C^*_i,
D_i^*$
in $U_i
\cup F_i$.
(This implies that $\bar{D}_i=D_i$ in the case when the $d_k$ lie on
$F_\infty$.)
Then $\bar{C}_i, \bar{D}_i$ are disjoint from $\rho_i^{-1}(P')$, and we have 
$$\frac{i(\bar{C}_i, \bar{D}_i) }{\mathrm{length}(\bar{C}_i)
\mathrm{length}(\bar{D}_i)} \rightarrow 0,$$
where $\length$ denotes the geodesic length with respect to $m_{F_i}$.
\item
If $c_k=d_k$, then we can take $C_k$ and $D_k$ to be equal.
\end{enumerate}
\end{lemma}

\begin{proof}
First we consider the case when  $c_i^*, d_i^*$
do not intersect $\epsilon$-Margulis tubes outside their core curves.
The infimum of the lengths of essential closed
curves on $F_i$ with respect to the metric induced from $M_i$ converge to that of $F_\infty$ with respect to the metric induced from $M$,
hence a lower bound of the lengths can be taken to be independent of $i$. 
On the other hand, since $\tilde{M}_i$ converges to $\tilde{M}$ geometrically, and since $\tilde{M}_i$
has a hyperbolic metric outside a compact set which converges to a compact set
of $\tilde{M}$ as $i \rightarrow \infty$, there is a negative constant
uniformly bounding
the sectional curvatures of the $\tilde{M}_i$ below.
%Also there is a universal positive lower bound for the lengths of essential
%simple closed curves on the $F_i$ since they converge to $F$ geometrically.
Since the constant $K$ in Lemma \ref{intersection} depends only on these two, we see
that we get the conclusion by setting $C_i=c_i, D_i=d_i$  in this case.

Now we shall consider the case when for any small $\epsilon >0$, there exists $i$
such that the
closed geodesic
$c_i^*$ or $d_i^*$ intersects an $\epsilon$-Margulis tube of $M_i$ outside the
core curve.
As in the proof of Lemma \ref{intersection}, we have only to consider $c_i^*$.
By extracting a subsequence with regard to $i$, we can assume that there
exists $\epsilon_i \rightarrow 0$  such that, the closed 
geodesic $c_i^*$, which lies in $M_i$, intersects
an $\epsilon_i$-Margulis tube $T_i$ in $M_i$ outside the core curve for each $i$.

%Then, taking a subsequence, we can assume that $k(i)$ is constant, which we
%denote by $k_0$.
%If $\phi([c_{k_0}])$ represents a parabolic element, the length of $c_{k_0,i}^*$
%goes to $0$.
%Hence for sufficiently large $i$, the closed geodesic $c_{k_0,
%i}^*$ is contained in the Margulis tube as its axis, and cannot intersect
%other Margulis tubes.
%This contradicts our assumption.
%Therefore, $\phi([c_{k_0}])$ must represent a loxodromic element.
% Since $\phi_i(G)$ converges strongly, the closed geodesic $c_{k_0,i}^*$
%converges geometrically to the closed geodesic $c_{k_0,\infty}^*$ representing
%$\phi([c_{k_0}])$.
%Since $c_{k_0,i}^*$ intersects $\epsilon_i$-Margulis tube with $\epsilon_i
%\rightarrow 0$, this is impossible.
%Thus we have shown $k \rightarrow \infty$ as $i \rightarrow \infty$; hence we can
%assume that $k=i$.

Let $h_i: (\Sigma, \sigma_i) \rightarrow M_i$ be a pleated surface with boundary homotopic to
$f_i$ relative to $\mathrm{Fr} \Sigma$ which realises $c_i$ as $c_i^*$.
Let $s$ be a component of the frontier of $\Sigma$.
Then $h_i(s)$ is a closed geodesic whose length goes to $0$ as $i \rightarrow
\infty$ since $\phi(s)$ represents a parabolic class whereas $\phi_i(s)$ does not.
Therefore, $h_i(s)$ is the core curve of some Margulis tube for large $i$.
This means that each  component of the frontier of $\Sigma$ is  either disjoint from or contained in $h_i^{-1}(T_i)$.
By a simple computation of area using a polar coordinate, we see that for any
$A \in \reals$, there exists $i_0 \in \naturals$ such that for $i \geq i_0$,
if $h_i^{-1}(T_i)$ contains a component with inessential boundary components, then
one of the inessential boundary components bounds a disc on $\Sigma$ whose  area with respect to $\sigma_i$ is greater than $A$. 
We let such $A$  be greater than $\mathrm{Area}(\Sigma)$, which does not depend on $i$ by the Gauss-Bonnet formula.
Then all the components of  $h_i^{-1}(T_i)$ are incompressible on $\Sigma$ for $i \geq
i_0$, and we see in particular
that $h_i^{-1}(T_i)$ contains a simple closed curve $\gamma_i$ which is
essential on $\Sigma$.
(We regard even peripheral curves as essential here.)

%We first consider the case when $h_i(\gamma_i)$ is essential in $M_i$, that
%is, when $h_i(\gamma_i)$ is homotopic to an iteration of the axis of $T_i$.
Passing to a subsequence, we can assume that $[\gamma_i]$ converges to a
projective lamination $[L_\gamma]$ in $\PML(\Sigma)$.
If $i(L_c, L_\gamma)=0$, as in Lemma \ref{intersection}, $h_i(\gamma_i)$ is essential in $M_i$ for large $i$.
By defining $C_i$ to be
$\gamma_i$, we can apply the argument for the case when the closed geodesics
do not intersect Margulis tubes outside the core curves and get the inequality as we
wanted.
On the other hand, for any core curve $\delta$ of a component of $P'$, the
geodesic length of $\rho_i^{-1}(\delta)$ in $M_i$ goes to $0$.
Therefore by the same inequality, we see that the projection of $f_i(C_i)$ to
$F_i$ in $F_i \cup U_i$ is disjoint from $\rho_i^{-1}(P')$ for large $i$.

Next assume that $i(L_c, L_\gamma)>0$.
Then applying Claim \ref{arcs} in the proof of Lemma \ref{intersection}, we get
a piecewise geodesic simple closed curve $\delta_i$ on $(\Sigma,\sigma_i)$ representing a constant free homotopy class
with regard to $i$ such that the closed geodesic $\delta_i^*$ homotopic to $h_i(\delta_i)$ lies within a
distance going to $0$ from $h_i(\delta_i)$ and has length going to infinity.
This contradicts the fact that $G_i$ converges algebraically.

In the case when $d_k$ lies on $F_\infty$, we define $D_i$ to be the
projection of $D_i^*$ to $F_i$ in $U_i \cup F_i$.
Then $D_i$ converges to a
projective lamination $[\mu_d]$ with $i(L_d, \mu_d)=0$ since $d_i$ is
homotopic to $d_i^*$ in $U_i \cup F_i$, and we can use the same argument as in
the proof of Lemma \ref{intersection} taking into account the fact that
$\tilde{M}_i$ converges to $\tilde{M}$ geometrically.
%V?rifier si cet argument est valide.
\end{proof}

\subsection{Connectedness of $\partial C' \setminus P'$ and strong convergence}
Now we return to the general situation of the case (2) when $C$ is a handlebody and $\partial C \setminus
P$ is connected.
Since in this case, $n=0$ and $m=1$ in the statement of Theorem \ref{main}, we
denote the only one $\mu_j$, which is $\mu_1$, by $\mu$, and $\Sigma_j$ by
$\Sigma$.
Also since $\Gamma$ is free and topologically tame, $M_\Gamma$ is
homeomorphic to an open handlebody and $\Phi$ is homotopic to a homeomorphism.
{\em Therefore, we can assume that $\Phi$ is a homeomorphism.}

\begin{lemma}
\label{connected}
Let $(C',P')$ be a nice compact core of $(M_\Gamma)_0$ which is adequate with respect to
a branched covering.
Then, $\partial C' \setminus P'$ is also connected.
\end{lemma}

\begin{proof}
Recall that, by Lemma \ref{unrealised}, we have a sequence of pleated surfaces $f_k: \Sigma
\rightarrow (M_\Gamma)_0$ homotopic to $\Phi|\Sigma$ relative to $P$
tending to an end $e$ of $(M_\Gamma)_0$, which realise $\{w_k \gamma_k\}$
converging to $\mu$ as $k \rightarrow \infty$.
Let $\Sigma'$ be a component of $\partial C' \setminus  P'$ facing this end $e$.
Since $C'$ is nice, the component of the complement of $C'$ in $(M_\Gamma)_0$ facing $\Sigma'$
is homeomorphic to $\Sigma' \times (0,\infty)$ with $\Sigma' \times \{0\}$ identified with $\Sigma'$, and  the image
$f_k(\Sigma)$ is contained in $\Sigma' \times
(0,\infty)$ for every $k$ after passing to a subsequence.
This implies that $f_k(\Sigma)$ is homotoped  within $\Sigma'
\times (0,\infty)$ into $\Sigma' \times \pt$.
By a standard argument as in \S4.E of \cite{OhM}, in which  a subsequence of  $\{f_k\}$ is extended to a family of  pleated
surfaces (allowed to be not constantly curved but negatively curved) realising an half-open arc tending to $\mu$, we can show that the $f_k$ are homotopic in $\Sigma' \times (0,\infty)$ after passing to a subsequence. 
Since $(f_k)_\#\pi_1(\Sigma)$ carries the entire $\pi_1(C') \cong \Gamma$, this
implies that $\pi_1(\Sigma')$ also carries the entire $\pi_1(C')$.
This means that $(C',P'')$ is either a relative compression body having
$\Sigma'$ as its exterior boundary or a product $I$-bundle as a pair, where $P''$ is the union of the components of $P'$ intersecting $\Fr \Sigma'$.

If $(C',P'')$ is a relative compression body with empty interior boundary, then
$\Sigma' =
\partial C'
\setminus P'$ and we are done.
Therefore, we can assume that $(C',P'')$ is either a relative compression body
with non-empty interior boundary or homeomorphic to $\Sigma' \times [-1,0]$.
In either case, $\partial C' \setminus P'$ has a component with negative Euler
characteristic other than $\Sigma'$.
On the other hand, we can show the following lemma which contradicts this fact.
Once we prove this lemma, we shall reach a contradiction and the proof of Lemma \ref{connected} will be completed.

\begin{lemma}
\label{Euler}
If $\partial C' \setminus P'$ is not connected, then $\partial C' \setminus P'$ cannot have a component with negative Euler characteristic other than $\Sigma'$.
\end{lemma}
\begin{proof}[Proof of Lemma \ref{Euler}]
Let $\pi: \Sigma' \times [0,\infty) \rightarrow \Sigma'$ be the projection to the first
factor.
Since the $f_k$ are all homotopic in $\Sigma' \times (0,\infty)$, their projections $\pi \circ f_k$ are all homotopic.
Note that $f_k$ can be extended to a map from $\partial C$ to $\partial C'$ by taking each component of $P$ to an annulus on $\mathrm{Fr} (M_\Gamma)_0$.
Therefore, we can take a map $p: \partial C \rightarrow \Sigma'$ with $p(P) \subset P''$ whose restriction to $\Sigma$ is  homotopic to all the $\pi \circ f_k$ relative to $P'$.
Furthermore this map can be extended to a map $\bar p$ from $C$ to $C'$ since
both of them are handlebodies and $\Phi$ is a homotopy equivalence.

Suppose first  that $p|\Sigma$ is not a degree-$0$
map to $\Sigma'$ relative to $P''$.
Then $f_k(\Sigma)$ represents a non-trivial second homology class in $\Sigma'
\times \reals$ relative to $\partial (M_\Gamma)_0$.
Recall that $\{f_k\}$ extends to a continuous family of (negatively curved) pleated surfaces $\rho: \Sigma \times [0,\infty) \rightarrow \Sigma' \times [0,\infty)$ such that $\rho(\ , k) =f_k$ and $\rho|\Sigma \times \{t\}$ tends to the end $e$ as $t \rightarrow \infty$.
Then the argument in Claim 3 of \S4.E in \cite{OhM} shows that $\rho$ can be deformed so that $\rho|\Sigma \times (t_0, \infty)$ is a homeomorphism onto a neighbourhood of $e$, and in particular that $f_k$ is homotopic to a homeomorphism to $\Sigma' \times \{t\}$ in $\Sigma' \times (0,\infty)$ for large $k$.
By extending the homeomorphism to $P$, we see that $\Sigma' \cup P''$ must be a closed surface homeomorphic to $\partial C$ then.
%It follows that the same argument as that for the case (1) works, i.e., Claim
%\ref{homologically
%non-trivial} is valid and $f_k|\Sigma$ is homotopic in $\Sigma' \times \reals$ to
%a homeomorphism to $\Sigma'$.
This implies in particular that $\Sigma'$ is the only component of $\partial
C' \setminus P'$ contradicting our assumption, and we are done.

Therefore, we can assume that $p|\Sigma$ has degree $0$ relative to $P''$ from now on.
%Recall that by Lemma \ref{parabolic}, $\Phi|P$ is an embedding into $P'$.
Since $p(P) \subset P''$, the map $p$ itself is also a degree-$0$ map from $\partial C$ to $\Sigma \cup P''$.
%We shall apply the simple loop conjecture proved by Gabai \cite{Ga}.
%Note that the simple loop conjecture for a surface with boundary is false in
%general.
%Still, we can show that this holds in the present case where the boundary is
%mapped homeomorphically,
%by using the relative version, i.e., simple arc conjecture (Theorem 3.1  of
%\cite{Ga}).
Then the simple loop conjecture proved by Gabai \cite{Ga} gives us  an essential  simple closed curve $c$ on $\partial C$ such that $p(c)$ inessential on $\Sigma' \cup P''$.
Since  $(C', P'')$ is either a relative compression body with exterior boundary $\Sigma'$ or  a trivial $I$-bundle over $\Sigma'$, the union of annuli $P''$ is just a collar neighbourhood of $\Fr \Sigma'$.
Therefore the degree-$0$ map $p$ from $\partial C$ to $\Sigma' \cup P''$ can be regarded as a map to $\Sigma'$.
Since $\Phi$ is a homotopy equivalence, in this situation, $c$ must bound a compressing disc in $C$.
Performing compression of $\partial C$ along this disc, we get a surface $\Sigma_1$ which contains the complement of a thin neighbourhood of $c$ in $\Sigma$ and is disconnected if the compressing disc is separating, and a map $p':\Sigma_1 \rightarrow \Sigma'$ whose restriction to $\Sigma \cap \Sigma_1$ coincides with $p|\Sigma \cap \Sigma_1$.
Let $P_1$ be the union of annulus components among $\Sigma_1 \cap P$.
Suppose that there is a component $\bar \Sigma_1$ of $\Sigma_1$ such that $p'|\bar \Sigma_1$ has non-zero degree as a map from $(\bar \Sigma_1, P_1 \cap \bar \Sigma_1)$ to $(\Sigma', P'')$.
Then we can also compress $\rho$ to get a  map $\rho_1 : \bar \Sigma_1 \times [0,\infty) \rightarrow \Sigma' \times [0,\infty)$ such that $\rho(\ , t)$ tends to $e$ as $t \rightarrow \infty$.
Repeating the previous argument for the non-zero degree case, we see that $\rho_1((\bar \Sigma_1 \setminus \Int P_1) \times (t_0, \infty))$ is a homeomorphism to a neighbourhood of $e$.
This implies that $p'_\#(\pi_1(\bar \Sigma_1))$ carries the fundamental group of $\Sigma'$ hence the entire $\pi_1(C')$.
This is a contradiction since $\bar \Sigma_1$ cannot carry the entire fundamental group of $C$.
Thus we have shown that $p'$ has degree $0$ on each component of $\Sigma_1$.

Therefore, using this argument repeatedly and extending the map to compressing discs, we can extend $p$ to a map from $C$ to $\Sigma'$.
Let $\iota$ be the inclusion of $\partial C$ to $C$.
What has been shown above  implies that $\mathrm{Ker} \iota_\# \subset
\mathrm{Ker} p_\#$.
Since $\Phi_\#$ is an isomorphism, we see that $\mathrm{Ker} p_\# \subset \mathrm{Ker} \iota_\#$, hence $\mathrm{Ker} \iota_\#=\mathrm{Ker}p_\#$.

Now fix a hyperbolic metric on $\Sigma'$ which makes each component of the
frontier an open end.
We need to use the following fact.
\begin{claim}
\label{continuity}
We can define a continuous map $p_*$ from $\mathcal M(\Sigma')$ to the space of the geodesic currents on $\Sigma'$ such that 
\begin{enumerate}
\item
when $\mu$ is a simple closed geodesic in $\mathcal M(\Sigma')$, its image $p_*(\mu)$ is the geodesic current corresponding to the closed geodesic homotopic to $p(\mu)$, and
\item when $p(\mu)$ is homotopic to a measured lamination on $\Sigma'$, which we denote again by $p(\mu)$, the geodesic current $p_*(\mu)$ coincides with $p(\mu)$ regarded as a geodesic current.
\end{enumerate}
\end{claim}
\begin{proof}
Consider the universal covering $\rho : \hyperbolic^2\rightarrow \Sigma'$ and let
$\Lambda$ be the limit set in $S^1_\infty$ of the covering translations by $\pi_1(\Sigma')$.
Recall that a geodesic current on $\Sigma'$ is defined to be a measure
on $\mathcal{G}=(\Lambda\times \Lambda \setminus\Delta)/\integers_2$ which is invariant under $\pi_1(\Sigma')$, where $\Delta$
denotes
the diagonal and $\integers_2$ acts as the interchange of the two factors.
(Refer to Bonahon \cite{Bo} and \cite{BoI}.)

Let  $\mu$ be a measured lamination contained in the Masur domain of $\Sigma$.
By Lemma \ref{extension}, $\mu \cup C^1$ lies in $\mathcal D(C)$.
Let $\tilde{\partial C}$ be the covering of $\partial C$ associated to the kernel of
the homomorphism $\iota_\# : \pi_1(\partial C) \rightarrow \pi_1(C)$ induced from the
inclusion $\iota$.
%Then each leaf $l$ of $\mu$ is mapped to a leaf $p(l)$ which is lifted to an
%open arc ending at two distinct points of $\Lambda$.
Note that the map $p$ is covered by a map between the coverings
$\tilde{p}: \tilde{\partial C} \rightarrow \hyperbolic^2$ since $\mathrm{Ker}
\iota_\#$ is contained in $\mathrm{Ker} p_\#$.

%Since $\pi_1(C)$ is isomorphic to $\pi_1(C')$ by $\Phi_\#$, any element of
%$\pi_1(\Sigma)$ that is not in $\mathrm{Ker} \iota_\#$ is mapped to a
%non-trivial element by $p_\#$.
Moreover, since $\mathrm{Ker} \iota_\#=\mathrm{Ker}p_\#$, the map
$\tilde{p}$ is proper and takes the ends of $\tilde{\partial C}$ to $\Lambda$
injectively.
Since $\mu \cup C^1$ lies in $\mathcal D(C)$, by Lemmata 3.1 and 3.3 in Lecuire \cite{LecI}, each lift of a leaf of $\mu$ connects two distinct points of the ends.
It follows that the map $\tilde{p}$ takes a lift of each leaf of $\mu$ to an open
arc  in $\hyperbolic^2$ ending
at two distinct points of $\Lambda$.
The ends of $\tilde{\partial C}$ can be regarded as embedded on the Riemann sphere
as the limit set $\Lambda_0$
of a Schottky group.
(See p.11 of Otal \cite{OtT} and \S 4 of Kleineidam-Souto \cite{KS}.)
Then the observation above shows that $\tilde{p}$
induces a continuous embedding from $\Lambda_0$, identified with the set of the ends of $\tilde{\partial C}$, to $\Lambda$.
% since the ends of $\tilde{\partial C}$ are embedded into $\Lambda$.
The transverse measure of $\mu$ defines a measure on $(\Lambda_0 \times \Lambda_0 \setminus \Delta_0)/\integers_2$ invariant under the covering translation group, where $\Delta_0$ denotes the diagonal of $\Lambda_0 \times \Lambda_0$.
We push forward this measure to $(\Lambda \times \Lambda \setminus \Delta)/\integers_2$ and takes the sum of all of its distinct translates by $\pi_1(\Sigma')$.
This geodesic current is defined to be $p_*(\mu)$.
%The transverse measure of the preimage of $m$ in $\hyperbolic^2$ obtained by lifting that of
%$\mu$ induces a measure $\mu$ on $(\Lambda_0 \times \Lambda_0\setminus
%\Delta)/\integers_2$.
%Summing up the measure obtained by pushing forward $m$ and its translates by $\pi_1(\Sigma')/\mathrm{Im} p_\#$, we get a measure on $(\Lambda \times \Lambda \setminus \Delta)/\integers_2$ which is invariant under $\pi_1(\Sigma')$.
%We define this measure to be $p_*(\mu)$.
We can easily check this coincides with the usual definition of geodesic currents corresponding to closed curves when $\mu$ is a simple closed curve in $\mathcal{M}(\Sigma)$.
The continuity of $p_*$ then  follows from the fact that $\Lambda_0$ is embedded into $\Lambda$.
It is obvious from our construction that in the case when $p(\mu)$ is a measured lamination on $\Sigma'$, if we regard $p(\mu)$ as a geodesic current, then it coincides with our geodesic current $p_*(\mu)$.
\end{proof}

Now, we return to the proof of Lemma \ref{Euler}.
Let $\lambda$ be a measured  lamination representing the ending lamination of the end $e$ facing $\Sigma'$.
Let $\{s_k d_k\}$ be a sequence of weighted simple closed curves converging to $\lambda$
such that the closed geodesic $d_k^*$ homotopic to $d_k$ in $\Sigma' \times [0,\infty)$
tends to the end $e$ as $k \rightarrow \infty$.

%Montrer que cette projection est homotope ? un plongement, utilisant le lemma%%
%d'intersection.
Regard  $w_kp(\gamma_k)$ as
a weighted closed geodesic with respect to the hyperbolic metric as we fixed before Claim \ref{continuity}.
Then $w_kp(\gamma_k)$ converges to some geodesic current on $\Sigma'$ by Claim \ref{continuity}.
Since the closed geodesic $f_k (\gamma_k)$, which is homotopic to $p(\gamma_k)$
in $\Sigma'
\times [0,\infty)$, 
tends to the end $e$ in $\Sigma'
\times [0,\infty)$, by Lemma
\ref{intersection},
there are sequences of simple closed curves $\{C_k\}$ on $\Sigma$ and
$\{D_k\}$ on $\Sigma'$ with conditions in the statement of Lemma
\ref{intersection} such that
 $$\lim_{k\rightarrow \infty} i(\frac{p_*(C_k)}{\mathrm{length}(p_*(C_k))}, \frac{p_*(C_k)}{\mathrm{length}(p_*(C_k)
)})= 0$$ and
 $$\lim_{k \rightarrow \infty}  i(\frac{p_*(C_k)}{\mathrm{length}(p_*(C_k))}, \frac{D_k}{\mathrm{length}(D_k)})
= 0.$$

As in Lemma \ref{intersection}, $C_k/\mathrm{length}(C_k)$ converges to a
measured lamination $\mu'$ having null-intersection number with $\mu$.
Since $\mu$ is arational, this means that $\mu'$ has the same support as
$\mu$.
Since $\mu'$ is also contained in the Masur domain, by Claim \ref{continuity},
$\{p_*(C_k)/\mathrm{length}(C_k)\}$ converges to the geodesic current $p_*(\mu')$.
We should also note that
there is a constant $L$ such that $\mathrm{length}(p_*(C_k)) \leq
\mathrm{length}(p(C_k)) \leq L \mathrm{length}(C_k)$, where $C_k$ is assumed
to be a closed geodesic.
Note that in contrast, $p(C_k)$ may not be a closed geodesic and $\mathrm{length}(p(C_k))$ is not the geodesic length.
Therefore we have \begin{eqnarray*}i(p_*(\mu'), p_*(\mu'))= \lim_{k \rightarrow \infty} i(\frac{p_*(C_k)}{\mathrm{length}(C_k)},
\frac{p_*(C_k)}{\mathrm{length}(C_k)}) \\\leq L^2\lim_{k\rightarrow \infty} i(\frac{p_*(C_k)}{\mathrm{length}(p_*(C_k))}, \frac{p_*(C_k)}{\mathrm{length}(p_*(C_k)
)})= 0.\end{eqnarray*}
This means that $p_*(\mu')$ is a measured lamination by Proposition 17 in Bonahon \cite{BoI}.
Since $\mu$ has the same support as $\mu'$, this implies that $p_*(\mu)$ is
also a measured lamination with the same support as $p_*(\mu')$ by our definition of $p_*$.
On the other hand, $\{D_k/\mathrm{length}(D_k)\}$ converges to a
measured lamination $\lambda'$ such that $i(\lambda, \lambda')=0$.
Since $\lambda$ is arational, this means that $\lambda'$ has the same support
as $\lambda$.
The equation above implies that $i(p_*(\mu'), \lambda')=0$, hence
$i(p_*(\mu),\lambda)=0$.
Since $\lambda$ is arational, it follows that $p_*(\mu)$ has the same support
as $\lambda$.

%Recall that $f_k$ is extended to a map from $C$ which is  homotopic to a
%homeomorphism $h: C \rightarrow C'$ since both $C$ and $C'$ are handlebodies.
%What was proved above means  that the support of $h(\mu)$ is homotopic to that of
% $\lambda$ in $C'$.
%Since $\lambda$ is arational (for it represents an ending lamination), $\pi
%\circ  \iota \circ p(\mu)$ is homotopic to an arational lamination.
%This is possible only when $\iota (\Sigma'')$ is homotopic to a homeomorphism
%to $\Sigma'$.
%By considering the Euler characteristic as before, this implies that $(C',P')
%\cong
%(\Sigma' , \partial \Sigma') \times I$.
%Since both $C$ and $C'$ are handlebodies and $\Phi|C$ is homotopic to a
%homeomorphism to $C'$, there is a homeomorphism $f: C \rightarrow C'$
%homotopic to $f_k$ in $\hyperbolic^3/\Gamma$.

Recall that we denote the union of core curves of $P$ by $C^1$.
Let $h: C
\rightarrow C'$ be a homeomorphism inducing $\phi$ between
$\pi_1(C)\cong G$ and $\pi_1(C') \cong \Gamma$, which we know to exist because
both $C$ and $C'$ are handlebodies.
Recall that $p$ can be extended to a map from $C$ to $C'$ which is homotopic to
$h$ in $C'$ as was shown before Claim \ref{continuity}. 
By Lemma \ref{extension}, we see that  $h(\mu \cup C^1)$ is contained in
$\mathcal{D}(C')$.

Consider a convex cocompact hyperbolic metric  on
$\mathrm{Int} C'$, and denote the convex cocompact 3-manifold by $N$.
Since $h(\mu \cup C^1)$ lies in $\mathcal{D}(C')$, it is realised by
a pleated surface $g: \partial C' \rightarrow N$ homotopic to (a perturbation into $\Int C'$ of) the inclusion such that $g_*$ embeds the lift
of $h(\mu)$ into the unit tangent bundle  $T_1(N)$.
(See Theorem 5.5  Lecuire \cite{Lecuire}.)
Since $p$ extends to a map homotopic to $h$, we see that $p_*(\mu)$, which has turned out to be a measured lamination, is
homotopic to $h(\mu)$ in $C'$.
Therefore the support of  $h(\mu)$ is homotopic to that of $\lambda$ in $C'$, and if
we forget the transverse measures, 
the image of $g\circ
h(\mu)$ coincides with the realisation of $\lambda$ by a pleated surface $g':
\partial C' \rightarrow N$ homotopic to the inclusion.

We see that  leaves of $h(\mu)$ isolated from one side are identified with those
of $\lambda$ since both $g$ and $g'$ induce embeddings into the unit tangent bundle $T_1(N)$ of $h(\mu)$
and $\lambda$.
Therefore the complementary regions of $\mu$ correspond one-to-one to 
those of $\lambda$ with the number of sides preserved.
Recall that every complementary region of $h(\mu)$ except for the ones containing  
components of $h(C^1)$
is simply connected.
Therefore  we see that $\lambda$ cannot a complementary region which has  negative Euler characteristic by simply calculating the area. 
\end{proof}
By the fact which we showed just before Lemma \ref{Euler}, we are lead to a contradiction  whether $p|\Sigma$ has degree $0$ or not, if $(C',P'')$ is a relative compression body with  non-empty interior boundary or $(C', P'') \cong (\Sigma', \partial
\Sigma') \times I$.
Therefore, we see that the only possibility is that $\partial C'
\setminus P'=\Sigma'$, which is connected.
We have thus completed the proof of Lemma \ref{connected}
\end{proof}

%Since $c_k$ is not homotopic into $P'$ for large $k$, we see that $(p\circ \iota)^{-1}(c_k)$
%is the boundary of essential annulus for large $k$.
%Note  that $\gamma_k \cup P$ lies in $\mathcal{D}(\partial C)$ by Lemma
%\ref{extension} measured lamination in $\mathcal{D}(\partial C)$.
%The support of a limit of essential annuli cannot be
%contained in a component of a measured lamination in $\mathcal{D}(\partial C)$
%as was shown in Lemma 3.5 in \cite{Lecuire}.
%This is a contradiction.

The following  derives from the covering theorem of Canary \cite{CaT}.

\begin{lemma}
\label{strong}
In the present case (2) of the alternatives in p.\pageref{two cases}, the convergence of $\{(G_i,\phi_i)\}$ to $(\Gamma, \phi)$ is strong.
\end{lemma}

\begin{proof}
%Let $(C', P')$ be a relative core of $(\hyperbolic^3/\Gamma)_0$ as above which
%is a handlebody.
%Hence there is a homeomorphism $h: C \rightarrow C'$ homotopic to $\Phi|C$.
%Then $h(P)$ is homotopic to $P'$
By taking a subsequence, we can assume that $\{\phi_i(G)\}$ converges
geometrically 
to a Kleinian group $G_\infty$ containing $\Gamma$.
Let $q : M_\Gamma \rightarrow M_{G_\infty}$ be the
covering associated to the inclusion $\Gamma \subset G_\infty$.
Since  $(M_\Gamma)_0$ has only one end by Lemma \ref{connected}, which is
topologically tame, $q$ must be finite-sheeted by Canary's covering theorem \cite{CaT}.
This implies that $\Gamma = G_\infty$ by the argument in \S9.3  in Thurston  
\cite{ThL} (see Lemma 2.3 in Ohshika
\cite{OhQ} for a detailed proof).
\end{proof}

\subsection{Proof of Proposition \ref{ending lamination} for the case when $\Sigma$ is null-homologous}

Now we are ready to start the proof of Proposition \ref{ending lamination} in the case (2) of p.\pageref{two cases}.

\begin{proof}[Proof of Proposition \ref{ending lamination} in the case (2)]
Consider a nice compact core $(C',P')$ of $(M_\Gamma)_0$ which is
adequate with respect to a branched covering $\tilde{M} \rightarrow
M_\Gamma$ by a negatively curved manifold obtained by Canary's construction.
Then  $\partial C' \setminus P'$ is liftable, and by Lemma \ref{connected},  is connected.
We denote the complement of $C'$, which has a product structure, by $\partial
C' \times (0,\infty)$ identifying $\partial C' \times \{0\}$ with $\partial C'$.

Set $M_i = M_{\phi_i(G)}$.
Let $x$ be a point in $\hyperbolic^3$, and set its images in $M_i$ and
$M_\Gamma$
under the universal covering projections to be  basepoints $x_i$ and $x_\infty$.
Since $\{(G_i,\phi_i)\}$ converges to $(\Gamma,\phi)$ strongly by Lemma \ref{strong}, there is a  $(K_i, R_i)$-approximate
isometry
$\rho_i : B_{R_i}(M_i, x_i) \rightarrow B_{K_i R_i}(M_\Gamma,
x_\infty)$ with $R_i
\rightarrow \infty$ and $K_i \rightarrow 1$, which can be assumed to be a
diffeomorphism to its image.
Furthermore, $\rho_i^{-1}(C')$ is a compact core of
$M_i$ and $\Phi_i$ can be homotoped so that $\Phi_i|C$ is a homeomorphism onto $\rho_i^{-1}(C')$ for
every $i$ if we extract a subsequence.
Also, $\rho_i \circ \Phi_i|C$ is homotopic to $\Phi|C$ in
$M_\Gamma$.

 Take a generator system $g_1, \dots g_s$ of $G$, and for each $i$, let $b_i$ be a
bouquet in $M_i$ consisting of the geodesic loops representing $\phi_i(g_1), \dots ,
\phi_i(g_s)$ based at $x_i$.
We note that the boundary of the convex core of $M_i$ is incompressible outside $b_i$.
The sequence of bouquets $\{b_i\}$ converges geometrically to a bouquet $b_\infty$ in
$M_\Gamma$ based at $x_\infty$ representing $\phi(g_1) \dots ,
\phi(g_s)$.
We can choose an adequate compact core $(C',P')$ which contains the bouquet $b_\infty$.
Then $\partial C'$ is incompressible outside $b_\infty$.
Also, passing to a subsequence and enlarging $C'$ in its regular neighbourhood, we can assume that $\rho_i^{-1}(C')$ contains $b_i$. 
%By enlarging $C'$ if necessary, we can assume the distance between $\partial
%C'$ and $b_\infty$ is greater than the double of the maximum of the diameters
%of pleated surfaces homotopic to $\partial C'$ realising laminations in
%$\mathcal{D}(C')$ which intersect $b_\infty$.
%(Refer to Lecuire \cite{Lecuire} for the fact that there is an upper bound for
%such diameters.)

\begin{lemma}
\label{adequate}
By enlarging $C'$ further within $(M_\Gamma)_0$ preserving its adequacy, we can make $\rho_i^{-1}(C')$ also an adequate
compact core for any large $i$ (with respect to the covering $p_i :\tilde M_i \rightarrow M_i$).
\end{lemma}

\begin{proof}
Since $C'$ is adequate, we can assume that $\rho_i^{-1}(C')$ contains a branching locus of the covering $p_i$ for every $i$ by taking a subsequence.
What remains to show is that the complement of $\rho_i^{-1}(C')$ has a product structure.

 Fix some simple closed curve $d$ contained in $\mathcal{D}(C)$.
 We assume that $d$ is not contained in $C^1$.
Since the bending lamination, which we denote by $\beta_i$, of the convex core
$C(M_i)$
of $M_i$
is contained in $\mathcal{D}(C)$ and $\mathcal{D}(C)$ is arcwise connected
(Proposition 4.2 in Lecuire \cite{Lecuire}),
we can connect $d$ and $\beta_i$ by an arc in $\mathcal{D}(C)$.
This gives rise to a continuous family of pleated surfaces and negatively 
curved pleated surfaces in $M_i$ realising a homotopy $H_i$ between a pleated
surface realising $d$ and the boundary of the convex core $C(M_i)$.
(This fact can be shown by the argument in \S4.E of \cite{OhM}, simply replacing
the Masur domain by $\mathcal{D}$ since it was shown that every lamination in
$\mathcal{D}$ is realised in $M_i$ by Theorem 5.1 in \cite{Lecuire}.)
Since pleated surfaces realising $d$ in the $M_i$ converge geometrically to one
in $M_\Gamma$ realising $d$, there is a uniform bound for both
their diameters and the
distances from
them to the $\rho_i^{-1}(C')$.
By enlarging $C'$, we can assume that $\rho_i^{-1}(C')$ contains the
pleated surface realising $d$ for every $i$.

For any sufficiently large $i$, the homotopy $H_i$ contains a sub-homotopy $H_i'$ between
a pleated
surface touching the boundary of  $\rho_i^{-1}(C')$ and
the boundary
of the convex core of $M_i$, whose image is disjoint from the bouquet $b_i$.
We can see as follows that there is an upper bound independent of $i$ for the diameters modulo the $\epsilon$-Margulis tubes in $M_i$ of the pleated surfaces (and negatively curved pleated surfaces) constituting
$H_i'$.
Suppose, seeking a contradiction, that such an upper bound does not exist.
Then, there exist a (negatively curved) pleated surface $f_i$ constituting $H'_i$ and  a sequence of
simple closed curves $\delta_i$ on $\partial C$ such that
$\mathrm{length}(f_i(\delta_i))  \rightarrow 0$ and $f_i(\delta_i)$ is not contained in an $\epsilon$-Margulis tube.
Then $f_i(\delta_i)$ must be 
null-homotopic.
Since $f_i$ is homotopic to $\partial C(M_i)$  outside $b_i$, and any compressing disc of
$C(M_i)$ intersects $b_i$, we see that $f_i(\delta_i)$ bounds a singular disc
with area going to $0$
which intersects $b_i$ essentially.
By taking an annular neighbourhood of $\delta_i$ on $\partial C$ with height going to $\infty$, consisting of
circles of small lengths all of which bound singular discs with area going to $0$,
we see that the length
of some arc $b_i$ has length going to $\infty$ as $i \rightarrow \infty$.
(Refer to the proof of Lemma 6.2 in \cite{OhM} for a similar argument.)
This is a contradiction.

In a similar way, we can show that the distance between $b_i$ and $\partial C(M_i)$ goes to $\infty$ as follows.
We refer the reader to the proof of Theorem 2.1 in  \cite{OhQ} for a detailed account of the argument which we shall give briefly below.
Suppose, seeking a contradiction, that the distance is bounded after taking a subsequence.
There is a positive lower bound for the lengths of compressing curves on $\partial C(M_i)$ since any compressing disc must intersect $b_i$ essentially, and we can argue as in the previous paragraph.
Since the total length of $b_i$ is bounded as $i \rightarrow \infty$, if the distance between $b_i$ and $\partial C(M_i)$ is bounded, so is the distance from $x_i$ to $\partial C(M_i)$.
Therefore, choosing a basepoint within a bounded distance from $x_i$, the surface $\partial C(M_i)$ converges geometrically to a pleated surface in $M_\Gamma$ as $i \rightarrow \infty$.
Although the limit surface may not be embedded, it can be approximated by an embedded surface arbitrarily closely because the $\partial C(M_i)$ are embedded.
The approximating surfaces are separating in $(M_\Gamma)_0$ since $\partial C(M_i)$ is separating.
Recall that by Lemma \ref{connected}, $(M_\Gamma)_0$ has only one end, which is geometrically infinite.
Therefore, there are closed geodesics in a non-compact component of  the complement of an embedded surface approximating the limit surface, arbitrarily far from the limit surface.
Take such a closed geodesic $c_\infty$ sufficiently far from the limit surface.
We can pull back $c_\infty$ to $M_i$ by $\rho_i^{-1}$ for large $i$ as a quasi-geodesic converging to $c_\infty$ geometrically.
It follows that  the closed geodesic  homotopic to $\rho_i^{-1}(c_\infty)$ is also far from $\partial C(M_i)$, and is not contained in $C(M_i)$.
This contradicts the fact that every closed geodesic lies in $C(M_i)$.

%Take some positive constant $K$, which we shall specify later.
%Let $f: \partial C \rightarrow M_i$ be a pleated surface constituting $H_i'$
%which touches the $K$-neighbourhood of $\rho_i^{-1}(C')$ such that $\rho_i
%\circ f$ is contained in $\hyperbolic^3/\Gamma \setminus C'$.
%Since $C'$ is adequate, $\hyperbolic^3/\Gamma \setminus C'$ has a product
%structure $\partial C' \times (0,\infty)$.
%If $\rho_i \circ f$ is compressible in $\partial C' \times (0,\infty)$, it can
%be compressed in an arbitrary $\partial C' \times [a,b]$ containing the image
%of $\rho_i \circ f$;
%hence $f$ must be also compressible in $M_i \setminus b_i$.
%This is a contradiction since $f$ is homotopic outside $b_i$ to the boundary of
%$C(M_i)$.
%Therefore, $\rho_i \circ f$ is incompressible in $\partial C' \times
%(0,\infty)$, hence is homotopic
%to $
%\partial C'$ outside $C'$.

Since $\partial C(M_i)$ is incompressible in $M_i \setminus b_i$, so is
every surface constituting $H_i'$.
By the standard technique using Freedman-Hass-Scott \cite{FHS}, we see that
$H_i'$ gives rise to a product structure
$\partial C \times [0,1]$ starting from a surface within a uniformly bounded
distance from $\rho_i^{-1}(C')$ and ending at $\partial C(M_i)$ such that $\partial C \times \{t\}$ is incompressible in $M_i \setminus b_i$ and has diameter modulo the $\epsilon$-Margulis tubes uniformly bounded from above.
To be more precise, we can homotope $H_i'$ in the complement of $b_i$ to an embedding $H''_i: \partial C \times [0,1] \rightarrow C(M_i)$ such that $H''_i(\partial C \times \{0\})$ lies within the $R$-neighbourhood of $\rho_i^{-1}(C')$, $H''_i(\partial C \times \{1\})=\partial C(M_i)$, and the diameter of $H''_i(\partial C \times \{t\})$ modulo the $\epsilon$-Margulis tubes is bounded by a constant independent of $t$ and $i$.
We can fix a positive number $K$ so that for every large $i$ there exists  $F_i=H_i''(\ ,t_i)$ such
that $F_i(\partial C)$ intersects
the $K$-neighbourhood of $\rho_i^{-1}(C')$ whereas $H_i'(\ ,[t_i,1])$ is disjoint
from the $K/2$-neighbourhood of $\rho_i^{-1}(C')$.
(To show such $K$ exists, we use the facts that the diameters modulo the $\epsilon$-Margulis
tubes of surfaces constituting $H_i'$ are uniformly bounded, which was proved
above and that the distance from $x_i$ to $\partial C_i$ goes to $\infty$.)
Then for large $i$, its push-forward $\rho_i\circ F_i(\partial C)$ is contained in
$\partial C' \times  (0,\infty)$.

If $\rho_i \circ F_i$ is compressible in $\partial C' \times (0,\infty)$, it can
be compressed in an arbitrary $\partial C' \times [a,b]$ containing the image
of $\rho_i \circ F_i$;
hence $F_i$ must be also compressible in $M_i \setminus b_i$, which is a contradiction.
Therefore, $\rho_i \circ F_i$ is incompressible in $\partial C'
\times (0,\infty)$.
Since the diameters of the $\rho_i \circ F_i$ modulo the $\epsilon$-thin part are bounded above and $F_i(\partial C)$ intersect the $K$-neighbourhood of $\rho_i^{-1}(C')$, there exists $L>0$ such that $\partial C' \times (0, L)$ contains all the
$\rho_i \circ F_i( \partial C)$.
Since $F_i$ is homotopic to $H_i''(\ , t)$ in the complement of $\rho_i^{-1}(C')$, the product structure induced from $\rho_i^{-1}(\partial C'
\times (0,L))$ and that of $
\partial C \times (0,1)$ are isotopic outside $\rho^{-1}(C')$.
This gives rise to a
product structure between $\partial \rho_i^{-1}(C')$ and the boundary of $C(M_i)$.

On the other hand, since $M_i$ is convex cocompact (note that $M_i$ does not have
parabolic elements in this case (2)),
$M_i
\setminus C(M_i)$ has a product structure homeomorphic to $
\partial C \times \reals$, where $\partial C(M_i) \times \{-\infty\}$ is identified with $\partial C(M_i)$.
By pasting these product structures, we get a product structure of the complement of $\rho_i^{-1}(C')$.
Therefore $\rho_i^{-1}(C')$ is an adequate compact core for large $i$.
\end{proof}

{\em We isotope $\Phi$ so that $\Phi(C)=C'$, where $C'$ is enlarged as above.}
%%%%%

Now, let $g_i : (\partial C , m_i)\rightarrow M_i$ be a pleated surface homotopic to
$\Phi_i|\partial C$
realising $\mu \cup C^1$.

\begin{lemma}
\label{diameters bounded}
The diameters of the pleated surfaces $g_i$ modulo the $\epsilon$-Margulis tubes are bounded above by a constant depending only on $\epsilon$.
Furthermore the distance between $\rho_i^{-1}(C')$ and $g_i(\partial C)$ goes to
infinity as $i \rightarrow \infty$.
\end{lemma}

\begin{proof}
Suppose, seeking a contradiction, that the diameters of the $g_i$ modulo the $\epsilon$-Margulis tubes are not uniformly bounded.
Then, by the same argument as before, there is a simple closed curve $d_i$ on $\partial C \setminus C^1$ whose length with respect to $m_i$ goes to $0$ such that $g_i(d_i)$ is null-homotopic.
Since $\mu$ is contained in the Masur domain of $\Sigma=\partial C \setminus C^1$, it must intersect a measured lamination $\lambda$ to whose projective class  $[d_i]$ converges in $\PML(\Sigma)$.
Since the length of $d_i$ with respect to $m_i$ goes to $0$, the same argument as in Lemma \ref{unrealised} implies that the length of $\mu$ with respect to $m_i$ goes to infinity.
On the other hand,  the length of $\mu$ with respect to $\nu_i$, which was defined in Definition \ref{nu_i}, is bounded.
Recall that every meridian on $\partial C$ intersects $\mu\cup C^1$; hence
its length with respect to $\nu_i$ goes to $\infty$.
Then by the result of Canary \cite{Ca},
the length of $g_i(\mu)$ in $M_i$ is bounded as $i \rightarrow \infty$.
This is a contradiction, and we have shown the first statement.

%Next note that $g_i(C^1)$ is contained in a Margulis tube for sufficiently large $i$ since the length of $\Phi_i(C^1)$ goes to $0$ as $i \rightarrow \infty$.
%Since $M_i$ converges geometrically to $\hyperbolic^3/\Gamma$, and there are only finitely many Margulis tubes within bounded distance from the cusp neighbourhood corresponding to $\Phi(C^1)$, 
We shall next show the second statement of our lemma.
As was shown in the proof of  Lemma \ref{unrealised}, the image $g_i(\mu)$ cannot go into $\epsilon_i$-Margulis tubes with $\epsilon_i \rightarrow 0$.
If the image of $g_i$ stays within bounded distance from
$\rho_i^{-1}(C')$, it converges geometrically to a pleated surface realising $\mu$ in $M_\Gamma$ since $\{\phi_i(G)\}$ converges to $\Gamma$ strongly.
This contradicts Lemma \ref{unrealised}.
Therefore the distance between $\rho_i^{-1}(C')$ and $g_i(\partial C)$ goes to infinity as $i \rightarrow \infty$.
\end{proof}

This lemma implies in particular that for
sufficiently large $i$, the image of $g_i$ is contained in
the complement $U_i$ of $\rho_i^{-1}(C')$, which has a product structure
homeomorphic
 to $
\partial C
\times (0,\infty)$ by Lemma \ref{adequate}.

\begin{lemma}
\label{incompressible}
The pleated surface $g_i$ is incompressible outside $\rho_i^{-1}(C')$.
\end{lemma}

\begin{proof}
Suppose that  $g_i$ is compressible in the complement of $\rho_i^{-1}(C')$, which has a product structure by Lemma \ref{adequate} as $U_i \cong \partial C
\times (0,\infty)
\subset M_i$.
Then we can extend $g_i$
to a map from a compression body whose exterior boundary is identified with $
\partial C$, as can be seen by homotoping the map  into $\partial C \times \pt$ and using the simple loop conjecture proved by Gabai \cite{Ga}.
Since $\partial C' \times (0,\infty)$ does not contain an immersed
incompressible surface with genus less than that of $\partial C'\cong \partial C$, this
compression body must be a handlebody, hence can be identified with $C$.
(Since $g_i$ is homotopic to $\Phi_i|C$, any compression of $g_i$ can be done
within $C$.)
Since the compression can be performed in any $\partial C' \times [a,b]$
containing the image of $g_i$, there is an extension $\hat{g}_i: C \rightarrow
M_i$
of $g_i$ whose image has distance from $\rho_i^{-1}(C')$ going to $\infty$.
 Let $\sigma$ be a spine of $C$, \ie a one-complex which is a deformation retract
of $C$.
Note that there is a retraction $r: C \rightarrow \sigma$ such that for each non-vertex point $x$ in $\sigma$, its preimage $r^{-1}(x)$ is a meridian of $C$.
%Let $\pi: M \setminus C' \cong \partial C' \times (0,\infty) \rightarrow \partial
%C'$ be the projection to the first factor as before.
By Lemma \ref{adequate}, we can  consider the projection $\pi_i$ of
$\bar U_i \cong\partial C' \times [0,\infty)$ to $\partial \rho_i^{-1}(C')$.
We also consider the projection of $\mu \cup C^1$ under $r$.
The transverse measure of $\mu \cup C^1$ defines a weight for each point of $\sigma$.
Since each non-vertex point of $\sigma$ corresponds to a meridian on $\partial C$ and $\mu
\cup C^1$ is contained in $\mathcal{D}(C)$,
the image $r(\mu \cup C^1)$ passes through every point of $\sigma$ with weight bounded below by a positive constant $\eta$.
Therefore, if $\pi_i \circ \hat{g}_i(\sigma)$ has an essential self-intersection, then
$i((\pi_i
\circ g_i)_*(\mu \cup C^1), (\pi_i \circ
g_i)_*(\mu \cup C^1)) \geq \eta^2$.
See Figure \ref{fig:retraction}.

\begin{figure}
\label{fig:retraction}
\includegraphics[height=7cm]{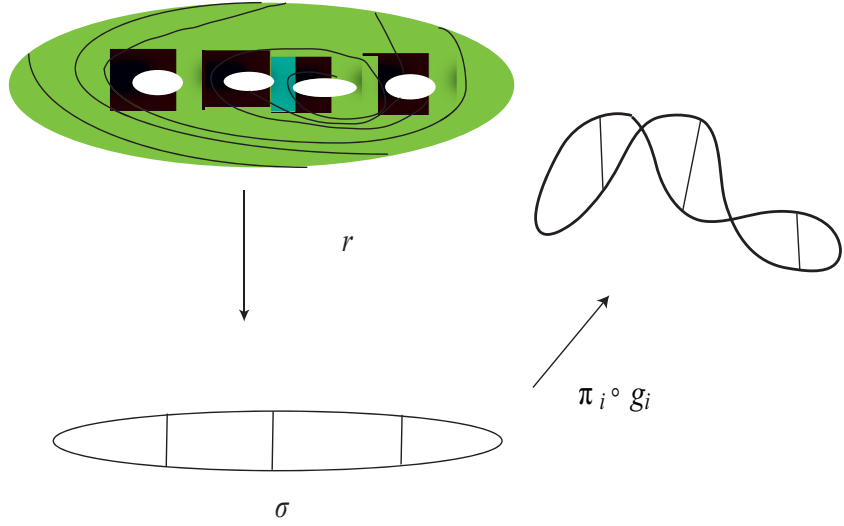}
\caption{The image of a lamination in $\mathcal{D}(C)$ passes through every point of $\sigma$.
Each self-intersection of $\pi_i \circ g_i(\sigma)$ contributes $\eta^2$ to $i((\pi_i \circ g_i)_*(\mu \cup C^1), (\pi_i \circ g_i)_* (\mu \cup C^1))$.}
\label{fig:spine}
\end{figure}

We shall prove that $\pi_i \circ \hat{g}_i(\sigma)$ has no self-intersection, hence is embedded.
For that, we shall show that {\em the homotopy class of $\rho_i \circ \pi_i \circ \hat g_i$ is independent of $i$ for large $i$}.
%For $i, j \in \naturals$, let $\Phi_{i,j}: M_i \rightarrow M_j$ be a homeomorphism such that $\Phi_{i,j} \circ \rho_i^{-1}=\rho_j^{-1}$ in the $R_i$-neighbourhood of $C'$ with $R_i \rightarrow \infty$.
%(It is obvious that such a homeomorphism exists after passing to a subsequence since $\rho_i^{-1}$ is  a homeomorphism from the $K_ir_i$-ball centred at $x_\infty$ with $K_ir_i \rightarrow \infty$.)
%Recall that $\phi_i(G)$ is a quasi-conformal deformation of $G$, corresponding to the marked conformal structure $\nu_i$ on $\partial C$.
We extend the conformal structures $\nu_i$ defined in Definition \ref{nu_i} to a continuous family $\nu_t$ for $t \in [1,\infty)$ such that $\nu_t$ escapes from any compact set of $\mathcal{T}(\partial C)$ and the length of $\mu \cup C^1$ with respect to $\nu_t$ is bounded as $t \rightarrow \infty$.
For each $t \in [i, i+1]$ with $i \in \naturals$, let $(G_t, \phi_t) \in AH(G)$ be the convex cocompact Kleinian group corresponding to $q(\nu_t)$, and set $M_t$ to be $M_{G_t}$ with homeomorphism $\Phi_t : M_G \rightarrow M_t$ induced from $\phi_t$.
Now, we shall show the existence of a homeomorphism $\bar \Phi_t: M_G \rightarrow M_t$  homotopic to $\Phi_t$ which has good properties as follows.

\begin{claim}
\label{homeomorphism}
For each $t \in [1,\infty)$, there is a homeomorphism $\bar\Phi_t: M_G \rightarrow M_t$ homotopic to $\Phi_t$ with the following properties.
\begin{enumerate}
%\item The submanifold $\bar\Phi_t(C)$ is an adequate compact core of $M_t$.
\item The diameters of $\bar \Phi_t(C)$ are bounded by a constant independent of $t$.
\item When $t$ is an integer $i$, we have $\bar \Phi_i|C=\rho_i^{-1}\circ \Phi|C$.
\item As $t \rightarrow t_0$ in $[1,\infty)$, the homeomorphism $\bar \Phi_t$ converges to $\bar\Phi_{t_0}$ geometrically.
\end{enumerate}
\end{claim}
\begin{proof}
We can take $C$ to be an adequate compact core.
We fix a homotopy class of simple closed curve $b$ in $\mathcal M(\Sigma)$.
Then it is obvious that $C^1 \cup b$ is  \lq \lq disc-busting" in Canary's sense, \ie intersects any compressing disc of $C$ essentially.  
When we consider a branched covering $p_t : \tilde M_t \rightarrow M_t$, its branching locus is always assumed to be the closed geodesic (or its perturbation if the closed geodesic is not simple) homotopic to $\Phi_t(b)$. 
We  take branched coverings which vary continuously with $t$ in the geometric topology.
Let $D_t$ be the infimum of the diameters of  all compact cores in $M_t$ that are adequate with respect to $p_t$.
We shall first show that $D_t$ is bounded from above by a constant independent of $t$.

Suppose, seeking a contradiction, that $D_t$ is not bounded.
Let $\{t_j\}$ be any monotone increasing sequence in $[1,\infty)$ such that $D_{t_j} \rightarrow \infty$.
Then by applying our result for $\{\phi_i\}$ in the previous section, regarding $\phi_{t_j}$ as $\phi_i$, we see that there is a subsequence $\{t'_j\}$ of $\{t_j\}$ such that $(G_{t'_j},\phi_{t'_j})$ converges  to $(\Gamma', \phi')$, which is either a convex cocompact group or a geometrically infinite group.
If $\Gamma'$ is convex cocompact, then $\{t'_j\}$ is bounded in $[1, \infty)$, $\Gamma'=G_{\lim t_j}$, and the convergence is strong.
If $\Gamma'$ is geometrically infinite, then $t'_j$ goes to $\infty$, and we can apply what we proved for $\Gamma$ up to now in this section.
In particular, Lemma \ref{strong} implies that the convergence to $\Gamma'$ is strong.

Set $M'=\hyperbolic ^3/\Gamma'$.
Let $\bar\rho_{t'_j} : B_{r_i}(M_{t'_j}, x_{t'_j}) \rightarrow B_{K_i r_i}(M', x'_\infty)$ be an approximate isometry associated to the geometric convergence of $M_{t'_j}$ to $M'$.
(We put a bar on $\rho$ to distinguish it from the approximate isometries associated to the convergence of $M_i$ to $M_\Gamma$.)
Take an adequate compact core $(\bar C, \bar P)$ of $(M')_0$.
Suppose first that  $M'$ is convex cocompact.
Then $M_{t_j}$ is quasi-isometric to $M'$ and $\rho_{t'_j}$ can be defined on the entire $M_{t_j}$ as a bi-Lipschitz homeomorphism.
We consider the branched covering of $M'$ branched along the closed geodesic homotopic to $\Phi'(b)$ or its perturbation as usual.
Then $\bar\rho_{t'_j}^{-1}(\bar C)$ is also an adequate compact core for large $j$ since  $\bar\rho_{t'_j}$ is a homeomorphism defined on the entire $M_{t_j}$ and branching loci move continuously in the geometric topology.
This shows that $(M_{t_j})_0$ has an adequate core with diameter less than $D_{t_j}$ for large $j$, contradicting the definition of $D_{t_j}$.

We shall next show that $\Phi'(b)$ does not represent a parabolic class, whether $\Gamma'$ is convex cocompact or geometrically infinite.
We can assume that $\Phi'|C$ is a homeomorphism to $\bar C$ since both of them are handlebodies.
Suppose, seeking a contradiction, that $\Phi'(b)$ represents a parabolic class.
Then, we can assume that $\Phi'(b) \subset \bar P$.
Since $C^1 \cup b$ is disc-busting, the complement $\bar F$ of $\Phi'(b \cup C^1)(\subset \bar P)$ in $\partial \bar C$ is incompressible in $\bar C$.
Also by applying Lemma \ref{connected} to $\Gamma'$, we see that there is only one end of $(M')_0$.
Since this end faces an essential subsurface of $\bar F$ or $\bar F$ itself, which is the incompressible frontier of $\bar C$ in $(M')_0$, by Bonahon's theorem, every unrealisable lamination is homotopic to the ending lamination.
Therefore $\Phi'(\mu)$ must be homotopic to a measured lamination in $\partial C \setminus (C^1 \cup b)$, which represents the ending lamination.
By the same argument as in the proof of Lemma \ref{Euler}, by comparing the complementary regions of $\mu$ and a measured lamination in $\partial C \setminus (C^1 \cup b)$, we see that this is impossible.
Thus we have shown that $\Phi'(b)$ cannot represent a parabolic class, and also that  we can consider the branching covering of $M'$ branched along the closed geodesic homotopic to $\Phi'(b)$ or its perturbation even when $\Gamma'$ is geometrically infinite.

Now, in the case when $\Gamma'$ is geometrically infinite, Lemma \ref{adequate} tells us that we can enlarge $\bar C$ so that $\bar\rho_{t'_j}^{-1}(\bar C)$ is an adequate compact core for large $j$.
This shows that $M_{t'_j}$ contains an adequate compact core whose diameter is less than $D_{t'_j}$ for large $j$, which is a contradiction.

Thus we have proved that $D_t$ is uniformly bounded from above, whether $\Gamma'$ is geometrically infinite or not.
We have also shown above that there is a positive lower bound $\xi$ for the lengths of the closed geodesics homotopic to the $\Phi_t(b)$. 
%We should note that this also implies that there is a positive lower bound $\xi$ for the lengths of the closed geodesic in $M_t$ homotopic to $\Phi_t(b)$.
Let $D$ be a constant which bounds  from above both the $D_t+1$ for all $t \in [1,\infty)$ and the diameters of $C'$ and the $\rho_i^{-1}(C')$.
%For each $t$, we let $\hat C_t$ be an adequate compact core with diameter less than $D+1$, which is guaranteed to exist by our definition of $D$.

Next, we shall show that there is a constant $E$ independent of $t$ such that $M_t$ has an adequate compact core $\hat C_t$ with diameter less than $2D$ whose  $4D$-neighbourhood  can be isotoped into $\hat C_t$ passing only the $E$-neighbourhood of $\hat C_t$.
Suppose, seeking a contradiction, that such $E$ does not exist.
Then, there is monotone increasing sequence $\{t_j\}$ in $[1,\infty)$ such that any adequate compact core $\hat C_{t_j}$ of $M_{t_j}$ with diameter less than $2D$ has $4D$-neighbourhood which cannot be isotoped into $\hat C_{t_j}$ within the $E_j$-neighbourhood of $\hat C_{t_j}$, where $E_j$ goes to $\infty$ as $j \rightarrow \infty$.
As before, after passing to a subsequence, $\{(G_{t_j}, \phi_{t_j})\}$ converges strongly to $(\Gamma', \phi')$ which is either convex cocompact or geometrically infinite, and  $(M')_0$ has an adequate compact core $\bar C$.
By the same argument as above, we can enlarge $\bar C$ so that for sufficiently large $j$, the submanifold $\bar \rho^{-1}_{t_j}(\bar C)$ is an adequate compact core of $M_{t_j}$.
We can choose such an enlarged adequate compact core to have diameter bounded by $D$ by our definition of $D$.
Since $\bar C$ is an adequate compact core in $(M')_0$, there is a constant $E'$ such that the $5D$-neighbourhood of $\bar C$ can be isotoped into $\bar C$ passing only the $E'$-neighbourhood of $\bar C$.
Pulling this back to $M_{t_j}$ for large $j$, we see that the $4D$-neighbourhood of $\bar\rho_{t_j}^{-1}(\bar C)$ can be isotoped into $\rho_{t_j}^{-1}(\bar C)$ passing only the $2E'$-neighbourhood of $\bar\rho_{t_j}^{-1}(\bar C)$.
Since the diameter of $\bar C$ is bounded by $D$, we see that the diameter of $\bar \rho_{t_j}^{-1}(\bar C)$ is less than $2D$ for large $j$.
Since $E_j$ is greater than $2E'$ for large $j$, this is a contradiction.
Thus we have proved the existence of $E$ as above.

For $D$ as was given above, we take  $K_0 \in (1,2]$ so that for any $K\leq K_0$ and any $K$-bi-Lipschitz homeomorphism $f: N \rightarrow N'$ between hyperbolic 3-manifolds, if a closed geodesic $c$ in $N$ has length greater than $\xi$, which was defined to be a lower bound for the lengths of the closed geodesics homotopic to the $\Phi_t(b)$,  then $f(c)$ is contained in the $D$-neighbourhood of the closed geodesic homotopic to $f(c)$.
The existence of such $K_0$ follows easily from the well-known properties of quasi-geodesics.

Now, for each $i \in \naturals$, we subdivide $[i, i+1]$ into intervals $i=s_0< s_1 <  \dots<s_{m-1} < s_m=i+1$ in such a way that for any $t, t' \in [s_j, s_{j+1}]$,  the quasi-conformal deformation from $(G_t, \phi_t)$ to $(G_{t'}, \phi_{t'})$ induces (by extending quasi-conformal maps to quasi-isometries in $\hyperbolic^3$ in the standard way) a $(K_0)^{1/2}$-bi-Lipschitz homeomorphism $\Psi_{t,t'}$ from $M_t$ to $M_{t'}$ which is homotopic to $\Phi_{t'} \circ \Phi_t^{-1}$.
(The number $m$ may depend on $i$.)
Since $\{(G_t, \phi_t)\}$ is a family of quasi-conformal deformations which is continuous with respect to the parameter $t$,  any subdivision whose maximal width is sufficiently small serves as a subdivision as above.
Now for each $s_j\, (j\neq 0,m)$, we define  $\bar C_{s_j}$ to be an adequate compact core of $M_{s_j}$whose diameter is less than $D_{s_j}+1$, which is guaranteed to exist by our definition of $D_t$.
For $j=0,m$, we define $\bar C_{s_j}$ to be $\rho_i^{-1}(C')$ and $\rho_{i+1}^{-1}(C')$ respectively.
We define $\bar \Phi_{s_j}|C$ to be a homeomorphism from $C$ to $\bar C_{s_j}$ homotopic in $M_{s_j}$ to $\Phi_{s_j}|C$.
In the case when $j=0, m$ we choose the homeomorphism to be $\rho_i^{-1} \circ \Phi$.
Since both $C$ and $\bar C_{s_j}$ are adequate, their complements have product structures and we can extend $\bar \Phi_{s_j}|C$ to a homeomorphism $\bar \Phi_{s_j}: M_G \rightarrow M_{s_j}$.
Then we see that the diameter of $\bar \Phi_t(C)$ is bounded by $D$ for $t=s_0, \dots , s_m$, and (ii) in our claim holds for $\bar \Phi_t$ defined thus far (\ie when $t$ is a subdividing point).

We need to extend the family of homeomorphisms to all the parameters $t \in [1,\infty)$.
Now we further subdivide $[s_j, s_{j+1}]$ into three subintervals of the same length, $[s_j, s_j'], [s_j', s_j'']$ and $[s_j'', s_{j+1}]$.
We consider a $(K_0)^{1/2}$-bi-Lipschitz homeomorphism $\Psi_{s_j,t}$ for $t \in [s_{j-1}', s_j']$.
Since  this bi-Lipschitz homeomorphism $\Psi_{s_j,t}$ is induced from a quasi-conformal deformation from $(G_{s_j}, \phi_{s_j})$ to $(G_t, \phi_t)$, we can assume that it varies continuously with respect to $t$ in the geometric topology, and $\Psi_{s_j,t}$ converges to the identity geometrically as $t \rightarrow s_j$.
%Furthermore, we can assume that this family has the property that $\Psi_{s, t+t'}=\Psi_{t,t'} \circ \Psi_{s,t}$.
For each $t \in [s_j, s_j']$, we define $\bar\Phi_t$ to be $\Psi_{s_j,t}\circ \bar \Phi_{s_j}$, and for each $t\in [s_j'', s_{j+1}]$, we define $\bar\Phi_t$ to be $\Psi_{s_{j+1},t} \circ \bar \Phi_{s_{j+1}}$.
Now, we shall fill the gap between $\bar \Phi_{s_j'}$ and $\bar \Phi_{s_j''}$.

Since both $\bar C_{s_j}$ and $\bar C_{s_{j+1}}$ have diameters bounded by $D$, the compact cores      $\bar \Phi_{s_j'}(C)$ and $ \Psi_{s_{j+1},s_j'}\circ \bar \Phi_{s_{j+1}}(C)$ have diameters bounded by $2D$ (because $K_0^{1/2} \leq \sqrt 2 <2$).
Moreover, since $\bar C_{s_j}$ contains the closed geodesic $b^*_{s_j}$ homotopic to $\Phi_{s_j}(b)$, and $\Psi_{s_j, s_j'}$ is ($K_0^{1/2}$-, hence also) $K_0$-bi-Lipschitz, $\bar \Phi_{s_j'}(C)$, which contains $\Psi_{s_j, s_j'}(b^*_{s_j})$, is within distance $D$ from any point on the closed geodesic homotopic to $\Phi_{s_j'}(b)$.
By considering the same for $\Psi_{s_{j+1}, s_j'} \circ \bar \Phi_{s_{j+1}}(C)=\Psi_{s_{j+1}, s_j'}(\bar C_{s_{j+1}})$, we see that  $\Psi_{s_{j+1}, s_j'} \circ \bar \Phi_{s_{j+1}}(C)$ is within distance $2D$ from $\bar \Phi_{s_j'}(C)$, hence is contained in the $4D$-neighbourhood of $\bar \Phi_{s_j'}(C)$.
By our definition of $E$, there is an isotopy taking $\Psi_{s_{j+1}, s_j'}\circ \bar \Phi_{s_{j+1}}(C)$ into $\bar \Phi_{s_j'}(C)$ which passes only the $E$-neighbourhood of $\bar \Phi_{s_j'}(C)$.
Since the complement of $C$ has a product structure, we can extend this isotopy to an isotopy $H: M_G \times [s_j', s_j''] \rightarrow M_{s_j'}$ such that each $H(\ ,t)$ is a homeomorphism to $M_{s_j'}$, $H(\ ,s_j')=\bar \Phi_{s_j'}$,  $H(\ ,s_j'')=\Psi_{s_{j+1},s_j'}\circ \bar \Phi_{s_{j+1}}$, and $H(C,t)$ is contained in the $E$-neighbourhood of $\bar \Phi_{s_j'}(C)$.
Now, for $t \in [s_j', s_j'']$, we define $\bar \Phi_t$ to be $\Psi_{s_{j+1},t} \circ \Psi_{s_{j+1}, s_j'}^{-1} \circ H(\ , t)$.
Since $\Psi_{s_{j+1},t} \circ \Psi_{s_{j+1}, s_j'}^{-1}$ is $K_0$-bi-Lipschitz with $K_0 \leq 2$, we see that $\bar \Phi_t(C)$ has diameter bounded by $2(2D+E)$.
Thus we have completed the definition of $\bar\Phi_t$, which have the properties (i) and (ii).
The continuity with respect to the geometric topology is obvious from our definition of $\bar \Phi_t$.
\end{proof}

Now we return to the proof of Lemma \ref{incompressible}.
Let $g_t: \partial C \rightarrow M_t$ be a pleated surface homotopic to $\Phi_t$ which realises $\mu \cup C^1$.
Since the length of $\mu$ with respect to $\nu_t$ is bounded as $t \rightarrow \infty$ and $C^1 \cup \mu$ intersects every compressing curve, by the same argument as for $g_i$, we see the length of $g_t(\mu)$ is bounded as $t \rightarrow \infty$.
Since the diameter of $\bar \Phi_t(C)$ is bounded as $t \rightarrow \infty$ by Claim \ref{homeomorphism}, by the same argument as for $g_i$, it follows that the surface $g_t(\partial C)$ is disjoint from $\bar \Phi_t(C)$ for large $t$.
Note that both $g_t(\mu)$ and $\bar \Phi_t$ vary continuously with respect to the geometric topology.
Let $U_t \cong \partial C \times (0,\infty)$ be a parametrisation of the complement of $\bar \Phi_t(C)$.
Let $\pi_t: M_t \rightarrow \partial \bar \Phi_t(C)$ be the projection to the first factor in $\bar U_t \cong \partial C \times [0,\infty)$, where we identify $\partial C \times \{0\}$ with $\bar \Phi_t(\partial C)$.
Since the parametrisation and the projection are unique up to homotopy, we can assume that $U_t$ and $\pi_t$ vary continuously with respect to the geometric topology.
Since $g_t(\partial C)$ never touches $\bar \Phi_t(C)$, this implies that the homotopy class of $\bar \Phi_t^{-1} \circ \pi_t \circ g_t$ is independent of  $t$.
Recall that $\rho_i^{-1}|C'=\bar \Phi_i \circ \Phi^{-1}|C'$ for every $i \in \naturals$.
Therefore for $i \in \naturals$, we see that the homotopy class of $\rho_i \circ \pi_i \circ g_i$ as a map from $\partial C$ to $\partial C'$ is independent of $i$.
 %hence so is the image of $\hat{g}_t$.
%%Probablement on a besoin de d?tailler. Pourquois les diam?tres des surface%%
%%sont born?es.
%This implies that for large and near $t, t'$, the surfaces $\Phi_{t,t'} \circ
%\hat{g}_t$ and $\hat{g}_{t'}$ are homotopic outside $\rho_{t'}^{-1}(C')$.
%It follows that $\Phi_{t,t'} \circ \pi_t \circ \hat{g}_t$ is homotopic to
%$\pi_{t'}
%\circ \hat{g}_{t'}$.
%Therefore the homotopy classes of $\rho_t \circ \pi_t \circ \hat{g}_t$ on
%$
%\partial
% C'$ are constant for large $t$, hence for all $i$ up to extracting a
%subsequence.
Let $p :  \partial C \rightarrow \partial C'$ be a map homotopic to these maps.
If $\pi_i \circ g_i(\sigma)$ has essential self-intersection, then $i(p_*(\mu \cup C^1), p_*(\mu \cup C^1)) \geq \eta^2$ for a positive constant $\eta$ independent of $i$ as was shown before.

Take a sequence of weighted simple closed curves $\{r_k c_k\}$ on $\Sigma=\partial C
\setminus P$ converging to $\mu$.
By taking a subsequence of $\{r_k c_k\}$, we can assume that for each $i$, there is a pleated
surface $g'_i : \partial C \rightarrow M_i$ realising $C^1 \cup c_i$ which is
homotopic to $g_i$ in $U_i$, whose diameter modulo the $\epsilon$-Margulis tubes is
bounded as $i \rightarrow \infty$.
We can extend $g_i'$ to $\hat{g}_i' : C \rightarrow M_i$ by the same way as
for $g_i$.
%Then the double points of $\hat{g}_i'(\sigma)$ correspond to the same homotopy
%classes of meridians for $C$ as those of $\hat{g}_i(\Sigma)$.
%This implies that $i(\pi_i \circ g_i'(c_i \cup C^1), \pi_i \circ g_i'(c_i \cup
%C^1)) \geq \eta^2/2$, hence $i(p(c_i \cup C^1), p(c_i \cup C^1)) \geq
%\eta^2/2$.
By Lemma \ref{uniform intersection}, there is a sequence of weighted simple
closed curves $R_k C_k$ on $\partial C \setminus P$ converging to $\mu'$ with
$i(\mu,
\mu')=0$ such that $i(R_i \pi_i \circ g'_i(C_i), R_i \pi_i \circ g'_i(C_i))
\rightarrow 0$.
Since $g_i'$ is homotopic to $g_i$ in $U_i$, we get $i(R_ip(C_i), R_ip(C_i))
\rightarrow 0$, which implies $i(p_*(\mu'), p_*(\mu'))=0$, hence also $i(p_*(\mu),
p_*(\mu))=0$.
Since $\pi_i \circ g'_i(C_i)$ is disjoint from $\Phi_i(C^1)$ as was shown in
Lemma
\ref{uniform intersection}, we have $i(p_*(\mu'), p(C^1))=0$.
These show that $i(p_*(\mu \cup C^1), p_*(\mu \cup C^1))=0$.
%This contradicts Lemma \ref{uniform intersection} since the distance from
%$\rho_i^{-1}(C')$ and the image of $g_i$ goes to $\infty$ as $i \rightarrow
%\infty$.
%(Note that Lemma \ref{uniform intersection} can be applied to measured
%laminations by approximating them by weighted simple closed curves.)
Thus $\sigma$ can be embedded on $\partial C'$ by $p$ if we move $p$ by a homotopy on $\partial C'$.

%Since the length of $\phi_t(\mu)$ is bounded as $t \rightarrow \infty$, for
%large $t$, the realisation $g_t(\mu)$ is disjoint from $\rho_t^{-1}(C')$,
%hence so is the image of $\hat{g}_t$.
%%Probablement on a besoin de d?tailler. Pourquois les diam?tres des surface%%
%%sont born?es.
%This implies that for large and near $t, t'$, the surfaces $\Phi_{t,t'} \circ
%\hat{g}_t$ and $\hat{g}_{t'}$ are homotopic outside $\rho_{t'}^{-1}(C')$.
%It follows that $\Phi_{t,t'} \circ \pi_t \circ \hat{g}_t$ is homotopic to
%$\pi_{t'}
%\circ \hat{g}_{t'}$.
%Therefore the homotopy classes of $\rho_t \circ \pi_t \circ \hat{g}_t$ on
%$
%\partial
% C'$ are constant for large $t$, hence for all $i$ up to extracting a
%subsequence.
%Let $p :  C \rightarrow \partial C'$ be a map homotopic to these maps.
We have also shown that $p_*(\mu)$, regarded as a geodesic current, is actually a measured lamination on $\partial C' \setminus P'$ since it has null-self-intersection.
On the other hand since $p|\sigma$ is embedding, $p_*(\mu)$ has a large complementary region, which is contained in $\partial C' \setminus (P' \cup p(\sigma))$.
Then by the same argument as in the proof of Lemma \ref{Euler}
considering the correspondence between the complementary regions of $p_*(\mu)$ and a measured lamination representing the 
ending lamination for the unique end of $(M_\Gamma)_0$,
we reach a contradiction.
Thus, we have shown that $g_i$ is incompressible in $U_i\cong \partial C \times (0,\infty)$.
\end{proof}

Since $g_i$ is incompressible in $U_i \cong \partial C \times (0,\infty)$, we
see that $\pi_i \circ g_i$ is homotopic a homeomorphism from $\partial C$ to
$\partial \bar C_i=\rho_i^{-1}(\partial C')$.
As was shown above, the maps $\rho_i \circ \pi_i \circ g_i$ are all homotopic,
which we denote by $h: \partial C \rightarrow \partial C'$.
Note that $h$ is homotopic in $C'$ to $\Phi|\partial C$, hence can be extended to a
homeomorphism from $C$ to $C'$.
Also since $g_i$ maps $C^1$ into Margulis tubes converging geometrically to cusps in $M_\Gamma$, we
can assume that $h(P)=P'$.
Take a weighted simple closed curve $\{s_kd_k\}$ on $\partial C' \setminus 
P'$ converging to
a measured lamination $\lambda$ representing the ending lamination of the unique end of $(M_\Gamma)_0$
such that the closed geodesic $d_k^*$ homotopic to $d_k$ in $U \cong (\partial
C' \setminus P')
\times [0,\infty)$ tends to the end.
Recall that we have a sequence of weighted simple closed curves $\{r_k c_k\}$
converging to $\mu$ with a pleated surface $g_i'$ realising $c_i \cup C^1$ in $M_i$.
Then by Lemma \ref{uniform intersection}, we have sequences of  simple
closed curves $C_k$ on $\partial C \setminus P$ and $D_k$ on $\partial C'
\setminus P'$ such that $[C_k]$ converges to $[\mu']$ with $i(\mu, \mu')=0$
and $[D_k]$ converges to $[\lambda']$ with $i(\lambda, \lambda')=0$, and we
have $\displaystyle i(\frac{\pi_i \circ g_i'(C_i)}{\mathrm{length}(\pi_i \circ g_i'(C_i))},
\frac{\rho_i^{-1}(D_i)}{\mathrm{length}(\rho_i^{-1}(D_i))}) \rightarrow 0$.
By pushing this forward by $\rho_i$, we get $\displaystyle i(\frac{h(C_i)}{\mathrm{length}(h(C_i))},
\frac{D_i}{\mathrm{length}(D_i)}) \rightarrow 0$, hence $i(h(\mu'), \lambda')=0$.
Since both $\mu$ and $\lambda$ are arational, $\mu'$ has the same support as $\mu$,
and $\lambda'$ has the same one as $\lambda$.
It follows that $h(\mu)$ has the same support as $\lambda$.
Since $h$ extends to a homeomorphism from $C$ to $C'$ taking $P$ to
$P'$ as was shown above, we have completed the proof of Proposition \ref{ending lamination} in the case (2).
%By pulling back $c_k^*$ by $\rho_i^{-1}$ and letting $k \rightarrow \infty$,
%which forces $i$ to go to $\infty$ in order to make $c_k^*$ lie in
%$B_{K_ir_i}(\hyperbolic^3/\Gamma, x_\infty)$,
%we see that $i(\Phi_i^{-1}\rho_i^{-1}(\lambda),
%\mu \cup C^1) \rightarrow 0$ by Lemma \ref{uniform intersection}.
%Since we homotoped the $\Phi_i$ so that  $\Phi_{i,i'} \circ \rho^{-1}_{i}|C'=
%\rho_{i'}^{-1}|C'$ for large $i,i'$, there is a homeomorphism $k: C
%\rightarrow C'$ which coincides with  $\rho_i \circ \Phi_i|C$ for large $i$.
%Then we have $i(\lambda, k(\mu \cup C^1))=0$.
%Let $C^\infty$ be the union of core curves of $P'$.
%Since  both $\mu \cup C^1$ and $\lambda \cup C^\infty$ are maximal, this is
%possible only when $k(\mu)$ has the same support as $\lambda$.
%This also implies that $k(P)=P'$ 
%and $k(\mu)$ represents an ending lamination for the end facing $\partial C'
%\setminus P'$.
Combining Lemma \ref{homologically non-trivial} with this, we have
completed the proof of Proposition \ref{ending lamination}.
\end{proof}

\section{Proof of the main theorem}
We shall complete the proof of Theorem \ref{main} in this section.
By Proposition \ref{ending lamination}, we have shown that $\Phi$ can be homotoped
so that $\Phi|P$ is an embedding into $P'$ and $\Phi|\Sigma_j$ is a
homeomorphism to a component of $\partial C' \setminus P'$ for $j =n+1, \dots
, m$ such that $\Phi|\mu_j$ represents the ending lamination for the end
facing the component, which is topologically tame.
It remains to deal with $\Phi|\Sigma_j$ for $j =1, \dots , n$.

Let $\Sigma$ be one of the $\Sigma_j$ for $j =1, \dots , n$.
Let $S_k$ be the boundary component of $C$ containing $\Sigma$, and $H$ a
subgroup of $G$ corresponding to the image of $\pi_1(S_k)$.
Let $H^\Sigma$ denote a subgroup of $H$ corresponding to 
the image of $\pi_1(\Sigma)$ in $\pi_1(C) \cong G$.
%Recall that the marked conformal structures at infinity $n_j^k$ on $\Sigma$
%were
%taken to
% converge to $m_j$.
%Let $\Omega_k$ be a component of the region of discontinuity of 
%$G$ covering the surface at infinity corresponding to $S_k$.
%Then the quasi-conformal deformation $(G_i, \phi_i)$ of $G$ is 
%realised by a quasi-conformal homeomorphism $h_i: S^2_\infty 
%\rightarrow S^2_\infty$ such that there is an open subset $\Omega$ 
%of $\Omega^k$
%on which $h_i|\Omega$ converges to an equivariant homeomorphism to an open
%set $\Omega'$ whose stabiliser is $\phi(H^\Sigma)$.
%(See Abikoff \cite{Ab}.)
%Since every frontier component of $\Sigma$ is mapped by $\Phi$ to a 
%closed curve representing a parabolic class, the surface 
%$\Omega'/\phi(H^\Sigma)$ is of finite type.
Then by the same argument as the proof of Lemma \ref{compression body} using the result of Abikoff, we see that there are an open subset $\Omega'$ of $\Omega_\Gamma$ invariant under $\phi(H^\Sigma)$, and  a frontier component of the convex core 
of $(M_\Gamma)_0$ facing $\Omega'/\phi(H^\Sigma)$.
Since $C'$ is a nice compact core, there is a component of $\partial C' \setminus P'$ ambient isotopic to this frontier component.
Therefore  the surface $\Phi|\Sigma$ is homotoped to this component 
of $\partial C' \setminus P'$ keeping $\Phi(\Fr \Sigma)$ in $P'$.
Repeating this argument for every one of the $\Sigma_j \ (j=1, \dots , n)$, and combining it  with Proposition \ref{ending lamination}, we see that $\Phi$ is
homotopic to a homeomorphism from $C$ to $C'$, hence from $M_G$ to
$M_\Gamma$, by Waldhausen's theorem \cite{Wa}.
Moreover, the ends of $(M_\Gamma)_0$ facing these 
components $\Sigma_1, \dots, \Sigma_n$ are geometrically finite and have conformal structure at infinity $m_1, \dots , m_n$.
This completes the proof of Theorem \ref{main}.

\bigskip  % gives some vertical space
\name{\noindent Ken'ichi Ohshika}

\address{
\noindent 
Department of Mathematics,
Graduate School of Science,
Osaka University,
Toyonaka, Osaka 560-0043,
Japan%
}

\email{\noindent ohshika@math.sci.osaka-u.ac.jp}

\classmark{57M50, 30F40}

\keywords{\noindent Kleinian group, Deformation space, End invariant}

\end{document}